\documentclass[11pt]{article}
\usepackage[margin=1in]{geometry}
\usepackage{amsmath,amssymb,amsthm,mathtools,booktabs,enumitem,multirow}
\usepackage[authoryear,round]{natbib}
\usepackage{xcolor,graphicx}
\usepackage[colorlinks=true,linkcolor=blue!50!black,citecolor=blue!50!black,urlcolor=blue!50!black]{hyperref}
\usepackage{setspace,placeins}
\newtheorem{theorem}{Theorem}
\newtheorem{proposition}{Proposition}
\newtheorem{corollary}{Corollary}
\newtheorem{lemma}{Lemma}[section]
\theoremstyle{definition}
\newtheorem{assumption}{Assumption}

\theoremstyle{remark}
\newtheorem{remark}{Remark}

\newcommand{\R}{\mathbb{R}}
\newcommand{\E}{\mathbb{E}}
\newcommand{\Pb}{\mathbb{P}}
\newcommand{\T}{\top}
\newcommand{\tr}{\operatorname{tr}}
\newcommand{\cF}{\mathcal{F}}
\newcommand{\toP}{\xrightarrow{\ \Pb\ }}
\newcommand{\toD}{\xrightarrow{\ d\ }}
\newcommand{\eps}{\varepsilon}
\newcommand{\Hspan}{\mathcal{H}^{+}}
\newcommand{\jb}[1]{\langle #1\rangle}

\title{Exact and Approximate Maximum Likelihood Estimation for Multivariate Diffusions under High-Frequency Sampling}
\author{Duy Nguyen\thanks{Marist University, Poughkeepsie, NY 12601, USA. Email: \texttt{nducduy@gmail.com}.}}

\begin{document}
\maketitle

\begin{abstract}
We establish the asymptotic theory of the exact maximum likelihood estimator (MLE) for discretely observed multivariate ergodic diffusions under high-frequency, long-span sampling, solving an open problem posed by  \citet{YCW2019}  and \citet{Li2013} . The drift may grow linearly, and the volatility is bounded away from zero and either bounded or proportional to $1+|x|$. When the sampling interval decays polynomially, the exact MLE is consistent, asymptotically normal and efficient, with rates $\sqrt{n\Delta_n}$ for drift and $\sqrt n$ for diffusion parameters. Under an additional condition on ergodic averages, and when the time span $n\Delta_n$ also grows polynomially, its moments converge as well. The proof never analyzes the exact likelihood directly: it shows that a closed-form density expansion of sufficiently high order is uniformly accurate in relative terms at every observed transition, so that the two log-likelihoods differ by $o_P(1)$ uniformly in the parameter. The same argument yields accuracy conditions under which an approximate likelihood is asymptotically equivalent to the exact one; in particular, the approximate MLE of \citet{YCW2019} converges to the true parameter. For the delta-expansion estimator of fixed order $K\ge1$, we derive the leading bias explicitly. It shows that $n\Delta_n^{K+1}\to0$ cannot be weakened when the bias constant is nonzero, and that the bias falls on the drift parameters for even $K$ and on the diffusion parameters for odd $K$. Monte Carlo experiments confirm the predicted biases.


\end{abstract}

\section{Introduction}

Multivariate diffusions are used throughout financial econometrics, in models of the term structure, of stochastic volatility, and of commodity and credit markets. Their parameters are usually estimated by maximum likelihood from discretely sampled data. The difficulty is that the transition density, and hence the likelihood, is rarely known in closed form.

Following \citet{AitSahalia2002,AitSahalia2008}, a large literature replaces the unknown density with an explicit small-time expansion and maximizes the resulting approximate likelihood. The delta expansion of \citet{YCW2019} is a recent member of this family. It is derived from an It\^o--Taylor argument, and its coefficients follow from a recursion that involves only derivatives of the drift and diffusion coefficients.

The statistical justification of these estimators has a gap. \citet{YCW2019} show that their approximate estimator approaches the exact MLE as the sampling interval $\Delta$ tends to zero with the number of observations $n$ fixed (their Theorem~4.1). To conclude that it converges to the true parameter, one also needs the asymptotic theory of the exact MLE. On this point, \citet[page 264]{YCW2019} write that ``the corresponding results about the convergence of the true MLEs in multivariate diffusions are challenging and still open in the literature.'' \citet{Li2013} leaves the same question for future work. For scalar diffusions, the gap is largely filled by \citet{AitSahalia2002} and \citet{ChangChen2011}. Both work with the transformation to unit volatility, which in general exists only in one dimension. \citet{ChangChen2011} also obtain asymptotic normality of the exact MLE only when $n\Delta_n^3\to\infty$, so that the sampling interval shrinks more slowly than $n^{-1/3}$, and under high-level conditions on the expansion terms. The open question therefore concerns the multivariate case in the high-frequency regime, with the sampling interval allowed to shrink at any polynomial rate.

\paragraph{Setting.}
We work in the high-frequency, long-span regime: $n\to\infty$, $\Delta_n\to0$ and $n\Delta_n\to\infty$. This is the natural regime for the question. With $\Delta$ fixed, an expansion of fixed order has an error that does not vanish, so the approximate estimator converges to a pseudo-true value. With the time span $n\Delta_n$ fixed, the drift parameters cannot be estimated consistently \citep{GenonCatalotJacod1993,Gobet2001}. We consider ergodic diffusions whose drift is Lipschitz and may grow linearly, and whose volatility is bounded away from zero and either bounded or proportional to the size of the state, comparable to $1+|x|$ in every direction. In the second case, the coefficients must also vary on the scale of the state, in the sense of the derivative bounds in Assumption~\ref{A:lin}. This class contains multivariate Ornstein--Uhlenbeck processes, Gaussian affine term-structure models, nonlinear models with bounded volatility, and models whose volatility grows with the level of the process, such as Student-type diffusions with heavy-tailed stationary laws. It does not contain square-root models such as CIR and Heston, whose volatility vanishes at the boundary.

\paragraph{The main idea.}
The exact likelihood is not available in closed form, so we never analyze it directly. Instead, we compare it with the likelihood built from a delta expansion of high order $K$. If the expansion is accurate in relative terms at every observed transition, simultaneously for all parameter values, the two log-likelihoods differ by $o_\Pb(1)$ uniformly, and the exact MLE inherits the asymptotic behavior of the explicit contrast. Comparing the exact likelihood with an explicit approximation is an old idea \citep{DacunhaCastelleFlorensZmirou1986,AitSahalia2002}. What makes it work here is that a high-order expansion is accurate enough to absorb even the most extreme observation in the sample.

\paragraph{Contributions.}
The paper makes four main contributions and provides one technical tool of independent interest.
\begin{enumerate}[leftmargin=*,itemsep=3pt]
\item \emph{Asymptotic theory for the exact MLE} (Theorem~\ref{thm:mle} and Corollary~\ref{cor:lan}). The exact MLE is consistent and asymptotically normal whenever $\Delta_n$ decays polynomially, with rate $\sqrt{n\Delta_n}$ for drift and $\sqrt n$ for diffusion parameters. The same argument gives a proof of local asymptotic normality (LAN) without Malliavin calculus; LAN was first established by \citet{Gobet2002}. Hence the exact MLE is efficient. Under an additional assumption on the rate of ergodic averages, which we verify only for stationary Ornstein--Uhlenbeck processes, and when $n\Delta_n$ also grows polynomially, its moments also converge (Theorem~\ref{thm:moment}).
\item \emph{A transfer principle for approximate likelihoods} (Theorem~\ref{thm:general}). We give a high-level condition under which an approximate likelihood yields an estimator asymptotically equivalent to the exact MLE, together with two sufficient conditions on the error of the density or of the log-density. The conditions concern only the accuracy of the approximation, not how it is constructed. They are worst-case bounds, so the theorem is a theoretical transfer principle rather than a practical recipe. Delta expansions of high order satisfy them with any centering. This includes the expansion of \citet{Li2013}, which \citet{ChenWanYang2026} show to be the delta expansion with $\mu_0=0$, and the conditions are also met, for scalar diffusions, the expansion of \citet{LeeSongLee2014}, under a restrictive condition on the Lamperti-transformed drift.
\item \emph{The bias of the delta-expansion estimator} (Theorems~\ref{thm:cons} and \ref{thm:an}). For the estimator of fixed order $K$ that is used in practice, efficiency under $n\Delta_n^{K+1}\to0$ is expected from the theory of higher-order contrasts \citep{Kessler1997,UchidaYoshida2012}. What is new is an explicit bias term of order $\sqrt n\Delta_n^{(K+1)/2}$, given by the projection of the first omitted correction on the score. It shows that the condition cannot be weakened whenever this bias is nonzero. For $K\ge1$, the leading bias falls on the drift parameters when $K$ is even and on the diffusion parameters when $K$ is odd; for the Euler contrast, $K=0$, it falls on the diffusion parameters. This alternation could not be seen in classical expansions, which truncate in integer powers of $\Delta$ and so correspond to even $K$ only; it appears because the delta expansion is graded in powers of $\Delta^{1/2}$.
\item \emph{Volatility proportional to the state} (Proposition~\ref{prop:lin} and Theorem~\ref{thm:lin}). All results except Corollary~\ref{cor:others}(b) extend to volatility that is comparable to $1+|x|$ in every direction, when the drift and volatility vary on the scale of the state; moment convergence then requires a version of the rate assumption. This class contains heavy-tailed models, such as Student-type diffusions, whose stationary laws have only finitely many moments. The extension rests on the observation that such models look the same at every scale: the expansion error, measured on the local scale, is bounded uniformly in the starting point, and no moments of the stationary law are needed. Intermediate growth, such as $(1+|x|^2)^{1/4}$, and growth in some directions only are not covered.
\item \emph{A technical tool: a uniform bound on the expansion error} (Proposition~\ref{prop:unif}). We give a self-contained proof that the delta expansion of every order is accurate on the whole state space, with a polynomial weight in the starting point. \citet{YCW2019} prove such a bound for starting points in a compact set, bounded drift and $K\ge m$.
\end{enumerate}
Together, these results supply the missing half of the decomposition argument of \citet{YCW2019} and \citet{Li2013}, justify Wald inference based on the exact likelihood, and tell the user of a fixed-order expansion how large $K$ must be and which parameters suffer if it is too small. A Monte Carlo study confirms the theory: once $\kappa_{22}\Delta_n$, the fastest mean-reversion rate times the sampling interval, is below about $0.1$, the bias predicted by Theorem~\ref{thm:an} matches the simulated bias to within a few hundredths of a standard deviation.

\paragraph{Related literature.}
The exact MLE has received less attention than its approximations, because it is not available in closed form. For reducible diffusions it can be computed by Monte Carlo \citep{BeskosEtAl2006,BeskosPapaspiliopoulosRoberts2009}, but its asymptotic theory under high-frequency sampling is a separate question. \citet{Gobet2002} establishes LAN for multivariate elliptic ergodic diffusions, with the same rates and information matrix as ours. LAN identifies the efficiency bound but is silent about the MLE. Showing that the MLE attains the bound additionally requires global consistency and a quadratic approximation of the log-likelihood that is uniform near $\theta_0$, and both require control of the exact likelihood over the whole parameter space. Theorem~\ref{thm:mle}(a) supplies exactly this control. How our proof of LAN compares with Gobet's is discussed after Corollary~\ref{cor:lan}.

\citet{Sorensen2024} develops a general efficiency theory for approximate martingale estimating functions of scalar diffusions sampled at high frequency; see also \citet{JacodSorensen2018} and \citet{KesslerLindnerSorensen2012}. Applying a score-based theory of this kind to the exact MLE requires control of the $\theta$-derivatives of $\log p_\Delta$. One could try to obtain such control by differentiating the Duhamel representation behind Proposition~\ref{prop:unif} in $\theta$. Our argument shows that this is unnecessary: it uses only the size of $p_\Delta$, and never its derivatives in $\theta$. With $\Delta$ fixed, the exact MLE is covered by the general theory for ergodic Markov chains \citep{Billingsley1961}, subject to verifying its regularity conditions, as \citet{AitSahalia2002} does for scalar diffusions. We are not aware of a complete proof of consistency and asymptotic normality of the exact MLE for multivariate ergodic diffusions in the high-frequency regime.

For approximate estimators, the theory of contrast and quasi-likelihood estimators in this regime is well developed. \citet{FlorensZmirou1989} and \citet{Yoshida1992} study Euler-type contrasts, and \citet{Kessler1997} adds higher-order corrections for scalar diffusions. \citet{UchidaYoshida2012,UchidaYoshida2014} and \citet{KamataniUchida2015} develop adaptive and multi-step estimators for multivariate models within the quasi-likelihood analysis of \citet{Yoshida2011}, and \citet{Jacobsen2001,Jacobsen2002} studies small-$\Delta$ optimality of estimating functions. See also \citet{PilipovicSamsonDitlevsen2024} for splitting schemes and \citet{Ganguly2024} for estimators without density expansions.

Closest to our analysis of the delta-expansion estimator, \citet{IguchiBeskosGraham2025} build a likelihood proxy from a small-time expansion with Hermite-type corrections, for elliptic and hypo-elliptic models, and prove asymptotic normality under $n\Delta_n^3\to0$, the same rate condition as for $K=2$ in our Theorem~\ref{thm:an}. Their proxy is a Gaussian-type contrast with corrected moments rather than a delta expansion. Our results for the delta-expansion estimator cover every order $K$, give the bias explicitly, and show that the rate condition is sharp. On the other hand, we treat only elliptic models. Explicit small-$\Delta$ expansions of the bias of discretized estimating equations, including the Euler scheme, go back to \citet{AitSahaliaMykland2003,AitSahaliaMykland2004}. In principle, our bias constant could be obtained from their general formulas applied to the score of the delta-expansion contrast. What is new is its closed form for the delta expansion of every order in the multivariate case, its parity structure, and the resulting sharpness of the rate condition.

The rest of this paper is organized as follows: Section~\ref{sec:setup} recalls the delta expansion and defines the estimator. Section~\ref{sec:ass} states the assumptions, and Section~\ref{sec:main} presents the results. Section~\ref{sec:mc} reports the Monte Carlo evidence, and Section~\ref{sec:disc} concludes. Proofs are collected in the appendices.

\section{The delta expansion and the estimator}\label{sec:setup}

\paragraph{Model and notation.}
We observe an $m$-dimensional diffusion
\begin{equation}\label{eq:sde}
dX_t=\mu(X_t;\alpha)\,dt+\sigma(X_t;\beta)\,dW_t,\qquad \nu(x;\beta)=\sigma(x;\beta)\sigma(x;\beta)^\T ,
\end{equation}
at the times $t_i=i\Delta_n$, $i=0,\dots,n$, where $W$ is an $m$-dimensional Brownian motion. The parameter is $\theta=(\alpha,\beta)$, with true value $\theta_0=(\alpha_0,\beta_0)$. We write $X_i=X_{t_i}$, $\Delta=\Delta_n$, $\eps=\sqrt{\Delta}$ and $\jb{x}=(1+|x|^2)^{1/2}$, and we let $\cF_i$ be the $\sigma$-field generated by $(X_s)_{s\le t_i}$. We keep the subscript in asymptotic statements, which involve $n$, and drop it when a statement concerns a single transition. In Propositions~\ref{prop:unif} and \ref{prop:lin}, which are not asymptotic, $\Delta\in(0,1]$ is arbitrary. The transition density over one step is $p_\Delta(x'|x;\theta)$. The exact log-likelihood, conditional on $X_0$, is $\ell_n(\theta)=\sum_{i=1}^n\log p_\Delta(X_i|X_{i-1};\theta)$. Under the assumptions of Section~\ref{sec:ass}, $p_\Delta(x'|x;\theta)$ is strictly positive and jointly continuous in $(x,x',\theta)$ (Lemma~\ref{lem:reg}). Hence $\ell_n$ is continuous on the compact set $\Theta$, and a measurable maximizer $\hat\theta_n$ exists \citep[Lemma~2]{Jennrich1969}. All our results for $\hat\theta_n$ except moment convergence (Theorem~\ref{thm:moment}), which uses an exact maximizer, also hold for approximate maximizers, with $\ell_n(\hat\theta_n)\ge\sup_\Theta\ell_n-o_\Pb(1)$. We also refer to the delta-expansion estimator defined below as the approximate MLE of \citet{YCW2019}.

\paragraph{The delta expansion.}
\citet{YCW2019} expand $p_\Delta$ around the density of a Brownian motion with drift whose coefficients are frozen at the starting point. In its general form, the frozen drift is a constant $\mu_0$, which may depend on the starting point $x$ and on $\theta$. The subscript in $\mu_0$ does not refer to the true parameter: $\mu_0$ is evaluated at the same $\theta$ as the rest of the expansion. Let
\[
z=\frac{x'-x-\mu_0\Delta}{\sqrt\Delta},\qquad
q(x'|x;\theta)=\Delta^{-m/2}\phi\big(z;\nu(x;\beta)\big),
\]
where $\phi(\cdot;\Sigma)$ is the $N(0,\Sigma)$ density. The delta expansion of order $K$ is
\begin{equation}\label{eq:delta}
\begin{aligned}
p^{(K,\Delta)}(x'|x;\theta)&=q(x'|x;\theta)\Big(1+\sum_{k=1}^K\Delta^{k/2}C_k(z;x,\theta)\Big),\\
C_k&=\sum_{N=\lceil(k+1)/2\rceil}^{2k}\frac1{N!}\sum_{|h|=2N-k}w_{N,h}(x;\theta)\,H_h\big(z;\nu(x;\beta)\big),
\end{aligned}
\end{equation}
where $H_h(z;\Sigma)=(-1)^{|h|}\phi(z;\Sigma)^{-1}\partial^h\phi(z;\Sigma)$ are multivariate Hermite polynomials. The coefficients $w_{N,h}(x;\theta)$ come from a recursion that involves only $\mu$, $\nu$, their derivatives and $\mu_0$ (equations (18)--(20) of \citet{YCW2019}, restated in Appendix~\ref{app:prop}).

For estimation we use the centering $\mu_0=\mu(x;\alpha)$, which \citet{YCW2019} recommend and use in their numerical work. With this choice, the Gaussian kernel $q$ is the Euler density, and the first correction $C_1$ does not involve the drift parameters (Lemma~\ref{lem:struct}(c)). Unless stated otherwise, $z$, $q$ and $C_k$ refer to this centering. Other centerings are treated in Section~\ref{sec:main}.

\paragraph{The estimator.}
Far in the tails, the bracket in \eqref{eq:delta} can become negative, and then its logarithm is undefined. With probability tending to one, no observation falls there (Lemma~\ref{lem:max}). Still, an estimator should be defined on every sample. We therefore write $R(z;x,\theta)=\sum_{k=1}^K\eps^kC_k(z;x,\theta)$, with $R\equiv0$ for $K=0$, and for $K\ge0$ maximize the truncated contrast
\begin{equation}\label{eq:contrast}
\ell_n^{(K)}(\theta)=\sum_{i=1}^n\Big[\log q(X_i|X_{i-1};\theta)+\log\max\big\{1+R_i(\theta),\tfrac12\big\}\Big],\qquad R_i(\theta)=R\big(z_i(\alpha);X_{i-1},\theta\big),
\end{equation}
where $z_i(\alpha)=(X_i-X_{i-1}-\mu(X_{i-1};\alpha)\Delta)/\eps$. The delta-expansion estimator $\hat\theta^{(K)}_n=(\hat\alpha_n,\hat\beta_n)$ is any measurable maximizer of $\ell^{(K)}_n$ over $\Theta$. On the event
\[
\Omega_n=\Big\{\max_{1\le i\le n}\ \sup_{\theta\in\Theta}|R_i(\theta)|\le\tfrac14\Big\},
\]
the truncation is inactive, and $\ell^{(K)}_n$ is the approximate log-likelihood of \citet{YCW2019}. The two constants play different roles: the floor $\frac12$ keeps the contrast finite on every sample, while on $\Omega_n$, where $|R_i|\le\frac14$, the logarithm and its derivatives are bounded uniformly, which is what the proofs use. For $K=0$ we set $R\equiv0$, so $\hat\theta^{(0)}_n$ is the Euler (Gaussian quasi-) maximum likelihood estimator.

Two properties of \eqref{eq:delta} drive the analysis. First, the leading factor $q$ is the Euler density, so the delta-expansion contrast is the Euler contrast plus a correction of order $\eps$ per observation. This reduces consistency to a perturbation of the Euler theory. Second, every $C_k$ is a combination of Hermite polynomials of positive degree, which have mean zero under the local Gaussian law $N(0,\nu(x;\beta))$. Hence $p^{(K,\Delta)}$ integrates to one exactly, and the error of the expansion shows up as an explicit bias of the score. This is where the order $K$ enters.

\section{Assumptions}\label{sec:ass}

The assumptions of this section cover bounded volatility; Section~\ref{sec:lin} replaces Assumptions~\ref{A:smooth} and \ref{A:ell} by a version for volatility of linear growth. The assumptions are of three kinds, plus one used only for moments. Smoothness and ellipticity (\ref{A:smooth}--\ref{A:ell}) control the expansion, ergodicity (\ref{A:erg}) supplies laws of large numbers, and identifiability (\ref{A:id}) pins down the true parameter. A rate condition on ergodic averages (\ref{A:rate}) is used only for moment convergence.

\begin{assumption}[parameters and smoothness]\label{A:smooth}
$\Theta=\Theta_\alpha\times\Theta_\beta$, where $\Theta_\alpha\subset\R^{p_\alpha}$ and $\Theta_\beta\subset\R^{p_\beta}$ are closures of bounded open convex sets, and $\theta_0$ lies in the interior of $\Theta$. The drift depends only on $\alpha$ and the diffusion coefficient only on $\beta$. The functions $\mu(x;\alpha)$ and $\sigma(x;\beta)$ are $C^\infty$ in $(x,\theta)$. The function $\sigma$ and all its partial derivatives are bounded on $\R^m\times\Theta$. Every partial derivative of $\mu$ of order at least one in $x$ is bounded on $\R^m\times\Theta$, and $|\partial^b_\alpha\mu(x;\alpha)|\le C\jb{x}$ for every multi-index $b$.
\end{assumption}

Separating drift and diffusion parameters is standard \citep{Kessler1997,Gobet2002,UchidaYoshida2012}, and all models in \citet{YCW2019} have this structure. The two groups are estimated at different rates, and the separation makes the limiting information block-diagonal. Many models in which a parameter appears in both coefficients can be reparametrized to separate them. For instance, the Ornstein--Uhlenbeck process $dX_t=-\kappa X_t\,dt+\sqrt{2\kappa v}\,dW_t$, parametrized by its stationary variance $v$, separates with $\sigma^2=2\kappa v$ as the diffusion parameter. When a parameter genuinely enters both coefficients, as in the Pearson diffusions of \citet{FormanSorensen2008}, it is identified through the diffusion coefficient at rate $\sqrt n$, and the rate matrix and the information structure change; we leave this case for future work.

The conditions on $\mu$ say that the drift is Lipschitz in $x$ and grows at most linearly. They allow the Ornstein--Uhlenbeck drift $\mu(x;\alpha)=\kappa(\bar x-x)$, Gaussian affine models, and nonlinear drifts such as $-\kappa x+g(x)$ with $g$ bounded and smooth. They exclude superlinear drifts such as $-x^3$. The conditions on $\sigma$ exclude square-root volatilities. We need them for the heat kernel bounds of \citet{MenozziPesceZhang2021}, on which Proposition~\ref{prop:unif} rests. Convexity of $\Theta$ is used in mean-value expansions and in a Sobolev inequality for suprema over the parameter space.

\begin{assumption}[uniform ellipticity]\label{A:ell}
There are constants $0<\lambda_{\min}\le\lambda_{\max}<\infty$ with $\lambda_{\min}I\le\nu(x;\beta)\le\lambda_{\max}I$ for all $x$ and $\beta$. We write $\Lambda=\lambda_{\max}/\lambda_{\min}$.
\end{assumption}

Uniform ellipticity is Assumption~2.3 of \citet{YCW2019}. It gives a smooth transition density with Gaussian bounds and keeps the Hermite polynomials in \eqref{eq:delta} well defined. The ratio $\Lambda$ reappears in Theorem~\ref{thm:mle}, where it controls how extreme an observed increment can look relative to the Gaussian kernel $q$.

\begin{assumption}[ergodicity and moments]\label{A:erg}
Under $\theta_0$: (i) $X$ is positive Harris recurrent with invariant probability $\pi$, so that $T^{-1}\int_0^Tg(X_s)\,ds\to\int g\,d\pi$ almost surely, for every initial law and every $\pi$-integrable $g$; (ii) the initial law has finite moments of all orders, and $\sup_{t\ge0}\E|X_t|^p<\infty$ for every $p\ge1$.
\end{assumption}

This assumption turns sums over observations into integrals against $\pi$. For Harris recurrence and the associated ergodic theorems, see \citet{Harris1956}, \citet{MeynTweedie1993II} and \citet{KaspiMandelbaum1994}.

A convenient sufficient condition is mean reversion: $\langle\mu(x;\alpha_0),x\rangle\le-c|x|$ for $|x|\ge r$, together with $\E e^{a|X_0|}<\infty$ for some $a>0$. To see why it suffices, let
\[
\mathcal Af(x)=\langle\mu(x;\alpha_0),\nabla f(x)\rangle+\tfrac12\tr\big(\nu(x;\beta_0)\nabla^2f(x)\big)
\]
be the generator of $X$ under $\theta_0$, acting on functions $f\in C^2(\R^m)$, and let $V(x)=e^{a'|x|}$ for $|x|\ge r$, smoothed near the origin. For $|x|\ge r$, $\nabla V=a'V\,x/|x|$ and $\nabla^2V=V\big[a'^2xx^\T/|x|^2+a'(I-xx^\T/|x|^2)/|x|\big]$, so
\[
\mathcal AV(x)=V(x)\Big[a'\frac{\langle\mu(x;\alpha_0),x\rangle}{|x|}+\frac{a'^2}2\,\frac{x^\T\nu x}{|x|^2}+\frac{a'}{2|x|}\Big(\tr\nu-\frac{x^\T\nu x}{|x|^2}\Big)\Big].
\]
The first term is at most $-a'c$ by mean reversion, the second at most $a'^2\lambda_{\max}/2$ by Assumption~\ref{A:ell}, and the third lies between $0$ and $a'm\lambda_{\max}/(2|x|)$. Hence
\[
\mathcal AV\le V\big[-a'c+a'^2\lambda_{\max}/2+O(|x|^{-1})\big].
\]
For $a'<c/\lambda_{\max}$ the bracket is at most $-a'c/2$ once $|x|$ is large, so $\mathcal AV\le-\frac{a'c}2V+b\mathbf 1_{B}$ for a constant $b$ and a ball $B$. This is a Foster--Lyapunov drift condition. Together with the irreducibility implied by uniform ellipticity, this gives positive Harris recurrence and uniform exponential moments \citep{MeynTweedie1993,DownMeynTweedie1995}. Linear drifts $-\kappa x$ whose matrix $\kappa$ has a positive definite symmetric part satisfy the condition, and so does the bounded drift in Section~\ref{sec:mc}. For a linear drift $-\kappa x$ with a general stable $\kappa$, that is, with all eigenvalues of $\kappa$ having positive real part, the Lyapunov equation $\kappa^\T P+P\kappa=I$ has a unique symmetric positive definite solution $P$, and $V(x)=\exp\{a'(x^\T Px)^{1/2}\}$, smoothed near the origin, plays the same role for small $a'>0$. More generally, it does so for any drift with $\langle\mu(x;\alpha_0),Px\rangle\le-c|x|$ for large $|x|$ and some symmetric positive definite $P$; for $P\neq I$ this condition is not implied by the Euclidean one above.

\begin{assumption}[identifiability]\label{A:id}
(a) If $\nu(x;\beta)=\nu(x;\beta_0)$ for $\pi$-almost every $x$, then $\beta=\beta_0$.
(b) If $\mu(x;\alpha)=\mu(x;\alpha_0)$ for $\pi$-almost every $x$, then $\alpha=\alpha_0$.
(c) The matrices
\[
\Gamma_\alpha=\int(\partial_\alpha\mu)^\T\nu^{-1}(\partial_\alpha\mu)\,d\pi,\qquad
(\Gamma_\beta)_{jk}=\frac12\int\tr\big(\nu^{-1}\partial_{\beta_j}\nu\,\nu^{-1}\partial_{\beta_k}\nu\big)\,d\pi,
\]
evaluated at $\theta_0$, are positive definite. We write $\Gamma=\mathrm{diag}(\Gamma_\alpha,\Gamma_\beta)$.
\end{assumption}

Parts (a) and (b) ensure that the limiting contrasts in the consistency proof have unique maximizers. Part (c) is the corresponding local condition, and $\Gamma$ is the Fisher information. Information about the drift accumulates with the time span $n\Delta_n$ rather than with the number of observations, which is why $\alpha$ is estimated at the slower rate.

\begin{assumption}[rate of ergodic averages]\label{A:rate}
For every $p\ge1$ and every $C^1$ function $g$ with $|g(x)|+|\nabla g(x)|\le C\jb{x}^r$ for some $C$ and $r$, $\sup_{T\ge1}T^{p/2}\,\E\big|T^{-1}\int_0^Tg(X_s)\,ds-\int g\,d\pi\big|^p\le C_{p,r}C^p$, where $C_{p,r}$ depends only on $p$ and $r$.
\end{assumption}

The expectation is taken under the initial law of the data. This assumption is standard in the quasi-likelihood analysis of ergodic diffusions \citep{Yoshida2011,UchidaYoshida2012}. One route to it runs as follows. Under the drift condition below Assumption~\ref{A:erg}, the process is exponentially ergodic in a weighted norm \citep{DownMeynTweedie1995}, so the stationary process is exponentially $\beta$-mixing. Moment inequalities for averages of mixing processes \citep[Section~1.4]{Doukhan1994} then give the bound in Assumption~\ref{A:rate} for functions with polynomial growth. For the stationary Ornstein--Uhlenbeck processes of Section~\ref{sec:mc}, all of these ingredients are standard. For other initial laws, the bound has to be transferred from the stationary case, which requires a coupling argument that we do not pursue.

\section{Main results}\label{sec:main}

This section is organized around the exact MLE. Section~\ref{sec:unif} provides the tool: a bound on the error of the delta expansion on the whole state space. Section~\ref{sec:exact} states the main result, the asymptotic theory of the exact MLE, and Section~\ref{sec:other} extends it to arbitrary density approximations that come with a quantitative error bound. Section~\ref{sec:approx} analyzes the delta-expansion estimator of fixed order; its consistency and quadratic expansion are also the main ingredients in the proof of the result on the exact MLE. Section~\ref{sec:moment} treats moments, and Section~\ref{sec:lin} extends all results to volatility of linear growth. Throughout, we assume that $\Delta_n\le Cn^{-a}$ for some $a>0$; Remark~\ref{rem:poly} explains why. Since $n\Delta_n\to\infty$, it is necessary that $a<1$.

\subsection{A uniform bound on the expansion error}\label{sec:unif}

Theorem~3.3 of \citet{YCW2019} bounds the error of the delta expansion for starting points in a compact set, bounded drift and orders $K\ge m$. That is not enough for our purpose. An ergodic process visits every region of the state space, and every transition contributes to the likelihood, so we need a bound on all of $\R^m$. With an unbounded drift the bound cannot be uniform in the starting point, but a polynomial weight suffices.

\begin{proposition}[uniform expansion error]\label{prop:unif}
Suppose Assumptions~\ref{A:smooth}--\ref{A:ell} hold, and let $\mu_0=\mu_0(x;\theta)$ be any centering with $|\mu_0(x;\theta)|\le C\jb{x}$. For every integer $K\ge1$ there are constants $C_K$ and $r_K$ such that, for all $\Delta\in(0,1]$, $x,x'\in\R^m$ and $\theta\in\Theta$,
\[
\big|p^{(K,\Delta)}(x'|x;\theta)-p_\Delta(x'|x;\theta)\big|\le C_K\jb{x}^{r_K}\Delta^{(K+1-m)/2}.
\]
\end{proposition}

Since $q$ is of order $\Delta^{-m/2}$ in the bulk of the distribution, the proposition says that the relative error of the expansion is of order $\Delta^{(K+1)/2}$ where the data lie.

The proof (Appendix~\ref{app:prop}) is self-contained apart from the heat kernel bounds of \citet{MenozziPesceZhang2021}. These are Gaussian bounds for the transition density and its gradient, centered at the deterministic flow of the drift. We write the error as a time integral by a Duhamel argument, show that the integrand vanishes near the starting point to an order that compensates for the singularity of the Hermite kernels, and integrate. The weight $\jb{x}^{r_K}$ comes from the linear growth of the drift.

\subsection{The exact MLE}\label{sec:exact}

We now turn to the question left open by \citet{YCW2019}. As explained in the introduction, the delta expansion serves here as a tool of proof rather than of computation. Drift and diffusion parameters are estimated at different rates, collected in the matrix
\[
D_n=\mathrm{diag}\big(\sqrt{n\Delta_n}\,I_{p_\alpha},\ \sqrt n\,I_{p_\beta}\big),
\]
where $I_p$ is the $p\times p$ identity matrix and $p_\alpha$, $p_\beta$ are the dimensions of $\alpha$ and $\beta$ (Assumption~\ref{A:smooth}). We write $S_n$ for the normalized martingale part of the Euler score at $\theta_0$, a sum of martingale differences with $S_n\toD N(0,\Gamma)$. Lemma~\ref{lem:clt} gives the precise definition and shows that the martingale parts of the scores of all delta-expansion contrasts differ from $S_n$ by $o_\Pb(1)$.

\begin{theorem}[exact MLE]\label{thm:mle}
Suppose Assumptions~\ref{A:smooth}--\ref{A:id} hold, $n\Delta_n\to\infty$, and $\Delta_n\le Cn^{-a}$ for some $a\in(0,1)$ and $C<\infty$. Let $\hat\theta_n$ be a measurable maximizer of $\ell_n$, or more generally a measurable $\hat\theta_n$ with $\ell_n(\hat\theta_n)\ge\sup_\Theta\ell_n-o_\Pb(1)$. Then:
\begin{enumerate}[label=(\alph*),nosep]
\item for every integer $K$ with $K+1>2(1+m\Lambda)/a$, $\sup_{\theta\in\Theta}|\ell_n(\theta)-\ell^{(K)}_n(\theta)|\toP0$;
\item $\hat\theta_n\toP\theta_0$ and $D_n(\hat\theta_n-\theta_0)=\Gamma^{-1}S_n+o_\Pb(1)\toD N\big(0,\mathrm{diag}(\Gamma_\alpha^{-1},\Gamma_\beta^{-1})\big)$;
\item for every $K\ge1$ with $n\Delta_n^{K+1}\to0$, $D_n(\hat\theta^{(K)}_n-\hat\theta_n)\toP0$.
\end{enumerate}
\end{theorem}

Part (a) is the heart of the argument. By Proposition~\ref{prop:unif}, the relative error of the expansion at an observation is of order $\Delta^{(K+1)/2}$ divided by the Gaussian kernel. The kernel is smallest at the most extreme standardized increment, which is of order $\sqrt{\log n}$, and this costs at most a factor $n^{(1+\delta)m\Lambda}$ for any $\delta>0$. A high enough order makes the total error over $n$ observations vanish. The threshold $2(1+m\Lambda)/a$ is crude, since it guards against the most extreme increment and the least favorable parameter value at once. This is harmless: the order $K$ in (a) is a device of the proof and is never computed in practice.

The proof has two inputs. The consistency and quadratic expansion of the explicit contrast are stated in Section~\ref{sec:approx} (Theorems~\ref{thm:cons} and \ref{thm:an}), which can be read independently of Sections~\ref{sec:exact}--\ref{sec:other}. Proposition~\ref{prop:unif} makes the explicit contrast $\ell^{(K)}_n$ a uniformly accurate surrogate for $\ell_n$, and the analysis of $\ell^{(K)}_n$ in Section~\ref{sec:approx} supplies consistency and a uniform quadratic expansion. Neither input involves derivatives of $p_\Delta$ in $\theta$.

Part (b) needs no rate condition on $\Delta_n$ beyond polynomial decay. Part (c) complements Theorem~4.1 of \citet{YCW2019}, which concerns the limit $\Delta\to0$ with $n$ fixed: the approximate and exact estimators agree at the scale of their own sampling error, for every order at which the approximate estimator is asymptotically unbiased. Together, (b) and (c) show that both terms of the decomposition
\[
\hat\theta^{(K)}_n-\theta_0=(\hat\theta^{(K)}_n-\hat\theta_n)+(\hat\theta_n-\theta_0)
\]
used by \citet{YCW2019} vanish at the optimal rates.

The uniform approximation in (a), combined with the quadratic expansion of the explicit contrast, also gives the LAN property.

\begin{corollary}[local asymptotic normality and efficiency]\label{cor:lan}
Under the assumptions of Theorem~\ref{thm:mle}, for every bounded sequence $u_n\in\R^{p_\alpha+p_\beta}$,
\[
\ell_n(\theta_0+D_n^{-1}u_n)-\ell_n(\theta_0)=u_n^\T S_n-\tfrac12u_n^\T\Gamma u_n+o_\Pb(1),\qquad S_n\toD N(0,\Gamma).
\]
Consequently, the exact MLE and, when $n\Delta_n^{K+1}\to0$, the delta-expansion estimator of order $K$ are asymptotically efficient in the sense of the H\'ajek--Le Cam convolution and local asymptotic minimax theorems.
\end{corollary}

The LAN property itself is not new. \citet{Gobet2002} proves it with the same rates and information matrix, using Malliavin calculus and requiring only $\Delta_n\to0$; his conditions and ours are not nested. Our proof avoids Malliavin calculus, relying instead on the heat kernel bounds behind Proposition~\ref{prop:unif}, at the price of polynomial decay of $\Delta_n$. Efficiency follows because both estimators are asymptotically linear in the central sequence $S_n$ \citep{Jeganathan1982,LeCamYang2000}. The rate matrix $D_n$ has different rates in the two blocks, but it is deterministic, so the standard theory for LAN experiments applies.

\paragraph{Consequences for inference.}
Theorem~\ref{thm:mle} and Corollary~\ref{cor:lan} justify Wald inference based on the exact likelihood. A consistent estimate of $\Gamma$ is obtained by plugging $\hat\theta_n$ into sample averages: $\hat\Gamma_\alpha=n^{-1}\sum_i(\partial_\alpha\mu)^\T\nu^{-1}(\partial_\alpha\mu)$ and $(\hat\Gamma_\beta)_{jk}=(2n)^{-1}\sum_i\tr(\nu^{-1}\partial_{\beta_j}\nu\,\nu^{-1}\partial_{\beta_k}\nu)$, evaluated at $(X_{i-1};\hat\theta_n)$. Its consistency follows from Lemma~\ref{lem:lln} and continuity. Wald confidence sets based on $\hat\theta_n$ and $\hat\Gamma$ therefore have asymptotically correct coverage. The uniform quadratic expansion \eqref{eq:laq} also suggests the usual chi-square limits for likelihood-ratio statistics. Because the two blocks are estimated at different rates, a general statement requires care with the restricted estimator, and we do not pursue it here. Theorem~\ref{thm:moment} below adds convergence of moments, which justifies mean squared error approximations.

\subsection{Other density approximations}\label{sec:other}

The proof of Theorem~\ref{thm:mle} uses the delta expansion only through its accuracy. The next theorem replaces it with an arbitrary approximation. Part (a) gives a high-level condition, and part (b) gives two sufficient conditions that can be checked directly, one on the error of the density and one on the error of the log-density.

\begin{theorem}[approximate likelihoods]\label{thm:general}
Suppose the assumptions of Theorem~\ref{thm:mle} hold. Let $\tilde p_n(x'|x;\theta)$ be real-valued and measurable in $(x,x',\theta)$, and let $\tilde\ell_n(\theta)=\sum_i\log\tilde p_n(X_i|X_{i-1};\theta)$, with $\log u=-\infty$ for $u\le0$.
\begin{enumerate}[label=(\alph*),nosep]
\item Suppose that $\sup_{\theta\in\Theta}|\tilde\ell_n(\theta)-\ell_n(\theta)|=o_\Pb(n\Delta_n)$ and that $\sup_{|\theta-\theta_0|\le\eta}|\tilde\ell_n(\theta)-\ell_n(\theta)|\toP0$ for some $\eta>0$. Then every $\tilde\theta_n$ with $\tilde\ell_n(\tilde\theta_n)\ge\sup_\Theta\tilde\ell_n-o_\Pb(1)$ satisfies $D_n(\tilde\theta_n-\hat\theta_n)\toP0$. In particular, $D_n(\tilde\theta_n-\theta_0)\toD N(0,\mathrm{diag}(\Gamma_\alpha^{-1},\Gamma_\beta^{-1}))$, and $\tilde\theta_n$ is efficient.
\item The conditions in (a) hold, with $\sup_{\theta\in\Theta}|\tilde\ell_n-\ell_n|\toP0$, if one of the following holds for constants $C$, $r$ and $\varrho$ that do not depend on $n$:
\begin{enumerate}[label=(\roman*),nosep]
\item \emph{density error:} for some $\gamma>(1+m\Lambda)/a$ and all $x,x',\theta$,
\[
\big|\tilde p_n(x'|x;\theta)-p_\Delta(x'|x;\theta)\big|\le C\jb{x}^r\Delta_n^{\gamma-m/2};
\]
\item \emph{log-density error:} $\tilde p_n>0$, and for some $\gamma>1/a$, $\varrho\ge0$ and all $x,x',\theta$,
\[
\big|\log\tilde p_n(x'|x;\theta)-\log p_\Delta(x'|x;\theta)\big|\le C\jb{x}^r\Delta_n^{\gamma}\big(1+|x'-x|^\varrho\Delta_n^{-\varrho/2}\big).
\]
\end{enumerate}
\end{enumerate}
\end{theorem}

\paragraph{Reading part (a).} Globally, the approximation need only be good enough to preserve consistency, which is why an error of order $o_\Pb(n\Delta_n)$ suffices; the stronger requirement $o_\Pb(1)$ is needed only near $\theta_0$. Part (a) is a direct consequence of the argmax argument behind Theorem~\ref{thm:mle}. Where $\tilde\ell_n$ is not continuous in $\theta$, the supremum and the probabilities in (a) are understood as outer ones.

\paragraph{Reading part (b).} Condition~(i) is the form in which Proposition~\ref{prop:unif} delivers accuracy: for the delta expansion of order $K$ it holds with $\gamma=(K+1)/2$ and reproduces the threshold of Theorem~\ref{thm:mle}(a), no more. Both conditions are worst-case bounds and should not be read as practical accuracy targets. Condition~(i) is very demanding, because an absolute error must be small compared with the density at the most extreme observed increment; in the bulk of the distribution it amounts to a relative accuracy of order $n^{-(1+m\Lambda)}$. Condition~(ii) works in relative terms from the start, and its threshold $1/a$ amounts roughly to a relative accuracy of order $n^{-1}$ per transition.

\paragraph{What is not covered.} Since (a) needs $o_\Pb(1)$ only near $\theta_0$, a local version of (i) with a smaller threshold should be possible. Stochastic approximations, such as simulated likelihoods, would require conditions on averaged errors and uniform control in $\theta$. We pursue neither, so Theorem~\ref{thm:general} should be read as a theoretical transfer principle rather than a practical recipe.

\paragraph{Relation to other density approximations.}
Section~5 of \citet{YCW2019} relates the delta expansion to earlier approximations. For a scalar diffusion with unit volatility, to which any scalar diffusion can be reduced by the Lamperti transform, the delta expansion with $\mu_0=0$ and order $2K$ coincides with the expansion of \citet{LeeSongLee2014} \citep[Proposition~5.1]{YCW2019}. That expansion, which its authors also call a delta expansion, rearranges the Hermite series of \citet{AitSahalia2002} by collecting terms in powers of $\Delta^{1/2}$.

In the multivariate case, \citet{YCW2019} rewrite the Malliavin-calculus expansion of \citet{Li2013} in terms of Hermite polynomials. They verify by symbolic computation that, in one and two dimensions, it agrees with the delta expansion with $\mu_0=0$ at every order they checked. \citet[Theorem~2]{ChenWanYang2026} settle the general case: in any dimension and at any order, the expansion of \citet{Li2013} has the same formula as the delta expansion with $\mu_0=0$, with a Hermite-polynomial expansion as the intermediate step. The two are therefore one expansion, written in two ways. The expansions of \citet{AitSahalia1999,AitSahalia2008} and \citet{Choi2013,Choi2015} are of a different kind: they approximate the log-density in integer powers of $\Delta$ and are not algebraically related to \eqref{eq:delta}. Theorem~\ref{thm:general} and Corollary~\ref{cor:others} below give conditions under which these approximations inherit the asymptotic theory of the exact MLE.

The next corollary applies the theorem to these expansions.

\begin{corollary}[other expansions]\label{cor:others}
Suppose the assumptions of Theorem~\ref{thm:mle} hold.
\begin{enumerate}[label=(\alph*),nosep]
\item For any centering $\mu_0(x;\theta)$ that is measurable in $(x,\theta)$ with $|\mu_0(x;\theta)|\le C\jb{x}$, in particular $\mu_0=0$, the delta expansion of order $K$ with $K+1>2(1+m\Lambda)/a$ satisfies condition~(i) of Theorem~\ref{thm:general}.
\item Let $m=1$ and let $F(x;\beta)=\int_0^xdu/\sigma(u;\beta)$ be the Lamperti transform. Suppose the drift $\mu_Y(y;\theta)$ of $Y=F(X;\beta)$ satisfies the spatial conditions on $\mu$ in Assumption~\ref{A:smooth}, uniformly in $\theta$: it is smooth in $y$, $|\mu_Y(y;\theta)|\le C\jb{y}$, and its spatial derivatives of every order $\ge1$ are bounded. (The separation of drift and diffusion parameters is not needed here, since $\mu_Y$ depends on both $\alpha$ and $\beta$.) Here $\Lambda$ is the ellipticity ratio of Assumption~\ref{A:ell} for $X$. Let $\tilde p_n(x'|x;\theta)=p_Y^{\mathrm{LSL}}(F(x')|F(x);\theta)/\sigma(x';\beta)$, where $p_Y^{\mathrm{LSL}}$ is the expansion of \citet{LeeSongLee2014} for $Y$ that coincides, by Proposition~5.1 of \citet{YCW2019} and in the indexing used there, with the delta expansion of $Y$ of order $2K'$ and centering $\mu_0=0$. If $2K'+1>2(1+\Lambda)/a$, then $\tilde p_n$ satisfies condition~(i) of Theorem~\ref{thm:general}.
\end{enumerate}
\end{corollary}

Part (a) covers the expansion of \citet{Li2013} as well. By the equivalence theorem of \citet{ChenWanYang2026}, it is the delta expansion with $\mu_0=0$, so at a sufficiently high order it satisfies condition~(i), and the approximate MLE of \citet{Li2013} inherits the asymptotic theory of the exact MLE in every dimension. Part (a) also shows that the centering does not matter for the equivalence with the exact MLE. It does matter at a fixed, low order: with $\mu_0=0$, the drift enters the corrections at order $\eps$ rather than $\eps^2$, and our consistency proof would have to be redone.

The condition on the drift of $Y$ in (b) is restrictive. Its derivative involves $\mu\sigma'/\sigma$, so with a linearly growing drift it essentially requires $\sigma'(x)=O(1/|x|)$. It holds, for instance, for $\mu(x;\alpha)=-\alpha x$ and $\sigma(x;\beta)=\beta(1+0.3\tanh x)$.

The log-density expansions of \citet{AitSahalia2008} and \citet{Choi2013,Choi2015} are natural candidates for condition~(ii). To our knowledge, however, no error bound of that form, uniform in the starting point and the parameter, is available for them.

\subsection{The delta-expansion estimator}\label{sec:approx}

We now analyze the delta-expansion estimator of a fixed order $K$, the estimator that is computed in practice. Its consistency and quadratic expansion are the tools behind Theorem~\ref{thm:mle}, and its bias shows how large $K$ must be.

\paragraph{Consistency.}
Consistency holds for every order. It depends only on the structure of the corrections, not on how accurately they approximate the true density, so it does not use Proposition~\ref{prop:unif}.

\begin{theorem}[consistency]\label{thm:cons}
Suppose Assumptions~\ref{A:smooth}--\ref{A:erg} and \ref{A:id}(a,b) hold, $n\Delta_n\to\infty$ and $\Delta_n\le Cn^{-a}$ for some $a>0$. Then for every $K\ge0$, $\Pb(\Omega_n)\to1$ and $\hat\theta^{(K)}_n\toP\theta_0$.
\end{theorem}

The proof (Appendix~\ref{app:cons}) separates the two time scales. Normalized by $n$, the contrast identifies the diffusion parameters, and the corrections are negligible. Normalized by $n\Delta_n$, the contrast identifies the drift parameters. At this scale, the corrections are as large as the signal and vanish only on average, because they are Hermite polynomials with mean zero under the local Gaussian law. This cancellation requires the diffusion parameters to be estimated consistently, so the proof treats them first.

\paragraph{Asymptotic distribution.}
For the asymptotic distribution, the order of the expansion matters, because it determines the bias of the score. With $D_n$ as in Section~\ref{sec:exact}, let
\[
\xi_n=\sqrt n\,\Delta_n^{(K+1)/2}.
\]
The bias is the projection of the first omitted correction $C_{K+1}$ on the Euler scores. For $K\ge1$ and $Z\sim N(0,\nu(x;\beta_0))$, let
\begin{equation}\label{eq:bias}
\begin{aligned}
b_\alpha(x)&=\E\big[(\partial_\alpha\mu)^\T\nu^{-1}Z\;C_{K+1}(Z;x,\theta_0)\big],\\
b_{\beta_j}(x)&=\tfrac12\E\big[\big(Z^\T\nu^{-1}\partial_{\beta_j}\nu\,\nu^{-1}Z-\tr(\nu^{-1}\partial_{\beta_j}\nu)\big)\,C_{K+1}(Z;x,\theta_0)\big],
\end{aligned}
\end{equation}
with all coefficients evaluated at $(x;\theta_0)$, and let $\bar B^{(K)}=\int(b_\alpha,b_\beta)\,d\pi$.

\begin{theorem}[asymptotic expansion]\label{thm:an}
Suppose Assumptions~\ref{A:smooth}--\ref{A:id} hold, $K\ge1$, $n\Delta_n\to\infty$ and $\Delta_n\le Cn^{-a}$ for some $a>0$. Then
\[
D_n\big(\hat\theta^{(K)}_n-\theta_0\big)=\Gamma^{-1}\big(S_n+\xi_n\bar B^{(K)}\big)+o_\Pb(1+\xi_n),
\]
where $S_n\toD N(0,\Gamma)$ is the normalized martingale part of the score (Lemma~\ref{lem:clt}). Moreover, $\bar B^{(K)}_\beta=0$ if $K$ is even and $\bar B^{(K)}_\alpha=0$ if $K$ is odd. The same expansion holds for the Euler estimator, $K=0$, with $\xi_n=\sqrt n\Delta_n$ and $\bar B^{(0)}:=\bar B^{(1)}$, whose drift block is zero.
\end{theorem}

The theorem describes all sampling regimes at once: the martingale term $S_n$ is always of order one, while the bias term $\xi_n\bar B^{(K)}$ may vanish, stay bounded or diverge, as Corollary~\ref{cor:regimes} makes explicit. The remainder is of smaller order because Proposition~\ref{prop:unif}, applied at order $K+1$, controls the second omitted term. The following consequences are immediate.

\begin{corollary}[sampling regimes]\label{cor:regimes}
Under the assumptions of Theorem~\ref{thm:an}, with $\xi_n^2=n\Delta_n^{K+1}$ for $K\ge1$ and $\xi_n^2=n\Delta_n^2$ for $K=0$:
\begin{enumerate}[label=(\alph*),nosep]
\item if $\xi_n\to0$, then $D_n(\hat\theta^{(K)}_n-\theta_0)\toD N\big(0,\mathrm{diag}(\Gamma_\alpha^{-1},\Gamma_\beta^{-1})\big)$;
\item if $\xi_n\to\xi\in(0,\infty)$, then $D_n(\hat\theta^{(K)}_n-\theta_0)\toD N\big(\xi\,\Gamma^{-1}\bar B^{(K)},\mathrm{diag}(\Gamma_\alpha^{-1},\Gamma_\beta^{-1})\big)$;
\item if $\xi_n\to\infty$ and $\bar B^{(K)}\neq0$, then $D_n(\hat\theta^{(K)}_n-\theta_0)/\xi_n\toP\Gamma^{-1}\bar B^{(K)}$, so $|D_n(\hat\theta^{(K)}_n-\theta_0)|\toP\infty$.
\end{enumerate}
\end{corollary}

Part (a) is the efficient limit: the rates are optimal, the two blocks are asymptotically independent, and the covariance attains the efficiency bound (Corollary~\ref{cor:lan}). Parts (b) and (c) show that, for the full parameter vector, the condition $\xi_n\to0$, that is, $n\Delta_n^{K+1}\to0$ for $K\ge1$, cannot be weakened whenever $\bar B^{(K)}\neq0$. In part (c), the expansion remains local, since $\xi_n/\sqrt{n\Delta_n}=\Delta_n^{K/2}\to0$ for $K\ge1$ and $\xi_n/\sqrt n=\Delta_n\to0$ for $K=0$. The bias constant is nonzero in all our examples: in closed form for $K=1$ and $K=2$ in the Ornstein--Uhlenbeck model below, and numerically for $K=3$. Each additional order relaxes the condition by a factor $\Delta_n$. For $K=1$, as for the Euler contrast ($K=0$), the condition is $n\Delta_n^2\to0$.

\paragraph{Comparison with the scalar case.}
Corollary~\ref{cor:regimes}(a) agrees with the scalar theory of \citet{ChangChen2011}. Their expansion of order $J$ in powers of $\Delta$ corresponds to $K=2J$ here, and for $\Delta_n\to0$ they show that the approximate MLE has the limit of the exact MLE when $J\ge2$, $n\Delta_n^3\to\infty$ and $n\Delta_n^{2J+1}\to0$. With $K=2J$, the last condition is exactly $\xi_n\to0$. Because expansions in integer powers of $\Delta$ correspond to even $K$ only, their first omitted term is always odd in $z$, and their leading bias always falls on the drift. The alternation between the two blocks for odd $K$ appears only because the delta expansion is graded in powers of $\Delta^{1/2}$. The present results remove the lower bound $n\Delta_n^3\to\infty$, apply in any dimension, and, through parts (b) and (c) and the parity statement, describe what happens when the condition fails.

\paragraph{Which parameters are biased.}
The parity statement in Theorem~\ref{thm:an} says on which parameters the leading bias falls. Let $K\ge1$. When $K$ is even, $C_{K+1}$ is a combination of Hermite polynomials of odd degree. These are orthogonal to the diffusion score, which is quadratic in $Z$, so the leading bias affects only the drift. When $K$ is odd, the roles are reversed, and the leading bias affects only the diffusion parameters. For the Euler contrast, $K=0$, the leading bias is that of $K=1$ and falls on the diffusion parameters. For the block that is spared, the score bias is of smaller order (Remark~\ref{rem:parity}). We expect that the weaker condition $n\Delta_n^{K+2}\to0$ then suffices for that block, but we do not prove this block-wise statement here.

\paragraph{An example.}
In the Ornstein--Uhlenbeck model $dX_t=-\kappa X_t\,dt+\mathrm{diag}(s_1,s_2)\,dW_t$ of Section~\ref{sec:mc}, $C_1=0$ and
\[
C_2=\tfrac12\big[z^\T\nu^{-1}E_1\nu^{-1}z-\tr(\nu^{-1}E_1)\big],\qquad E_1=-\tfrac12(\kappa\nu+\nu\kappa^\T).
\]
This gives $\bar B^{(1)}_{s_j}=-\kappa_{jj}/s_j\neq0$, so that $(\Gamma^{-1}\bar B^{(1)})_{s_j}\Delta_n=-\frac12\kappa_{jj}s_j\Delta_n$. This is the familiar downward bias of Euler volatility estimates, which understate $s_j$ by a fraction of about $\frac12\kappa_{jj}\Delta_n$.

\begin{remark}[bias orders by block]\label{rem:parity}
The proof gives the order of the score bias in each block. Define $\gamma$ by requiring that the bias contribution to the normalized score in that block is of order $\sqrt n\,\Delta_n^{\gamma}$. The orders are:
\begin{center}
\begin{tabular}{@{}lccccc@{}}
\toprule
$K$ & 0 & 1 & 2 & 3 & 4\\\midrule
drift, $\gamma_\alpha$ & $\mathbf{3/2}$ & $3/2$ & $\mathbf{3/2}$ & $5/2$ & $\mathbf{5/2}$\\
diffusion, $\gamma_\beta$ & $\mathbf{1}$ & $\mathbf{1}$ & $2$ & $\mathbf{2}$ & $3$\\\bottomrule
\end{tabular}
\end{center}
Bold entries are exact (given a nonzero bias constant); the others are upper bounds, up to factors that grow more slowly than any power of $n$ (Lemma~\ref{lem:bias}).

One caution applies to the non-bold entries. The normalized Hessian has off-diagonal blocks of order $\eps$, so the leading bias in one block also enters the estimate of the other block at order $\sqrt n\Delta_n^{(K+2)/2}$. Heuristically, the table therefore also gives the orders of the estimator bias, but in the non-leading block the constant combines two effects; Section~\ref{sec:mc} shows an instance.

The $K=0$ column is the classical Euler case, with the familiar conditions $n\Delta_n^3\to0$ for the drift and $n\Delta_n^2\to0$ for the diffusion parameters (Appendix~\ref{app:an}). Relative to the Euler contrast, $K=2$ improves only the diffusion parameters, and $K=3$ is needed to improve the drift.
\end{remark}

\begin{remark}[polynomial decay of $\Delta_n$]\label{rem:poly}
The condition $\Delta_n\le Cn^{-a}$ plays two roles. First, with a linearly growing drift, it controls the largest state in the sample: $\max_i|X_i|$ grows more slowly than any power of $n$, so $\eps\max_i\jb{X_i}^r\to0$. Second, in Theorem~\ref{thm:mle}, it makes a fixed-order expansion accurate enough at the most extreme observation.

For the delta-expansion estimator the condition costs nothing in practice: efficiency requires $n\Delta_n^{K+1}\to0$ anyway, which implies polynomial decay. For the exact MLE it is a limitation of the method. If $\Delta_n$ decays more slowly than every power of $n$, for instance $\Delta_n=1/\log n$, no fixed order suffices. LAN holds without the condition under the (different) assumptions of \citet{Gobet2002}, and we conjecture that Theorem~\ref{thm:mle} holds without it under Assumptions~\ref{A:smooth}--\ref{A:id}. In applications the condition is mild. Daily data, $\Delta_n=1/252$, with $n$ between $10^3$ and $10^5$ correspond to $\Delta_n=n^{-a}$ with $a$ between $0.48$ and $0.80$.
\end{remark}

\subsection{Moment convergence}\label{sec:moment}

Convergence in law does not imply convergence of moments. Moments are needed, for instance, to justify bias and mean squared error approximations. Because $\Theta$ is compact, $|D_n(\tilde\theta-\theta_0)|\le C\sqrt n$ for every $\Theta$-valued $\tilde\theta$. Moment convergence therefore follows once large deviations from $\theta_0$ have probability that decays faster than any power of $n$, and Assumption~\ref{A:rate} delivers this.

\begin{theorem}[moment convergence]\label{thm:moment}
Suppose Assumptions~\ref{A:smooth}--\ref{A:rate} hold and $c\,n^{-a_1}\le\Delta_n\le C\,n^{-a}$ for some $0<a\le a_1<1$. Let $N\sim N(0,\mathrm{diag}(\Gamma_\alpha^{-1},\Gamma_\beta^{-1}))$. Then, for every continuous $f$ with $|f(u)|\le C(1+|u|^p)$ for some $p$:
\begin{enumerate}[label=(\alph*),nosep]
\item $\E f\big(D_n(\hat\theta_n-\theta_0)\big)\to\E f(N)$, for an exact maximizer $\hat\theta_n$;
\item $\E f\big(D_n(\hat\theta^{(K)}_n-\theta_0)\big)\to\E f(N)$ for every $K\ge1$ with $n\Delta_n^{K+1}\to0$.
\end{enumerate}
\end{theorem}

The lower bound on $\Delta_n$ makes the time span $n\Delta_n$ grow polynomially, so that the drift contrast concentrates with polynomial probability bounds. This is needed both for the delta-expansion estimator and for the consistency step of the exact MLE. The proof (Appendix~\ref{app:moment}) repeats the arguments for Theorems~\ref{thm:cons}--\ref{thm:mle}, with every convergence in probability replaced by an $L^p$ bound. It is simpler than the quasi-likelihood analysis of \citet{Yoshida2011}, because compactness of $\Theta$ makes a large deviation inequality unnecessary.

\subsection{Volatility of linear growth}\label{sec:lin}

Assumptions~\ref{A:smooth} and \ref{A:ell} require the volatility to be bounded. This excludes models in which volatility rises with the level of the process, such as $dX_t=-\kappa X_t\,dt+s\sqrt{1+X_t^2}\,dW_t$, a member of the Pearson family of \citet{FormanSorensen2008} whose stationary law is a scaled Student distribution. This subsection shows that all results extend to volatility that is proportional to the size of the state: it grows like $1+|x|$ in every direction and is bounded away from zero. The two classes are not nested, because the new assumption also asks the drift to vary on the scale of the state.

\begin{assumption}[proportional volatility]\label{A:lin}
$\Theta$ and the parameterization are as in Assumption~\ref{A:smooth}, $\mu$ and $\sigma$ are $C^\infty$ in $(x,\theta)$, and for all multi-indices $a$ (in $x$) and $b$ (in $\theta$),
\[
|\partial^a_x\partial^b_\theta\mu(x;\alpha)|+|\partial^a_x\partial^b_\theta\sigma(x;\beta)|\le C_{a,b}\jb{x}^{1-|a|}.
\]
Moreover, there are constants $0<\lambda_{\min}\le\lambda_{\max}<\infty$ with $\lambda_{\min}\jb{x}^2I\le\nu(x;\beta)\le\lambda_{\max}\jb{x}^2I$ for all $x$ and $\beta$, and we write $\Lambda=\lambda_{\max}/\lambda_{\min}$.
\end{assumption}

The derivative bounds say that $\mu$ and $\sigma$ vary on the scale $\jb{x}$, which is the scale of the volatility itself. With $\jb{x}=(1+|x|^2)^{1/2}$, the example $\sigma(x)=s\sqrt{1+x^2}$ has $\nu=s^2\jb{x}^2$, so $\Lambda=1$. Linear drifts and $\sigma(x;\beta)=\beta\sqrt{1+|x|^2}\,\Sigma$, with $\Sigma$ a fixed invertible matrix, satisfy them. So do $\sigma(x;\beta)=\sqrt{1+|x|^2}\,\Sigma(x;\beta)$, where $\Sigma$ is smooth, bounded and uniformly elliptic with $|\partial^a_x\partial^b_\beta\Sigma|\le C\jb{x}^{-|a|}$ for all $a$ and $b$. The assumption excludes volatility that vanishes somewhere, as in CIR and Heston models, and volatility that grows in some directions but not in others, such as $\mathrm{diag}(\sqrt{1+x_1^2},\sqrt{1+x_2^2})$.

The key observation is that the model looks the same at every scale. Near a starting point $x$, the process moves on the scale $\jb{x}\sqrt\Delta$, and the rescaled process $(X_t-x)/\jb{x}$ has bounded, uniformly elliptic coefficients near the origin, with constants that do not depend on $x$. The delta expansion is equivariant under this affine change of variables. Hence the expansion error, measured relative to the local Gaussian scale, is bounded uniformly in the starting point, and no polynomial weight is needed.

\begin{proposition}[scale-invariant expansion error]\label{prop:lin}
Suppose Assumption~\ref{A:lin} holds, and let $\mu_0=\mu_0(x;\theta)$ be any centering with $|\mu_0(x;\theta)|\le C\jb{x}$. For every integer $K\ge1$ there is a constant $C_K$ such that, for all $\Delta\in(0,1]$, $\theta\in\Theta$, $x\in\R^m$ and $x'$ with $|x'-x|\le\jb{x}/8$,
\[
\big|p^{(K,\Delta)}(x'|x;\theta)-p_\Delta(x'|x;\theta)\big|\le C_K\big(\jb{x}^2\Delta\big)^{-m/2}\,\Delta^{(K+1)/2}.
\]
\end{proposition}

The bound is local, but observed transitions satisfy $|X_i-X_{i-1}|\le C\jb{X_{i-1}}\sqrt{\Delta_n\log n}$ with probability tending to one, so all of them fall in the region where it applies. The proof (Appendix~\ref{app:lin}) rescales and localizes the process, applies Proposition~\ref{prop:unif} to a modified process with bounded coefficients, and bounds the effect of the localization with the parabolic Harnack inequality.

Table~\ref{tab:linvol} illustrates the proposition for the scalar model above with $\kappa=1$ and $s=0.5$. Its exact transition density is computed by solving the Kolmogorov forward equation for the Lamperti transform $Y=\operatorname{asinh}(X)/s$, which has unit volatility and the bounded drift $-(\kappa/s+s/2)\tanh(sY)$. The maximal relative error of the delta expansion over $|x'-x-\mu(x)\Delta|\le3\sigma(x)\sqrt\Delta$ does not grow with the starting point. Between the two smallest values of $\Delta$, dividing $\Delta$ by four divides the error by $4.3$, $9.8$ and $17.6$ for $K=1,2,3$, close to the factors $4$, $8$ and $16$ implied by the rates $\Delta$, $\Delta^{3/2}$ and $\Delta^2$. At $\Delta=0.01$ the expansion is not yet in its asymptotic regime at three local standard deviations.

\begin{table}[t]
\centering\small
\caption{Maximal relative error of the delta expansion of order $K$ over $|x'-x-\mu(x)\Delta|\le3\sigma(x)\sqrt\Delta$, for $dX_t=-X_t\,dt+0.5\sqrt{1+X_t^2}\,dW_t$ started at $x$. The benchmark is a numerical solution of the forward equation; successive grid refinements change it by at most $10^{-5}$ in relative terms, and entries below $10^{-5}$ are not resolved. At $\Delta=0.01$ the window has half-width $3\sigma(x)\sqrt\Delta=0.15\jb{x}$, so $|x'-x|$ reaches about $0.16\jb{x}$ for $x\ge10$, slightly beyond the region of Proposition~\ref{prop:lin}, whose constant $\frac18$ is not sharp. At $x=0$ the odd-order corrections vanish by symmetry, so $K=2$ and $K=3$ coincide.}\label{tab:linvol}
\begin{tabular}{@{}llccc@{}}
\toprule
$\Delta$ & $x$ & $K=1$ & $K=2$ & $K=3$\\\midrule
\multirow{4}{*}{0.01} & 0 & $2.3\times10^{-2}$ & $3.1\times10^{-4}$ & $3.1\times10^{-4}$\\
 & 10 & $8.8\times10^{-2}$ & $5.5\times10^{-2}$ & $1.1\times10^{-2}$\\
 & 100 & $8.9\times10^{-2}$ & $5.5\times10^{-2}$ & $1.1\times10^{-2}$\\
 & 1000 & $8.9\times10^{-2}$ & $5.5\times10^{-2}$ & $1.1\times10^{-2}$\\
\addlinespace
\multirow{4}{*}{0.0025} & 0 & $5.7\times10^{-3}$ & $2.1\times10^{-5}$ & $2.1\times10^{-5}$\\
 & 10 & $1.9\times10^{-2}$ & $4.2\times10^{-3}$ & $4.9\times10^{-4}$\\
 & 100 & $1.9\times10^{-2}$ & $4.3\times10^{-3}$ & $5.0\times10^{-4}$\\
 & 1000 & $1.9\times10^{-2}$ & $4.3\times10^{-3}$ & $5.0\times10^{-4}$\\
\addlinespace
\multirow{4}{*}{0.000625} & 0 & $1.4\times10^{-3}$ & $<10^{-5}$ & $<10^{-5}$\\
 & 10 & $4.5\times10^{-3}$ & $4.3\times10^{-4}$ & $2.8\times10^{-5}$\\
 & 100 & $4.5\times10^{-3}$ & $4.4\times10^{-4}$ & $2.8\times10^{-5}$\\
 & 1000 & $4.5\times10^{-3}$ & $4.4\times10^{-4}$ & $2.8\times10^{-5}$\\
\bottomrule
\end{tabular}
\end{table}

With Proposition~\ref{prop:lin} in place of Proposition~\ref{prop:unif}, the rest of the theory carries over. Because every quantity is now naturally measured relative to $\jb{x}$, the polynomial weights disappear, and so does the need for moments of the stationary law. This matters here: in the scalar example, the stationary law has a Student-type tail and only finitely many moments.

\begin{theorem}[proportional volatility]\label{thm:lin}
Suppose Assumptions~\ref{A:lin}, \ref{A:erg}(i) and \ref{A:id} hold, with an arbitrary initial law, $n\Delta_n\to\infty$, and $\Delta_n\le Cn^{-a}$ for some $a\in(0,1)$. Then:
\begin{enumerate}[label=(\alph*),nosep]
\item Theorems~\ref{thm:mle}, \ref{thm:cons} and \ref{thm:an}, and Corollaries~\ref{cor:lan} and \ref{cor:regimes}, hold as stated, with $\Lambda$ from Assumption~\ref{A:lin};
\item Theorem~\ref{thm:general}(a) holds as stated, and Theorem~\ref{thm:general}(b) holds when conditions (i) and (ii) are replaced by
\begin{enumerate}[label=(\roman*$'$),nosep]
\item $|\tilde p_n(x'|x;\theta)-p_\Delta(x'|x;\theta)|\le C(\jb{x}^2\Delta_n)^{-m/2}\Delta_n^\gamma$,
\item $|\log\tilde p_n(x'|x;\theta)-\log p_\Delta(x'|x;\theta)|\le C\Delta_n^\gamma\big(1+|x'-x|^\varrho\jb{x}^{-\varrho}\Delta_n^{-\varrho/2}\big)$,
\end{enumerate}
both required only for $|x'-x|\le\jb{x}/8$, with the same thresholds on $\gamma$; Corollary~\ref{cor:others}(a) holds as stated;
\item Theorem~\ref{thm:moment} holds, under its two-sided condition $c\,n^{-a_1}\le\Delta_n\le Cn^{-a}$, if in addition Assumption~\ref{A:rate} holds for functions $g$ with $|g(x)|+\jb{x}|\nabla g(x)|\le C$.
\end{enumerate}
\end{theorem}

Corollary~\ref{cor:others}(b) is not claimed under Assumption~\ref{A:lin}. The moment condition in Assumption~\ref{A:erg}(ii) is not needed, because all functions that are averaged over the sample are bounded once increments are measured on the local scale. The information matrices in Assumption~\ref{A:id} are finite for the same reason: $\nu^{-1/2}\partial_\alpha\mu$ and $\nu^{-1}\partial_{\beta_j}\nu$ are bounded. The proof (Appendix~\ref{app:lin}) goes through the appendices and lists the changes.

\section{Monte Carlo evidence}\label{sec:mc}

The experiments address three questions. Are the limits in Theorems~\ref{thm:an} and \ref{thm:mle} accurate at realistic sample sizes? Does the bias of the delta-expansion estimator have the size and block structure that Theorem~\ref{thm:an} predicts? And what happens when the sampling interval is not small relative to the time scale of the model?

\subsection{Designs and implementation}

\paragraph{Ornstein--Uhlenbeck designs.}
Designs~1--3 use the bivariate Ornstein--Uhlenbeck model
\[
dX_t=-\kappa X_t\,dt+\mathrm{diag}(s_1,s_2)\,dW_t,\qquad \kappa=\begin{pmatrix}\kappa_{11}&0\\\kappa_{21}&\kappa_{22}\end{pmatrix},\qquad s_1=s_2=1.
\]
This is the bivariate OU specification of \citet{YCW2019}, with diffusion parameters added. Its transition density is Gaussian and can be computed exactly by the matrix-exponential method of \citet{VanLoan1978}, so the exact MLE is available as a benchmark. Paths are simulated exactly from this Gaussian transition, starting from the stationary law, so the OU designs carry no discretization error. The model satisfies Assumptions~\ref{A:smooth}--\ref{A:id}, since the symmetric part of $\kappa$ is positive definite in all three designs. Because $\nu$ is constant, $C_1\equiv0$, and $K=1$ coincides with the Euler contrast. Closed forms for $C_2$, $C_3$ and $C_4$ are given in Appendix~\ref{app:mc}. The three designs differ as follows:
\begin{itemize}[nosep]
\item Design~1 uses the parameter values of \citet{YCW2019}, $\kappa=(5,1,10)$, and $\Delta_n=n^{-0.6}$, so that $n\Delta_n^2\to0$.
\item Design~2 has slow mean reversion, $\kappa=(0.5,0.1,1)$, and $\Delta_n=n^{-0.4}$. Here $n\Delta_n^2\to\infty$, which puts the Euler contrast outside its classical regime, while $\kappa\Delta_n$ stays small.
\item Design~3 uses the values of Design~1 with $\Delta_n=n^{-0.4}$. Now $\kappa_{22}\Delta_n$ falls only from $0.63$ at $n=1000$ to $0.21$ at $n=16000$, so the sampling interval is not small relative to the fastest time scale of the model.
\end{itemize}

\paragraph{A nonlinear design.}
Design~4 uses a model with state-dependent volatility,
\[
\mu(x)=\begin{pmatrix}-\alpha_1\tanh x_1\\ \alpha_2\tanh x_1-\alpha_3\tanh x_2\end{pmatrix},\qquad
\sigma(x)=\begin{pmatrix}\beta_1(1+0.3\tanh x_2)&0\\0&\beta_2(1+0.3\tanh x_1)\end{pmatrix},
\]
with $(\alpha_1,\alpha_2,\alpha_3,\beta_1,\beta_2)=(1,0.5,1.5,0.7,0.9)$ and $\Delta_n=n^{-0.4}$. The drift is bounded and mean-reverting, and the volatility is bounded and uniformly elliptic. Since $\nu$ varies with the state, $C_1\neq0$, and every part of Lemma~\ref{lem:struct} is exercised. No closed-form density is available, so we compare the expansion orders $K=0,1,2$.

Paths are simulated by an Euler scheme with 40 sub-steps per observation interval, after a burn-in of 50 time units. The simulation error is of relative order $\Delta_n/40$. As an order-of-magnitude estimate, $\sqrt n\Delta_n/40\approx0.07$ standardized units for the diffusion parameters at $n=16000$. A direct check with 160 sub-steps, reported with the results, confirms that the effect is smaller still. The matrix $\Gamma$ is computed from a separate Euler path of $4\times10^5$ observations with spacing $0.05$ and 10 sub-steps, which spans $2\times10^4$ time units.

\paragraph{Implementation.}
The coefficients $C_k$ are computed in closed form, using for Design~4 the tensor formulas of Appendix~\ref{app:mc}. Every closed form was checked against a direct numerical implementation of the recursion of \citet{YCW2019} at randomly drawn $(x,z,\theta)$, and the two agree to machine precision. All estimators maximize the truncated contrast \eqref{eq:contrast}. We use L-BFGS-B, followed by a Nelder--Mead polishing step whenever L-BFGS-B does not report convergence.

To mimic practice, the optimizer does not start at the true value. For the OU designs it starts from the Euler estimate, which has a closed form by least squares. For Design~4 it starts from the Euler estimate, itself computed from a crude initial value ($\alpha=(0.5,0,1)$ and $\beta$ from realized variances). After polishing, every fit converged in all four designs. In Design~4, L-BFGS-B flagged two of the 900 fits, both for $K\ge1$ and none for the Euler estimator; polishing changed no entry of Table~\ref{tab:nl} at the reported precision. The optimizer is run without bounds, with the volatility parameters on the log scale, so the compactness of $\Theta$ is not imposed in the simulations.

\paragraph{Reported statistics.}
We report the Monte Carlo mean and standard deviation of the standardized errors $\sqrt{n\Delta_n}(\hat\alpha_j-\alpha_{0j})/\mathrm{sd}_j$ and $\sqrt n(\hat\beta_j-\beta_{0j})/\mathrm{sd}_j$, where $\mathrm{sd}_j=[(\Gamma^{-1})_{jj}]^{1/2}$ with $\Gamma$ evaluated at $\theta_0$. In Designs~1--3 the parameters are the entries of $\kappa$ and $s$. If the asymptotic theory is accurate, these statistics have mean $0$ and standard deviation $1$.

All estimators are computed on the same simulated paths. Their differences from the exact MLE therefore carry much less Monte Carlo noise than the errors themselves, and they isolate the effect of the approximation (Table~\ref{tab:paired}). Table~\ref{tab:cov} reports the coverage of nominal 95\% Wald intervals based on the plug-in estimate of $\Gamma$ at $\hat\theta$, which is what a practitioner would compute. Seeds differ across designs.

The results are organized as follows. We first show that the exact MLE is close to its Gaussian limit, then use paired differences to measure the approximation bias of each estimator in each block, and compare it with the prediction of Theorem~\ref{thm:an}. Designs~3 and 4 then probe a coarse sampling interval and state-dependent volatility.

\subsection{Results}

\paragraph{The exact MLE.}
In Designs~1 and 2, the exact MLE is close to its Gaussian limit, with coverage between 90\% and 98\% (Tables~\ref{tab:ou12} and \ref{tab:cov}). It has a small upward finite-span bias of $0.1$--$0.3$ standard deviations in the mean-reversion parameters $\kappa_{11}$ and $\kappa_{22}$. The paired differences in Table~\ref{tab:paired} separate this bias from the approximation error, and we focus on them below.

\begin{table}[!htbp]\centering\footnotesize
\caption{Standardized errors: Monte Carlo mean (standard deviation); the target is $0\,(1)$. Ornstein--Uhlenbeck model, 200 replications; the Monte Carlo standard error of a mean is about $0.07$.}\label{tab:ou12}
\setlength{\tabcolsep}{3.5pt}\begin{tabular}{@{}llccccc@{}}\toprule
$n$ ($n\Delta_n$) & estimator & $\kappa_{11}$ & $\kappa_{21}$ & $\kappa_{22}$ & $s_1$ & $s_2$\\\midrule
\multicolumn{7}{@{}l}{\emph{Design 1: $\kappa=(5,1,10)$, $\Delta_n=n^{-0.6}$}}\\
1000 (16) & exact MLE & $+0.15$\,(1.07) & $-0.05$\,(1.14) & $+0.19$\,(1.13) & $-0.17$\,(0.97) & $+0.10$\,(1.16) \\
 & Euler ($K=0$) & $-0.11$\,(0.99) & $-0.17$\,(1.03) & $-0.52$\,(0.96) & $-1.92$\,(0.88) & $-3.30$\,(0.97) \\
 & $K=2$ & $-0.05$\,(1.03) & $-0.13$\,(1.06) & $-0.07$\,(1.15) & $-0.16$\,(0.97) & $+0.13$\,(1.15) \\
 & $K=3$ & $+0.25$\,(1.15) & $+0.07$\,(1.27) & $+1.11$\,(1.59) & $-0.09$\,(0.98) & $+0.52$\,(1.21) \\
 & $K=4$ & $+0.17$\,(1.08) & $-0.03$\,(1.16) & $+0.36$\,(1.24) & $-0.17$\,(0.97) & $+0.09$\,(1.16) \\
\addlinespace
4000 (28) & exact MLE & $+0.28$\,(1.05) & $+0.06$\,(0.99) & $+0.17$\,(1.07) & $-0.02$\,(1.05) & $-0.01$\,(1.01) \\
 & Euler ($K=0$) & $+0.12$\,(1.01) & $-0.02$\,(0.96) & $-0.24$\,(0.99) & $-1.59$\,(1.02) & $-3.05$\,(0.95) \\
 & $K=2$ & $+0.12$\,(1.01) & $-0.03$\,(0.94) & $-0.23$\,(1.01) & $-0.01$\,(1.05) & $+0.03$\,(1.01) \\
 & $K=3$ & $+0.28$\,(1.05) & $+0.06$\,(1.00) & $+0.20$\,(1.09) & $+0.02$\,(1.05) & $+0.15$\,(1.02) \\
 & $K=4$ & $+0.28$\,(1.05) & $+0.06$\,(0.99) & $+0.18$\,(1.07) & $-0.02$\,(1.05) & $-0.01$\,(1.01) \\
\addlinespace
16000 (48) & exact MLE & $+0.17$\,(0.96) & $-0.05$\,(0.96) & $+0.22$\,(0.96) & $+0.09$\,(1.03) & $+0.14$\,(1.05) \\
 & Euler ($K=0$) & $+0.08$\,(0.95) & $-0.10$\,(0.95) & $-0.01$\,(0.93) & $-1.27$\,(1.01) & $-2.55$\,(1.04) \\
 & $K=2$ & $+0.08$\,(0.95) & $-0.09$\,(0.94) & $-0.01$\,(0.93) & $+0.09$\,(1.03) & $+0.15$\,(1.05) \\
 & $K=3$ & $+0.17$\,(0.96) & $-0.04$\,(0.96) & $+0.23$\,(0.96) & $+0.10$\,(1.03) & $+0.20$\,(1.06) \\
 & $K=4$ & $+0.17$\,(0.96) & $-0.04$\,(0.96) & $+0.23$\,(0.96) & $+0.08$\,(1.03) & $+0.14$\,(1.05) \\
\addlinespace
\multicolumn{7}{@{}l}{\emph{Design 2: $\kappa=(0.5,0.1,1)$, $\Delta_n=n^{-0.4}$}}\\
1000 (63) & exact MLE & $+0.14$\,(1.00) & $+0.01$\,(1.08) & $+0.27$\,(1.07) & $+0.13$\,(1.03) & $+0.02$\,(1.05) \\
 & Euler ($K=0$) & $+0.07$\,(0.97) & $-0.03$\,(1.01) & $+0.07$\,(1.00) & $-0.60$\,(1.00) & $-1.42$\,(1.00) \\
 & $K=2$ & $+0.07$\,(0.97) & $-0.03$\,(1.02) & $+0.08$\,(1.01) & $+0.13$\,(1.03) & $+0.04$\,(1.05) \\
 & $K=3$ & $+0.14$\,(1.00) & $+0.01$\,(1.08) & $+0.29$\,(1.09) & $+0.14$\,(1.03) & $+0.09$\,(1.06) \\
 & $K=4$ & $+0.14$\,(1.00) & $+0.01$\,(1.08) & $+0.28$\,(1.07) & $+0.13$\,(1.03) & $+0.02$\,(1.05) \\
\addlinespace
4000 (145) & exact MLE & $+0.22$\,(0.94) & $-0.14$\,(1.04) & $+0.13$\,(0.99) & $+0.03$\,(1.02) & $-0.10$\,(0.97) \\
 & Euler ($K=0$) & $+0.16$\,(0.93) & $-0.18$\,(1.03) & $-0.02$\,(0.95) & $-0.81$\,(1.00) & $-1.72$\,(0.94) \\
 & $K=2$ & $+0.16$\,(0.93) & $-0.17$\,(1.01) & $-0.03$\,(0.95) & $+0.03$\,(1.02) & $-0.09$\,(0.97) \\
 & $K=3$ & $+0.22$\,(0.94) & $-0.14$\,(1.04) & $+0.13$\,(0.99) & $+0.04$\,(1.02) & $-0.06$\,(0.98) \\
 & $K=4$ & $+0.22$\,(0.94) & $-0.14$\,(1.04) & $+0.14$\,(0.99) & $+0.03$\,(1.02) & $-0.10$\,(0.97) \\
\addlinespace
16000 (333) & exact MLE & $+0.21$\,(1.08) & $-0.01$\,(1.01) & $+0.10$\,(1.00) & $-0.03$\,(0.97) & $-0.01$\,(0.96) \\
 & Euler ($K=0$) & $+0.16$\,(1.07) & $-0.05$\,(0.97) & $-0.04$\,(0.98) & $-0.98$\,(0.96) & $-1.87$\,(0.94) \\
 & $K=2$ & $+0.16$\,(1.07) & $-0.04$\,(0.99) & $-0.03$\,(0.98) & $-0.03$\,(0.97) & $0.00$\,(0.96) \\
 & $K=3$ & $+0.21$\,(1.08) & $-0.01$\,(1.01) & $+0.10$\,(1.00) & $-0.03$\,(0.97) & $+0.02$\,(0.96) \\
 & $K=4$ & $+0.21$\,(1.08) & $-0.01$\,(1.01) & $+0.10$\,(1.00) & $-0.03$\,(0.97) & $-0.01$\,(0.96) \\
\bottomrule\end{tabular}
\end{table}

\paragraph{Designs 1 and 2: diffusion parameters.}
The Euler estimator of the diffusion parameters is badly biased. Relative to the exact MLE, its standardized bias for $s_2$ is between $-2.7$ and $-3.4$ in Design~1 and between $-1.4$ and $-1.9$ in Design~2 (Table~\ref{tab:paired}). Its Wald intervals for $s_2$ cover the true value only 9--67\% of the time.

The bias evolves as the theory predicts. Its order $\sqrt n\Delta_n$ is proportional to $n^{-0.1}$ in Design~1 and to $n^{0.1}$ in Design~2, which implies factors of $0.87$ and $1.15$ per quadrupling of $n$. For $s_2$, the observed factors are $0.89$ and $0.88$ in Design~1, and $1.13$ and $1.15$ in Design~2.

The $K=2$ expansion removes the diffusion bias. Its paired differences from the exact MLE are at most $0.04$, and its coverage is between 90\% and 98\%. The $K=3$ estimator does almost as well once $n\ge4000$, with differences of at most $0.16$, but at $n=1000$ in Design~1 its difference for $s_2$ is $0.43$. At that sample size the higher-order estimators are also more variable: the standard deviation of the $K=3$ estimate of $\kappa_{22}$ is $1.59$ in Design~1.

Why does $K=3$ do worse than $K=2$ here? For both orders, the score bias in the diffusion block is of order $\sqrt n\Delta_n^2$ and has the same constant. The difference comes, heuristically, from the cross-block effect noted in Remark~\ref{rem:parity}. For $K=2$, the drift bias feeds into the diffusion estimates through the cross-block Hessian, at order $\sqrt n\Delta_n^{(K+2)/2}=\sqrt n\Delta_n^2$, the same order as the score bias, and largely offsets it. We quantify this by computing, by Gauss--Hermite quadrature under the exact transition, the expected score and the maximizer of the expected contrast. In Design~1 at $n=16000$, the score bias alone implies a standardized bias of $0.055$ in $s_2$ for both orders, and the cross term contributes $-0.041$ for $K=2$ and essentially zero for $K=3$. This matches the simulated values of $0.01$ and $0.06$.

\begin{table}[!htbp]\centering\footnotesize
\caption{Paired differences from the exact MLE: Monte Carlo mean of (estimator minus exact MLE) in standardized units, with its Monte Carlo standard error in parentheses. Same replications as Table~\ref{tab:ou12}.}\label{tab:paired}
\setlength{\tabcolsep}{4pt}\begin{tabular}{@{}llccccc@{}}\toprule
$n$ & estimator & $\kappa_{11}$ & $\kappa_{21}$ & $\kappa_{22}$ & $s_1$ & $s_2$\\\midrule
\multicolumn{7}{@{}l}{\emph{Design 1}}\\
1000 & Euler ($K=0$) & $-0.26$\,(0.01) & $-0.12$\,(0.02) & $-0.71$\,(0.01) & $-1.75$\,(0.02) & $-3.40$\,(0.03) \\
 & $K=2$ & $-0.20$\,(0.01) & $-0.08$\,(0.01) & $-0.26$\,(0.01) & $+0.01$\,(0.00) & $+0.03$\,(0.00) \\
 & $K=3$ & $+0.09$\,(0.01) & $+0.12$\,(0.01) & $+0.92$\,(0.04) & $+0.08$\,(0.00) & $+0.43$\,(0.01) \\
 & $K=4$ & $+0.01$\,(0.00) & $+0.02$\,(0.00) & $+0.17$\,(0.01) & $0.00$\,(0.00) & $-0.01$\,(0.00) \\
\addlinespace
4000 & Euler ($K=0$) & $-0.15$\,(0.00) & $-0.08$\,(0.01) & $-0.41$\,(0.01) & $-1.57$\,(0.01) & $-3.04$\,(0.02) \\
 & $K=2$ & $-0.15$\,(0.00) & $-0.09$\,(0.00) & $-0.39$\,(0.01) & $+0.01$\,(0.00) & $+0.04$\,(0.00) \\
 & $K=3$ & $0.00$\,(0.00) & $0.00$\,(0.00) & $+0.04$\,(0.00) & $+0.04$\,(0.00) & $+0.16$\,(0.00) \\
 & $K=4$ & $0.00$\,(0.00) & $0.00$\,(0.00) & $+0.01$\,(0.00) & $0.00$\,(0.00) & $0.00$\,(0.00) \\
\addlinespace
16000 & Euler ($K=0$) & $-0.09$\,(0.00) & $-0.06$\,(0.01) & $-0.24$\,(0.00) & $-1.36$\,(0.01) & $-2.69$\,(0.01) \\
 & $K=2$ & $-0.09$\,(0.00) & $-0.05$\,(0.00) & $-0.24$\,(0.00) & $0.00$\,(0.00) & $+0.01$\,(0.00) \\
 & $K=3$ & $0.00$\,(0.00) & $0.00$\,(0.00) & $0.00$\,(0.00) & $+0.01$\,(0.00) & $+0.06$\,(0.00) \\
 & $K=4$ & $0.00$\,(0.00) & $0.00$\,(0.00) & $0.00$\,(0.00) & $0.00$\,(0.00) & $0.00$\,(0.00) \\
\addlinespace
\multicolumn{7}{@{}l}{\emph{Design 2}}\\
1000 & Euler ($K=0$) & $-0.07$\,(0.00) & $-0.04$\,(0.01) & $-0.20$\,(0.01) & $-0.72$\,(0.01) & $-1.44$\,(0.02) \\
 & $K=2$ & $-0.07$\,(0.00) & $-0.04$\,(0.00) & $-0.19$\,(0.01) & $0.00$\,(0.00) & $+0.02$\,(0.00) \\
 & $K=3$ & $0.00$\,(0.00) & $0.00$\,(0.00) & $+0.01$\,(0.00) & $+0.02$\,(0.00) & $+0.07$\,(0.00) \\
 & $K=4$ & $0.00$\,(0.00) & $0.00$\,(0.00) & $0.00$\,(0.00) & $0.00$\,(0.00) & $0.00$\,(0.00) \\
\addlinespace
4000 & Euler ($K=0$) & $-0.06$\,(0.00) & $-0.03$\,(0.01) & $-0.16$\,(0.00) & $-0.83$\,(0.01) & $-1.62$\,(0.01) \\
 & $K=2$ & $-0.06$\,(0.00) & $-0.03$\,(0.00) & $-0.16$\,(0.00) & $0.00$\,(0.00) & $+0.01$\,(0.00) \\
 & $K=3$ & $0.00$\,(0.00) & $0.00$\,(0.00) & $0.00$\,(0.00) & $+0.01$\,(0.00) & $+0.04$\,(0.00) \\
 & $K=4$ & $0.00$\,(0.00) & $0.00$\,(0.00) & $0.00$\,(0.00) & $0.00$\,(0.00) & $0.00$\,(0.00) \\
\addlinespace
16000 & Euler ($K=0$) & $-0.05$\,(0.00) & $-0.04$\,(0.01) & $-0.14$\,(0.00) & $-0.95$\,(0.01) & $-1.86$\,(0.01) \\
 & $K=2$ & $-0.05$\,(0.00) & $-0.03$\,(0.00) & $-0.14$\,(0.00) & $0.00$\,(0.00) & $+0.01$\,(0.00) \\
 & $K=3$ & $0.00$\,(0.00) & $0.00$\,(0.00) & $0.00$\,(0.00) & $+0.01$\,(0.00) & $+0.03$\,(0.00) \\
 & $K=4$ & $0.00$\,(0.00) & $0.00$\,(0.00) & $0.00$\,(0.00) & $0.00$\,(0.00) & $0.00$\,(0.00) \\
\bottomrule\end{tabular}
\end{table}

\paragraph{Designs 1 and 2: drift parameters.}
The drift parameters show the block structure of Remark~\ref{rem:parity} even more clearly. Once $n\ge4000$, the Euler and $K=2$ estimators have nearly identical drift biases relative to the exact MLE: $-0.41$ and $-0.39$ for $\kappa_{22}$ at $n=4000$ in Design~1, and $-0.24$ for both at $n=16000$. This is what the theory predicts, since $\gamma_\alpha=3/2$ for both orders.

The rate of decay also matches. For the Euler estimator in Design~1, the bias falls by a factor of $0.58$ per quadrupling of $n$, against the predicted $4^{-0.4}=0.57$ for an order of $\sqrt n\Delta_n^{3/2}$. In Design~2 it falls by factors of $0.79$ and $0.87$, against the predicted $4^{-0.1}=0.87$.

Going one order further removes the drift bias. The $K=3$ estimator, for which $\gamma_\alpha=5/2$, has paired differences of at most $0.04$ for $n\ge4000$, and the $K=4$ estimator agrees with the exact MLE to within $0.01$. The only sizeable departures from these patterns occur at $n=1000$ in Design~1, where $\kappa_{22}\Delta_n=0.16$ and higher-order terms are still visible.

In these two designs, the truncation in \eqref{eq:contrast} is rarely active. Even at $n=1000$, the average fraction of observations with $|R_i(\hat\theta)|>\frac14$, which we do not tabulate, is at most 3.9\% in Design~1 and 0.4\% in Design~2.

\paragraph{Predicted versus simulated bias.}
Theorem~\ref{thm:an} predicts the size of the approximation bias, not only its order. In the OU designs, the constants $\bar B^{(K)}$ can be computed by Gaussian quadrature. Table~\ref{tab:biaspred} compares the predicted bias $\xi_n(\Gamma^{-1}\bar B^{(K)})_j/\mathrm{sd}_j$ with the simulated paired differences, and the agreement is close. The table also reports the all-order bias at fixed $\Delta_n$, computed from the maximizer of the expected contrast under the exact transition, which includes every order of the expansion error but no finite-sample noise.

For the Euler estimator of the volatilities, the bias constant has the closed form $\bar B^{(1)}_{s_j}=-\kappa_{jj}/s_j$, so the predicted standardized bias is $-\sqrt{n}\Delta_n\kappa_{jj}/\sqrt2$. It is within $0.05$ of the simulated value in every case except $s_2$ at $n=1000$ in Design~1, where it is within $0.15$. For the drift bias of the $K=2$ estimator and the diffusion bias of the $K=3$ estimator, prediction and simulation differ by at most $0.02$ in Design~2 and for $n\ge4000$ in Design~1. The exception is again $n=1000$ in Design~1, most of all for $\kappa_{22}$, where the leading term predicts $-0.71$, the all-order bias is $-0.40$ and the simulation gives $-0.26$. About two-thirds of this discrepancy is due to terms of higher order in $\Delta_n$, which are visible because $\kappa_{22}\Delta_n=0.16$; the rest is finite-sample variation over a time span of only $n\Delta_n=16$.

Finally, the Euler drift bias in Table~\ref{tab:paired} agrees with the $K=2$ prediction, as the heuristic in Remark~\ref{rem:parity} suggests: for both orders, the leading drift bias is the projection of $C_3$ on the drift score.

\begin{table}[!htbp]
\centering\small
\caption{Approximation bias in Designs~1 and 2, in standardized units. ``pred'' is the leading bias $\xi_n(\Gamma^{-1}\bar B^{(K)})_j/\mathrm{sd}_j$ of Theorem~\ref{thm:an}; ``all'' is the all-order bias $D_n(\theta^*_K-\theta_0)/\mathrm{sd}$, where $\theta^*_K$ maximizes the expected contrast under the exact transition; ``MC'' is the Monte Carlo mean of the paired difference between the estimator and the exact MLE. Only the block that carries the leading bias is shown. In these designs the Euler estimator coincides with $K=1$.}\label{tab:biaspred}
\setlength{\tabcolsep}{3.5pt}
\begin{tabular}{@{}cr*{12}{r}@{}}
\multicolumn{14}{@{}l}{\emph{Panel A: diffusion parameters}}\\\toprule
&& \multicolumn{6}{c}{Euler ($K=0$)} & \multicolumn{6}{c}{$K=3$}\\\cmidrule(lr){3-8}\cmidrule(l){9-14}
&& \multicolumn{3}{c}{$s_1$} & \multicolumn{3}{c}{$s_2$} & \multicolumn{3}{c}{$s_1$} & \multicolumn{3}{c}{$s_2$}\\
Design & $n$ & pred & all & MC & pred & all & MC & pred & all & MC & pred & all & MC\\\midrule
1 & 1{,}000 & $-$1.77 & $-$1.71 & $-$1.75 & $-$3.54 & $-$3.32 & $-$3.40 & \phantom{$-$}0.09 & \phantom{$-$}0.08 & \phantom{$-$}0.08 & \phantom{$-$}0.38 & \phantom{$-$}0.41 & \phantom{$-$}0.43\\
 & 4{,}000 & $-$1.54 & $-$1.52 & $-$1.57 & $-$3.09 & $-$3.00 & $-$3.04 & \phantom{$-$}0.04 & \phantom{$-$}0.04 & \phantom{$-$}0.04 & \phantom{$-$}0.14 & \phantom{$-$}0.15 & \phantom{$-$}0.16\\
 & 16{,}000 & $-$1.34 & $-$1.33 & $-$1.36 & $-$2.69 & $-$2.65 & $-$2.69 & \phantom{$-$}0.01 & \phantom{$-$}0.01 & \phantom{$-$}0.01 & \phantom{$-$}0.05 & \phantom{$-$}0.06 & \phantom{$-$}0.06\\
\addlinespace
2 & 1{,}000 & $-$0.71 & $-$0.70 & $-$0.72 & $-$1.41 & $-$1.37 & $-$1.44 & \phantom{$-$}0.01 & \phantom{$-$}0.02 & \phantom{$-$}0.02 & \phantom{$-$}0.06 & \phantom{$-$}0.06 & \phantom{$-$}0.07\\
 & 4{,}000 & $-$0.81 & $-$0.80 & $-$0.83 & $-$1.62 & $-$1.60 & $-$1.62 & \phantom{$-$}0.01 & \phantom{$-$}0.01 & \phantom{$-$}0.01 & \phantom{$-$}0.04 & \phantom{$-$}0.04 & \phantom{$-$}0.04\\
 & 16{,}000 & $-$0.93 & $-$0.93 & $-$0.95 & $-$1.86 & $-$1.85 & $-$1.86 & \phantom{$-$}0.01 & \phantom{$-$}0.01 & \phantom{$-$}0.01 & \phantom{$-$}0.03 & \phantom{$-$}0.03 & \phantom{$-$}0.03\\
\bottomrule
\end{tabular}

\bigskip
\begin{tabular}{@{}cr*{9}{r}@{}}
\multicolumn{11}{@{}l}{\emph{Panel B: drift parameters}}\\\toprule
&& \multicolumn{9}{c}{$K=2$}\\\cmidrule(l){3-11}
&& \multicolumn{3}{c}{$\kappa_{11}$} & \multicolumn{3}{c}{$\kappa_{21}$} & \multicolumn{3}{c}{$\kappa_{22}$}\\
Design & $n$ & pred & all & MC & pred & all & MC & pred & all & MC\\\midrule
1 & 1{,}000 & $-$0.25 & $-$0.19 & $-$0.20 & $-$0.15 & $-$0.08 & $-$0.08 & $-$0.71 & $-$0.40 & $-$0.26\\
 & 4{,}000 & $-$0.14 & $-$0.14 & $-$0.15 & $-$0.09 & $-$0.08 & $-$0.09 & $-$0.41 & $-$0.39 & $-$0.39\\
 & 16{,}000 & $-$0.08 & $-$0.08 & $-$0.09 & $-$0.05 & $-$0.05 & $-$0.05 & $-$0.23 & $-$0.23 & $-$0.24\\
\addlinespace
2 & 1{,}000 & $-$0.06 & $-$0.06 & $-$0.07 & $-$0.04 & $-$0.04 & $-$0.04 & $-$0.18 & $-$0.17 & $-$0.19\\
 & 4{,}000 & $-$0.05 & $-$0.05 & $-$0.06 & $-$0.03 & $-$0.03 & $-$0.03 & $-$0.15 & $-$0.15 & $-$0.16\\
 & 16{,}000 & $-$0.05 & $-$0.05 & $-$0.05 & $-$0.03 & $-$0.03 & $-$0.03 & $-$0.13 & $-$0.13 & $-$0.14\\
\bottomrule
\end{tabular}
\end{table}

\paragraph{Design 3: a coarse sampling interval.}
When the sampling interval is not small relative to $1/\kappa_{22}$, the picture changes (Table~\ref{tab:ou3}). The exact MLE is still well behaved. The Euler estimator, however, is badly biased in both blocks: its standardized biases are about $-4$ for $\kappa_{22}$ and between $-11$ and $-17$ for $s_2$, and its Wald intervals for $\kappa_{22}$, $s_1$ and $s_2$ never contain the truth.

The expansions are now used well beyond their range of validity. At $n=1000$, between 57\% and 76\% of the observations have $|R_i(\hat\theta)|>\frac14$, and the floor $\frac12$ in \eqref{eq:contrast} binds for 5--7\% of them. Higher orders do not help, and they are also more variable, with standardized standard deviations of up to $1.94$. For $\kappa_{22}$, the odd order $K=3$ overshoots, with standardized biases of $+6.5$ to $+9.8$, and $K=4$ overshoots less ($+1.6$ to $+3.9$). The $K=2$ estimator, whose drift bias has the same leading order as Euler's, ends up the least biased ($-1.6$ to $-0.15$). As $n$ grows and $\kappa_{22}\Delta_n$ falls to $0.21$, the fraction of large corrections drops to about 5\%. The $K=2$ and $K=4$ estimators then recover the diffusion parameters, but the drift remains distorted for $K=3$ and $K=4$, and even $K=2$ still has a bias of $-0.78$ for $\kappa_{11}$, with coverage of 84\%.

This suggests a practical heuristic, which we have not calibrated beyond these designs. The fraction of observations with $|R_i(\hat\theta)|>\frac14$ is a cheap diagnostic, available in any application. When it is large, as in Design~3 at $n=1000$, the sampling interval is too coarse for the expansion, and raising the order does not help. Fractions of a few percent, as in Design~1 at $n=1000$, were harmless in our experiments.

\begin{table}[!htbp]\centering\footnotesize
\caption{Standardized errors: Monte Carlo mean (standard deviation); the target is $0\,(1)$. 200 replications. The last column is the average percentage of observations with $|R_i(\hat\theta)|>\frac14$.}\label{tab:ou3}
\setlength{\tabcolsep}{3.5pt}\begin{tabular}{@{}llcccccr@{}}\toprule
$n$ ($n\Delta_n$) & estimator & $\kappa_{11}$ & $\kappa_{21}$ & $\kappa_{22}$ & $s_1$ & $s_2$ & $|R_i|>\frac14$ (\%)\\\midrule
\multicolumn{8}{@{}l}{\emph{Design 3: $\kappa=(5,1,10)$, $\Delta_n=n^{-0.4}$ ($\kappa_{22}\Delta_n=0.63,\,0.36,\,0.21$)}}\\
1000 (63) & exact MLE & $+0.07$\,(1.13) & $-0.02$\,(1.15) & $+0.17$\,(1.41) & $0.00$\,(1.14) & $+0.02$\,(1.22) & -- \\
 & Euler ($K=0$) & $-1.75$\,(0.82) & $-1.00$\,(0.77) & $-4.53$\,(0.74) & $-6.24$\,(0.82) & $-11.05$\,(0.70) & -- \\
 & $K=2$ & $-0.41$\,(1.08) & $-0.59$\,(0.93) & $-1.60$\,(0.95) & $-1.28$\,(1.04) & $-3.61$\,(0.94) & 57.0 \\
 & $K=3$ & $+1.44$\,(1.50) & $+0.67$\,(1.52) & $+6.50$\,(1.81) & $-1.14$\,(1.08) & $-0.97$\,(1.21) & 76.1 \\
 & $K=4$ & $+0.64$\,(1.39) & $+0.25$\,(1.42) & $+2.94$\,(1.60) & $-0.93$\,(1.10) & $-3.34$\,(1.05) & 67.4 \\
\addlinespace
4000 (145) & exact MLE & $+0.07$\,(1.08) & $+0.02$\,(1.11) & $+0.07$\,(1.18) & $-0.01$\,(1.10) & $-0.02$\,(1.20) & -- \\
 & Euler ($K=0$) & $-1.57$\,(0.90) & $-0.87$\,(0.84) & $-4.31$\,(0.82) & $-7.56$\,(0.96) & $-14.04$\,(0.86) & -- \\
 & $K=2$ & $-0.32$\,(1.10) & $-0.26$\,(1.02) & $-0.15$\,(1.11) & $-0.56$\,(1.08) & $-1.79$\,(1.08) & 17.6 \\
 & $K=3$ & $+1.68$\,(1.36) & $+1.37$\,(1.49) & $+9.78$\,(1.94) & $-0.28$\,(1.11) & $+1.48$\,(1.26) & 47.2 \\
 & $K=4$ & $+0.62$\,(1.20) & $+0.66$\,(1.28) & $+3.91$\,(1.60) & $-0.29$\,(1.09) & $-0.94$\,(1.13) & 36.6 \\
\addlinespace
16000 (333) & exact MLE & $+0.08$\,(1.13) & $-0.05$\,(0.93) & $-0.05$\,(1.00) & $-0.08$\,(0.98) & $-0.05$\,(1.12) & -- \\
 & Euler ($K=0$) & $-1.38$\,(1.02) & $-0.89$\,(0.81) & $-4.02$\,(0.82) & $-9.02$\,(0.85) & $-17.11$\,(0.89) & -- \\
 & $K=2$ & $-0.78$\,(1.11) & $-0.38$\,(0.87) & $-0.70$\,(1.06) & $-0.14$\,(0.98) & $-0.15$\,(1.09) & 4.8 \\
 & $K=3$ & $+1.00$\,(1.28) & $+0.96$\,(1.10) & $+7.07$\,(1.54) & $+0.27$\,(1.00) & $+2.49$\,(1.18) & 5.8 \\
 & $K=4$ & $+0.24$\,(1.16) & $+0.23$\,(0.97) & $+1.56$\,(1.21) & $-0.12$\,(0.98) & $-0.12$\,(1.11) & 5.3 \\
\bottomrule\end{tabular}
\end{table}

\paragraph{Design 4: state-dependent volatility.}
Table~\ref{tab:nl} shows the same pattern in a model with state-dependent volatility, where $C_1\neq0$. The Euler estimator of the diffusion parameters is biased by $-1.2$ to $-1.9$ standardized units, and for $\beta_2$ the bias grows with $n$, as $\sqrt n\Delta_n\propto n^{0.1}$ predicts. Its Wald intervals for $\beta_2$ cover the true value in only 47--64\% of replications.

Adding the first correction ($K=1$) changes essentially nothing. This is what Remark~\ref{rem:parity} predicts: $\gamma_\beta=1$ for both $K=0$ and $K=1$, because $C_1$ is orthogonal to the leading diffusion score. The second-order expansion removes the bias, with standardized means between $-0.07$ and $0.18$, and restores the coverage of the diffusion parameters to 91--95\%. All remaining deviations are within two Monte Carlo standard errors.

The drift estimates are essentially the same for all three orders, with small finite-span biases of at most $0.23$. Without an exact benchmark, these cannot be separated from the expansion bias, which Remark~\ref{rem:parity} predicts to be of order $\sqrt n\Delta_n^{3/2}$, in standardized units, for all three orders. The correction terms remain small throughout: the average fraction of observations with $|R_i(\hat\theta)|>\frac14$ is at most 1.5\%.

Is the simulation scheme accurate enough for these conclusions? To check, we repeated the $n=16000$ experiment with 100 new replications, simulating each path twice from the same Brownian increments: once with 160 Euler sub-steps per observation interval and once with the 40 sub-steps used above. The paired differences in standardized errors (40 minus 160 sub-steps) are $0.026$ and $0.039$ for $\beta_1$ and $\beta_2$ under the Euler contrast, $0.028$ and $0.041$ under $K=2$, and at most $0.004$ for the drift parameters, with paired standard errors below $0.001$. These are below the rough estimate of $0.07$ given with the design and well below the Monte Carlo standard error of $0.10$, so the simulation scheme does not affect the comparison between expansion orders.

\begin{table}[!htbp]\centering\footnotesize
\caption{Standardized errors: Monte Carlo mean (standard deviation); the target is $0\,(1)$. Nonlinear model satisfying Assumptions~\ref{A:smooth}--\ref{A:id}, 100 replications; the Monte Carlo standard error of a mean is about $0.10$. The last column is as in Table~\ref{tab:ou3}.}\label{tab:nl}
\setlength{\tabcolsep}{3.5pt}\begin{tabular}{@{}llcccccr@{}}\toprule
$n$ ($n\Delta_n$) & estimator & $\alpha_1$ & $\alpha_2$ & $\alpha_3$ & $\beta_1$ & $\beta_2$ & $|R_i|>\frac14$ (\%)\\\midrule
\multicolumn{8}{@{}l}{\emph{Design 4: bounded nonlinear model, $\Delta_n=n^{-0.4}$}}\\
1000 (63) & Euler ($K=0$) & $+0.07$\,(0.97) & $-0.11$\,(0.98) & $-0.10$\,(1.05) & $-1.16$\,(0.98) & $-1.58$\,(0.91) & -- \\
 & $K=1$ & $+0.06$\,(0.97) & $-0.12$\,(0.97) & $-0.09$\,(1.06) & $-1.16$\,(0.98) & $-1.59$\,(0.91) & 0.6 \\
 & $K=2$ & $+0.11$\,(1.00) & $-0.12$\,(0.95) & $-0.02$\,(1.09) & $-0.07$\,(1.01) & $+0.01$\,(1.04) & 1.5 \\
\addlinespace
4000 (145) & Euler ($K=0$) & $-0.05$\,(1.01) & $-0.21$\,(0.96) & $-0.14$\,(0.94) & $-1.19$\,(0.96) & $-1.77$\,(0.89) & -- \\
 & $K=1$ & $-0.05$\,(1.01) & $-0.21$\,(0.96) & $-0.14$\,(0.93) & $-1.20$\,(0.96) & $-1.77$\,(0.89) & 0.3 \\
 & $K=2$ & $-0.03$\,(1.00) & $-0.23$\,(0.96) & $-0.13$\,(0.95) & $+0.03$\,(0.98) & $+0.06$\,(0.95) & 0.5 \\
\addlinespace
16000 (333) & Euler ($K=0$) & $-0.13$\,(1.07) & $+0.04$\,(0.98) & $-0.22$\,(1.02) & $-1.44$\,(1.03) & $-1.91$\,(1.05) & -- \\
 & $K=1$ & $-0.13$\,(1.07) & $+0.04$\,(0.98) & $-0.23$\,(1.02) & $-1.45$\,(1.03) & $-1.91$\,(1.05) & 0.1 \\
 & $K=2$ & $-0.13$\,(1.07) & $+0.04$\,(0.98) & $-0.22$\,(1.02) & $-0.05$\,(1.06) & $+0.18$\,(1.06) & 0.2 \\
\bottomrule\end{tabular}
\end{table}

\begin{table}[!htbp]\centering\footnotesize
\caption{Coverage (\%) of nominal 95\% Wald intervals based on the plug-in estimate of $\Gamma$, at the smallest and largest sample sizes. With 200 (100) replications, the Monte Carlo standard error of a coverage near 95\% is 1.5 (2.2) points.}\label{tab:cov}
\setlength{\tabcolsep}{5pt}\begin{tabular}{@{}llccccc@{}}\toprule
$n$ & estimator & \multicolumn{3}{c}{drift parameters} & \multicolumn{2}{c}{diffusion parameters}\\\midrule
\multicolumn{7}{@{}l}{\emph{Design 1}}\\
1000 & exact MLE & 94 & 94 & 92 & 96 & 90 \\
 & Euler ($K=0$) & 94 & 90 & 88 & 46 & 9 \\
 & $K=2$ & 94 & 95 & 92 & 96 & 90 \\
 & $K=3$ & 92 & 90 & 76 & 96 & 86 \\
 & $K=4$ & 94 & 92 & 88 & 96 & 90 \\
\addlinespace
16000 & exact MLE & 94 & 96 & 98 & 94 & 93 \\
 & Euler ($K=0$) & 95 & 96 & 98 & 75 & 25 \\
 & $K=2$ & 95 & 96 & 98 & 94 & 93 \\
 & $K=3$ & 94 & 96 & 98 & 94 & 93 \\
 & $K=4$ & 94 & 96 & 98 & 94 & 93 \\
\addlinespace
\multicolumn{7}{@{}l}{\emph{Design 2}}\\
1000 & exact MLE & 96 & 94 & 94 & 96 & 94 \\
 & Euler ($K=0$) & 96 & 94 & 96 & 88 & 67 \\
 & $K=2$ & 97 & 96 & 96 & 96 & 94 \\
 & $K=3$ & 96 & 94 & 94 & 95 & 94 \\
 & $K=4$ & 96 & 94 & 94 & 96 & 94 \\
\addlinespace
16000 & exact MLE & 96 & 96 & 96 & 95 & 97 \\
 & Euler ($K=0$) & 96 & 96 & 96 & 82 & 54 \\
 & $K=2$ & 96 & 96 & 96 & 95 & 97 \\
 & $K=3$ & 96 & 96 & 96 & 95 & 97 \\
 & $K=4$ & 96 & 96 & 96 & 95 & 97 \\
\addlinespace
\multicolumn{7}{@{}l}{\emph{Design 3}}\\
1000 & exact MLE & 90 & 93 & 86 & 91 & 88 \\
 & Euler ($K=0$) & 46 & 74 & 0 & 0 & 0 \\
 & $K=2$ & 91 & 88 & 56 & 71 & 4 \\
 & $K=3$ & 61 & 76 & 0 & 74 & 74 \\
 & $K=4$ & 81 & 82 & 26 & 80 & 8 \\
\addlinespace
16000 & exact MLE & 92 & 98 & 95 & 96 & 91 \\
 & Euler ($K=0$) & 67 & 86 & 0 & 0 & 0 \\
 & $K=2$ & 84 & 98 & 88 & 96 & 92 \\
 & $K=3$ & 78 & 82 & 0 & 94 & 34 \\
 & $K=4$ & 92 & 94 & 61 & 96 & 92 \\
\addlinespace
\multicolumn{7}{@{}l}{\emph{Design 4}}\\
1000 & Euler ($K=0$) & 97 & 95 & 93 & 75 & 64 \\
 & $K=1$ & 97 & 95 & 94 & 76 & 64 \\
 & $K=2$ & 97 & 95 & 94 & 94 & 95 \\
\addlinespace
16000 & Euler ($K=0$) & 91 & 96 & 94 & 71 & 47 \\
 & $K=1$ & 91 & 97 & 95 & 70 & 47 \\
 & $K=2$ & 91 & 96 & 95 & 94 & 91 \\
\bottomrule\end{tabular}
\end{table}

\FloatBarrier
\section{Conclusion}\label{sec:disc}

This paper develops the asymptotic theory of the exact MLE for multivariate diffusions observed at high frequency over a long span. For ergodic diffusions with a drift that may grow linearly and a volatility that is bounded away from zero and either bounded or proportional to the size of the state, the exact MLE is consistent, asymptotically normal and efficient whenever $\Delta_n$ decays polynomially, and, under an additional condition on ergodic averages and polynomial growth of $n\Delta_n$, its moments converge. This answers the question left open by \citet{YCW2019} and \citet{Li2013}, and it justifies Wald inference based on the exact likelihood for this class of models. The proof compares the exact likelihood with an explicit expansion of high order, and the same comparison gives conditions, in terms of approximation error, under which an approximate likelihood inherits the properties of the exact MLE. Together with the equivalence theorem of \citet{ChenWanYang2026}, this covers the expansion of \citet{Li2013}. For the delta-expansion estimator of fixed order, an explicit bias term shows that the condition $n\Delta_n^{K+1}\to0$ for efficiency cannot be weakened when the bias constant is nonzero, and identifies the parameters that the leading bias affects. The Monte Carlo evidence agrees with these predictions in size as well as in order, and shows that a higher order does not help when the sampling interval is coarse relative to the time scale of the model.

For future research, the most important is to cover square-root and other degenerate models, such as CIR and Heston, whose volatility vanishes at the boundary. This requires heat kernel bounds for degenerate diffusions and control of the boundary behavior. Volatility that grows in some directions but not in others, as in componentwise models, or at an intermediate rate, would require an anisotropic or non-proportional version of the rescaling in Section~\ref{sec:lin}. Further problems include:
\begin{itemize}[nosep]
\item drifts of superlinear growth, and time-inhomogeneous models, where no invariant law is available;
\item the exact MLE when $\Delta_n$ decays more slowly than any power of $n$ (Remark~\ref{rem:poly});
\item fixed-order theory for other expansions, such as the delta expansion with $\mu_0=0$ and log-density expansions, and a data-driven choice of $K$;
\item richer settings with irregular sampling, microstructure noise, jumps or latent factors.
\end{itemize}

\section*{Acknowledgment} The author would like to thank Nian Yang for sharing the code used in \citet{YCW2019}, which was very helpful for Section~\ref{sec:mc}.
Additionally, the author would like to thank all the co-authors of the paper \cite{{chau2024efficient}}, especially Thai Nguyen (Laval University, Canada) for many fruitful discussions while working on the paper, which led to insights that were essential in completing this paper.

\section*{Statement on AI use} In Section~\ref{sec:mc}, the author used the free version of Claude (Anthropic) to assist with code development in  discussions guided by the author's ideas. The author
checked the code, the results and the text, and takes full responsibility for the accuracy and
correctness of the mathematical arguments, the computational results and the content of the paper.

\appendix

\section{Preliminary lemmas}\label{app:prelim}

Throughout the appendices, $C$ and $r$ denote generic constants that may change from line to line and do not depend on $n$, $i$, $x$ or $\theta$. Recall $\jb{x}=(1+|x|^2)^{1/2}$, which satisfies $\jb{x'}\le\jb{x}+|x'-x|$ and $|\nabla\jb{x}|\le1$. A ``polynomial with weighted coefficients'' is a polynomial in $z$ whose coefficients $c(x;\theta)$ satisfy $\sum_{j\le3}|\partial^j_\theta c|+\sum_{j\le2}|\nabla_x\partial^j_\theta c|\le C\jb{x}^r$. The derivatives of order three in $\theta$ and the mixed derivatives are used only in Appendices~\ref{app:an} and \ref{app:moment}; for the coefficients of the $C_k$ they hold by the smoothness in Assumption~\ref{A:smooth}. We write $\Hspan(x,\beta)$ for the linear span of $\{H_h(\cdot;\nu(x;\beta)):|h|\ge1\}$. Each of its elements has mean zero under $N(0,\nu(x;\beta))$. Let $\zeta_i=(X_i-X_{i-1})/\eps$, so that $z_i(\alpha)=\zeta_i-\eps\mu(X_{i-1};\alpha)$, and let
\[
\mathcal L(\zeta,x;\theta)=\log\big(1+R(\zeta-\eps\mu(x;\alpha);x,\theta)\big)
\]
on $\{|R|\le\frac14\}$. Derivatives in $\theta$ are taken with the data $(\zeta,x)$ fixed. Compared with the case of bounded coefficients, the only change in this appendix is that every bound now carries a polynomial weight in the current state. Moment bounds make such weights harmless in averages, and Lemma~\ref{lem:max} shows that the largest state in the sample grows more slowly than any power of $n$.

The first lemma records the algebraic structure of the corrections $C_k$ that the consistency and bias arguments rely on: their polynomial form, their zero mean under the local Gaussian law, and the order in $\eps$ at which the drift parameters enter.

\begin{lemma}[structure of the corrections]\label{lem:struct}
Under Assumptions~\ref{A:smooth}--\ref{A:ell}:
\begin{enumerate}[label=(\alph*),nosep]
\item Each $C_k(\cdot;x,\theta)$ is a polynomial of degree at most $3k$ with weighted coefficients.
\item $C_k(\cdot;x,\theta)\in\Hspan(x,\beta)$ for every $k$. Consequently, $\int p^{(K,\Delta)}(x'|x;\theta)\,dx'=1$ for all $x$, $\theta$ and $\Delta$.
\item $C_1(z;x,\theta)$ does not depend on $\alpha$.
\item On $\{|R|\le\frac14\}$, with $z=\zeta-\eps\mu(x;\alpha)$, there are $G^{(1)},G^{(2)}\in\Hspan(x,\beta)$ with weighted coefficients, and functions $Q^{(1)}(z;x,\theta),Q^{(2)}(z;x,\theta)$ with $|Q^{(j)}|\le C\jb{x}^r(1+|z|^{6K+2})$, where $C$ does not depend on $(x,\theta)$, such that
\[
\partial_\alpha\mathcal L=\eps^2G^{(1)}(z;x,\theta)+\eps^3Q^{(1)}(z;x,\theta),\qquad
\partial^2_\alpha\mathcal L=\eps^2G^{(2)}(z;x,\theta)+\eps^3Q^{(2)}(z;x,\theta).
\]
Moreover, $|\partial_\alpha\mathcal L|\le C\eps^2\jb{x}^r(1+|z|^{3K})$, $|\mathcal L|+|\partial_\beta\mathcal L|\le C\eps\jb{x}^r(1+|z|^{3K})$, $|\partial^2_\beta\mathcal L|\le C\eps\jb{x}^r(1+|z|^{6K})$ and $|\partial_\alpha\partial_\beta\mathcal L|\le C\eps^2\jb{x}^r(1+|z|^{6K+2})$.
\end{enumerate}
\end{lemma}

Parts (a) and (b) show that the corrections are $O_\Pb(\eps)$ at the scale $1/n$. Part (d) gives the finer structure needed at the drift scale $1/(n\Delta)$.

\begin{proof}
(a) In \eqref{eq:delta}, $N\le2k$ gives $|h|=2N-k\le3k$. By Lemma~\ref{lem:w}(b) in Appendix~\ref{app:prop}, $|w_{N,h}(x;\theta)|\le C\jb{x}^N$. The coefficients $w_{N,h}(x;\theta)$ are polynomials in $\mu$, $\mu_0$, $\nu$ and their derivatives at $x$, so the same bound holds for their $\theta$-derivatives up to order three and for the mixed derivatives $\nabla_x\partial^j_\theta$, $j\le2$. Differentiating in $\alpha$ produces mixed derivatives $\partial^a_x\partial^b_\alpha\mu$ with $|a|\ge1$, which are bounded by Assumption~\ref{A:smooth}, so no extra powers of $\jb{x}$ arise beyond those already counted. The coefficients of $H_h(\cdot;\Sigma)$ are polynomials in the entries of $\Sigma^{-1}$, which are bounded and smooth in $\beta$ by Assumption~\ref{A:ell}.

(b) In \eqref{eq:delta}, $N\ge\lceil(k+1)/2\rceil$, so $|h|=2N-k\ge1$. The change of variables $x'\mapsto z$ gives $\int q\,H_h\,dx'=\int\phi(z;\nu)H_h(z;\nu)\,dz=(-1)^{|h|}\int\partial^h\phi(z;\nu)\,dz=0$.

(c) For $k=1$ the pairs $(N,|h|)$ in \eqref{eq:delta} are $(1,1)$ and $(2,3)$. With $\mu_0=\mu(x;\alpha)$, the recursion \eqref{eq:rec} gives $w_{1,e_j}(y)=\mu_j(y)-\mu_j(x)$, which vanishes at $y=x$. For $|h|=3$,
\[
w_{2,h}(y)=\mathcal Aw_{1,h}(y)+\sum_iA_iw_{1,h-e_i}(y)+\tfrac12\sum_{i,j}(\nu_{ij}(y)-\nu_{ij}(x))w_{1,h-e_i-e_j}(y).
\]
The first term is zero because $|h|>2$, and the last vanishes at $y=x$. In the middle term, $w_{1,h-e_i}$ has $|h-e_i|=2$, so it is a difference $\nu_{kl}(y)-\nu_{kl}(x)$, with factor $\frac12$ if $h-e_i=2e_k$ and factor $1$ if $h-e_i=e_k+e_l$ with $k\ne l$, since both ordered pairs contribute. The operator $A_if=(\mu_i(y)-\mu_i(x))f+\sum_j\nu_{ij}(y)\partial_jf$ maps it, at $y=x$, to $\sum_j\nu_{ij}(x)\partial_j\nu_{kl}(x)$, which depends on $\beta$ only. In tensor notation, $C_1=\frac14\sum_{a,b,i,j}\nu_{ab}\,\partial_a\nu_{ij}\,H_{bij}(z)$, where $H_{bij}=-\phi^{-1}\partial_b\partial_i\partial_j\phi$. In passing from the multi-index sum $\frac1{2!}\sum_{|h|=3}w_{2,h}H_h$ to the sum over ordered triples $(b,i,j)$, each multi-index $h$ appears $3!/h!$ times; this multiplicity and the pair factors above combine to the constant $\frac14$. We also checked the formula numerically against the recursion (Section~\ref{sec:mc}).

(d) As a function of $(\zeta,x,\theta)$, the parameter $\alpha$ enters $R=\sum_k\eps^kC_k(\zeta-\eps\mu(x;\alpha);x,\theta)$ in two ways. It enters through the argument, with $\partial_\alpha z=-\eps\,\partial_\alpha\mu$, and through the coefficients $w_{N,h}$. The Hermite polynomials depend on $\beta$ only. By (c), $\partial_\alpha(\eps C_1)=-\eps^2\nabla_zC_1\,\partial_\alpha\mu$. The generating function $\sum_hH_h(z;\Sigma)u^h/h!=\exp(u^\T\Sigma^{-1}z-\frac12u^\T\Sigma^{-1}u)$ yields $\partial_{z_j}H_h=\sum_k(\Sigma^{-1})_{jk}h_kH_{h-e_k}$. Since $C_1$ involves only $|h|=3$, we get $\nabla_zC_1\in\Hspan$. The explicit $\alpha$-derivative of $\eps^2C_2$ is $\eps^2$ times Hermite polynomials with $|h|\in\{2,4\}$; the term with $|h|=6$ is built from $\nu$ alone (Appendix~\ref{app:mc}) and does not depend on $\alpha$. The remaining contributions to $\partial_\alpha R$ carry at least $\eps^3$; these include the degree-one terms of $C_k$ for odd $k\ge3$. Writing $\partial_\alpha\mathcal L=\partial_\alpha R-R\,\partial_\alpha R/(1+R)$ gives the first display, with $G^{(1)}=\partial^{w}_\alpha C_2-\nabla_zC_1\,\partial_\alpha\mu$, where $C_2:=0$ if $K=1$. On $\{|R|\le\frac14\}$ we have $|1+R|\ge\frac34$. This gives the direct bound of degree $3K$ for $\partial_\alpha\mathcal L=\partial_\alpha R/(1+R)$, and the degree bound on $Q^{(1)}$ accounts for the product $R\,\partial_\alpha R/(1+R)$; similarly for $\partial^2_\alpha\mathcal L=[\partial^2_\alpha R\,(1+R)-\partial_\alpha R\,\partial_\alpha R^\T]/(1+R)^2$. Differentiating once more in $\alpha$ preserves the structure. The explicit derivative of $G^{(1)}$ stays in $\Hspan$, the term $-\eps\nabla_zG^{(1)}\partial_\alpha\mu$ is absorbed into $Q^{(2)}$, and so are the derivatives of the $\eps^3$ terms. The factors $\partial_\alpha\mu$ grow at most linearly and are absorbed into the weights. The remaining bounds follow from (a), because every term of $R$ carries at least one factor $\eps$, and every $\alpha$-derivative of $R$ carries at least $\eps^2$.
\end{proof}

The next lemma controls the largest state and the largest standardized increment in the sample. It is used to show that the truncation in \eqref{eq:contrast} is eventually inactive and, in Theorem~\ref{thm:mle}, to bound the worst-case relative error of the expansion.

\begin{lemma}[maximal states and increments]\label{lem:max}
Suppose Assumptions~\ref{A:smooth}--\ref{A:erg} hold, $\Delta_n\le Cn^{-a}$, and let $\delta>0$.
\begin{enumerate}[label=(\alph*),nosep]
\item For every $q\ge1$, $\Pb(\max_{0\le i\le n}|X_i|>n^\delta)\le C_qn^{1-q\delta}$.
\item There is $C$ such that, for $u>0$ and uniformly in $i$, on $\{|X_{i-1}|\le n^{\delta}\}$,
\[
\Pb\Big(\sup_{\alpha\in\Theta_\alpha}|z_i(\alpha)|\ge(1+C\Delta_n)u+C\eps n^{\delta}\,\Big|\,\cF_{i-1}\Big)\le2m\exp\Big(-\frac{u^2}{2m\lambda_{\max}}\Big).
\]
\item $\Pb\big(\max_{1\le i\le n}\sup_{\alpha}|z_i(\alpha)|^2\le2(1+\delta)\,m\lambda_{\max}\log n\big)\to1$, and $\Pb\big(\max_i\sup_\alpha|z_i(\alpha)|^2>c\log n\big)=O(n^{1-c/(2m\lambda_{\max})+\delta})+O(n^{-L})$ for every $c,L>0$.
\item $\max_i\sup_\theta|R_i(\theta)|\toP0$; in particular $\Pb(\Omega_n)\to1$.
\end{enumerate}
\end{lemma}

\begin{proof}
(a) follows from Markov's inequality, a union bound and Assumption~\ref{A:erg}(ii). For (b), let $x=X_{i-1}$ and $M_s=\int_{t_{i-1}}^s\sigma(X_u;\beta_0)\,dW_u$ for $s\in[t_{i-1},t_i]$. Since $|\mu(X_u)|\le|\mu(x)|+L|X_u-x|$ with $L$ the Lipschitz constant of $\mu$, Gronwall's lemma gives $\sup_s|X_s-x|\le e^{L\Delta}(\Delta|\mu(x)|+\sup_s|M_s|)$. Hence
\[
|\eps z_i-M_{t_i}|=\Big|\int_{t_{i-1}}^{t_i}(\mu(X_u)-\mu(x))\,du\Big|\le C\Delta\big(\Delta\jb{x}+\sup_s|M_s|\big),
\]
and $\sup_\alpha|z_i(\alpha)-z_i|\le\eps\sup_\alpha|\mu(x;\alpha)-\mu(x;\alpha_0)|\le C\eps\jb{x}$. Each coordinate $M^j$ is a continuous martingale with quadratic variation at most $\lambda_{\max}\Delta$, so the exponential inequality for continuous martingales with bounded quadratic variation \citep[Ch.~IV, Ex.~3.16]{RevuzYor1999} gives $\Pb(\sup_s|M_s|\ge\eps u\,|\,\cF_{i-1})\le2m\exp(-u^2/(2m\lambda_{\max}))$. This proves (b). For (c), take $\delta'<\min(\delta,a)/4$ and $u^2=2(1+\delta/2)m\lambda_{\max}\log n$ in (b). A union bound over $i$, together with (a) for $\delta'$, bounds the failure probability by $2mn^{-\delta/2}+o(1)$, because $\eps n^{\delta'}\to0$ and $(1+C\Delta)u\le\sqrt{2(1+\delta)m\lambda_{\max}\log n}$ for large $n$. The second claim follows in the same way, the term $O(n^{-L})$ coming from part (a) with $q$ large. For (d), Lemma~\ref{lem:struct}(a) gives $|R_i|\le\sum_k C\eps^k\jb{X_{i-1}}^r(1+|z_i(\alpha)|^{3k})$, which is at most $C\sum_k(\eps n^{r\delta'}(\log n)^{3/2})^k\to0$ on the events of (a) and (c).
\end{proof}

The bound in (c) is of the right order: for Gaussian increments, the maximum of $n$ squared standardized values grows like $2\log n$ times the largest variance. The factor $m$ comes from a crude union bound over coordinates. It enters the exponent of $n$, and hence the thresholds $(1+m\Lambda)/a$ in Theorems~\ref{thm:mle}(a) and \ref{thm:general}(b). A sharper chi-square-type bound would remove it, but the thresholds are not meant to be sharp, and we keep the simpler argument.

Over one step, the standardized increment is Gaussian up to an error of order $\eps$ in every $L^p$ norm. The following lemma makes this precise and replaces exact Gaussian computations throughout.

\begin{lemma}[one-step conditional moments]\label{lem:mom}
Under $\theta_0$ and Assumptions~\ref{A:smooth}--\ref{A:ell}, let $x=X_{i-1}$, $z_i=z_i(\alpha_0)$ and $\upsilon_i=(W_{t_i}-W_{t_{i-1}})/\eps$. Then:
\begin{enumerate}[label=(\alph*),nosep]
\item $z_i=\sigma(x;\beta_0)\upsilon_i+\rho_i$ with $\E[|\rho_i|^p\,|\,\cF_{i-1}]\le C_p\Delta^{p/2}\jb{x}^p$ for every $p\ge1$;
\item $|\E[z_i\,|\,\cF_{i-1}]|\le C\Delta^{3/2}\jb{x}$;
\item if $f$ satisfies $|f(z)-f(z')|\le c(1+|z|^r+|z'|^r)|z-z'|$, then
\[\Big|\E[f(z_i)|\cF_{i-1}]-\int f(z)\,\phi\big(z;\nu(x;\beta_0)\big)\,dz\Big|\le C_rc\,\eps\jb{x}^{r+1};\]
\item $\sup_{\alpha}|z_i(\alpha)-z_i|\le C\eps\jb{x}$, and $\E[\sup_\alpha|z_i(\alpha)|^p|\cF_{i-1}]+\E[\sup_{s\in[t_{i-1},t_i]}\jb{X_s}^p|\cF_{i-1}]\le C_p\jb{x}^p$.
\end{enumerate}
\end{lemma}

\begin{proof}
By Gronwall's lemma as in the proof of Lemma~\ref{lem:max} and the Burkholder--Davis--Gundy inequality, $\E[\sup_{s\in[t_{i-1},t_i]}|X_s-x|^p|\cF_{i-1}]\le C_p(\Delta^p\jb{x}^p+\Delta^{p/2})\le C_p\Delta^{p/2}\jb{x}^p$. Decompose
\[
z_i=\sigma(x)\upsilon_i+\eps^{-1}\int(\mu(X_s)-\mu(x))\,ds+\eps^{-1}\int(\sigma(X_s)-\sigma(x))\,dW_s.
\]
The second term is at most $\eps^{-1}\Delta\|\nabla\mu\|_\infty\sup_s|X_s-x|$, which has conditional $L^p$ norm $O(\Delta\jb{x})$. By BDG, the third term has conditional $L^p$ norm at most $C\eps^{-1}(\Delta\cdot\Delta\jb{x}^2)^{1/2}=C\eps\jb{x}$. This proves (a). For (b), $\E[z_i|\cF_{i-1}]=\eps^{-1}\int\E[\mu(X_s)-\mu(x)|\cF_{i-1}]\,ds$. By It\^o's formula, the integrand equals $\int_{t_{i-1}}^s\E[\mathcal A\mu(X_u)|\cF_{i-1}]\,du$, where $\mathcal A$ is the generator. Since $|\mathcal A\mu(y)|\le C\jb{y}$ under Assumption~\ref{A:smooth}, the integrand is at most $C\jb{x}(s-t_{i-1})$. For (c), the condition on $f$, H\"older's inequality and (a) bound $|\E[f(z_i)-f(\sigma(x)\upsilon_i)|\cF_{i-1}]|$ by $C_rc\,\eps\jb{x}^{r+1}$, and $\sigma(x)\upsilon_i\sim N(0,\nu(x;\beta_0))$ given $\cF_{i-1}$. Part (d) follows from the proof of Lemma~\ref{lem:max} and the moment bound above.
\end{proof}

Part (b) is the reason why the drift contrast in the consistency proof has a negligible centering: the Euler mean is correct to order $\Delta^2$, so the standardized increment has conditional mean $O(\Delta^{3/2})$ rather than $O(\eps)$.

The next three lemmas turn averages over observations into integrals against $\pi$, uniformly in the parameter. Lemma~\ref{lem:erg} is used for functions of the state, Lemma~\ref{lem:lln} in Step~1 of the consistency proof and in the Hessian bounds, and Lemma~\ref{lem:mart} at the drift scale.

\begin{lemma}[ergodic averages]\label{lem:erg}
Under Assumptions~\ref{A:smooth}--\ref{A:erg}, let $g:\R^m\to\R$ be $C^1$ with $|g(x)|+|\nabla g(x)|\le C\jb{x}^r$. If $\Delta_n\to0$ and $n\Delta_n\to\infty$, then $n^{-1}\sum_{i=1}^ng(X_{i-1})\toP\int g\,d\pi$, for every initial law allowed by Assumption~\ref{A:erg}.
\end{lemma}

\begin{proof}
Let $T_n=n\Delta_n$. Then
\[
\frac1n\sum_ig(X_{i-1})-\frac1{T_n}\int_0^{T_n}g(X_s)\,ds=\frac1{T_n}\sum_i\int_{t_{i-1}}^{t_i}\big(g(X_{i-1})-g(X_s)\big)\,ds.
\]
By the mean-value theorem, Assumption~\ref{A:erg}(ii) and Lemma~\ref{lem:mom}(d), its $L^1$ norm is $O(\eps)$. The continuous-time average converges almost surely by Assumption~\ref{A:erg}(i). Integrability of $g$ under $\pi$ follows from Assumption~\ref{A:erg}(ii) and Fatou's lemma.
\end{proof}

\begin{lemma}[uniform law of large numbers]\label{lem:lln}
Suppose Assumptions~\ref{A:smooth}--\ref{A:erg} hold, $\Delta_n\to0$ and $n\Delta_n\to\infty$. Let $\vartheta$ range over a compact set $\mathcal K\subset\R^q$, and let $f(x,z;\vartheta)$ be $C^1$ in $(x,\vartheta)$ and satisfy, for some $r$,
\[
|f|+|\partial_\vartheta f|+|\nabla_xf|\le C\jb{x}^r(1+|z|^r),\qquad|f(x,z;\vartheta)-f(x,z';\vartheta)|\le C\jb{x}^r(1+|z|^r+|z'|^r)|z-z'|.
\]
With $\bar f(x;\vartheta)=\int f(x,z;\vartheta)\phi(z;\nu(x;\beta_0))\,dz$,
\[
\sup_{\vartheta\in\mathcal K,\ \alpha'\in\Theta_\alpha}\Big|\frac1n\sum_{i=1}^nf\big(X_{i-1},z_i(\alpha');\vartheta\big)-\int\bar f(x;\vartheta)\,\pi(dx)\Big|\toP0 .
\]
\end{lemma}

\begin{proof}
By Lemma~\ref{lem:mom}(d) and the Lipschitz condition in $z$, replacing $z_i(\alpha')$ by $z_i=z_i(\alpha_0)$ changes the average by at most $C\eps n^{-1}\sum_i\jb{X_{i-1}}^{r+1}(1+\sup_\alpha|z_i(\alpha)|^r)=O_\Pb(\eps)$, uniformly. Fix $\vartheta$ and write $f_i=f(X_{i-1},z_i;\vartheta)$. Lemma~\ref{lem:mom}(c) gives $\E[f_i|\cF_{i-1}]=\bar f(X_{i-1};\vartheta)+O(\eps\jb{X_{i-1}}^{2r+1})$. The martingale $n^{-1}\sum_i(f_i-\E[f_i|\cF_{i-1}])$ has second moment at most $Cn^{-1}$ by Assumption~\ref{A:erg}(ii). The function $\bar f(\cdot;\vartheta)$ is $C^1$ with polynomially bounded gradient: differentiating under the integral, $\nabla_x\bar f=\int\nabla_xf\,\phi(z;\nu)\,dz+\int f\,\nabla_x\phi(z;\nu(x;\beta_0))\,dz$, where the first term is controlled by the bound on $\nabla_xf$ and the second because $\nu(\cdot;\beta_0)$ is smooth, bounded with bounded derivatives, and uniformly elliptic. When the lemma is applied to $\partial_\vartheta\bar f$ as well, as in Appendix~\ref{app:moment}, the same argument uses polynomial bounds on $\nabla_x\partial_\vartheta f$, which hold for all functions to which we apply it, since they are polynomials in $z$ with smooth coefficients. Lemma~\ref{lem:erg} then gives pointwise convergence. For uniformity, $|f_i(\vartheta)-f_i(\vartheta')|\le C|\vartheta-\vartheta'|\jb{X_{i-1}}^r(1+|z_i|^r)$, where $n^{-1}\sum_i\jb{X_{i-1}}^r(1+|z_i|^r)=O_\Pb(1)$, and the limit is Lipschitz in $\vartheta$. Convergence on a finite $\eta$-net of $\mathcal K$ therefore extends to all of $\mathcal K$.
\end{proof}

\begin{lemma}[uniform martingale bound]\label{lem:mart}
Let $g(x;\theta)$ be $\R^m$-valued and $C^1$ in $\theta$, with $|g|+|\partial_\theta g|\le C\jb{x}^r$, and let $M_n(\theta)=\sum_{i=1}^ng(X_{i-1};\theta)^\T(z_i-\E[z_i|\cF_{i-1}])$ with $z_i=z_i(\alpha_0)$. Under Assumptions~\ref{A:smooth}--\ref{A:erg}, $\E\sup_{\theta\in\Theta}|M_n(\theta)|^p\le C_pn^{p/2}$ for every $p\ge1$. Consequently,
\[
\sup_{\theta}\Big|\frac1{n\eps}\sum_{i=1}^ng(X_{i-1};\theta)^\T z_i\Big|=O_\Pb\big((n\Delta_n)^{-1/2}\big)+O_\Pb(\Delta_n).
\]
\end{lemma}

\begin{proof}
Let $p>p_\alpha+p_\beta$. Burkholder's inequality, Lemma~\ref{lem:mom}(a) and Assumption~\ref{A:erg}(ii) give $\E|M_n(\theta)|^p+\E|\partial_\theta M_n(\theta)|^p\le Cn^{p/2}$ uniformly in $\theta$. Since $\Theta$ is the closure of a bounded convex open set, the Sobolev embedding $W^{1,p}\hookrightarrow C^0$ \citep[Thm.~4.12]{AdamsFournier2003} gives $\E\sup_\theta|M_n|^p\le C\int_\Theta(\E|M_n|^p+\E|\partial_\theta M_n|^p)\,d\theta\le Cn^{p/2}$. Smaller $p$ follow by H\"older's inequality. The second statement then follows from Lemma~\ref{lem:mom}(b).
\end{proof}

Lemmas~\ref{lem:erg}--\ref{lem:mart} are standard in the high-frequency literature \citep[cf.][]{Kessler1997,Yoshida2011}. We include proofs because the uniformity over $(\alpha',\vartheta)$ in Lemma~\ref{lem:lln} is exactly what the treatment of the correction terms requires.

The last lemma of this appendix concerns the exact likelihood itself. It shows that the exact MLE is well defined, which the rest of the argument takes for granted.

\begin{lemma}[regularity of the exact density]\label{lem:reg}
Under Assumptions~\ref{A:smooth}--\ref{A:ell}, fix $\Delta\in(0,1]$ and choose for each $(x,\theta)$ the version of $p_\Delta(\cdot|x;\theta)$ that is continuous in $x'$. Then $p_\Delta(x'|x;\theta)$ is strictly positive and jointly continuous in $(x,x',\theta)$.
\end{lemma}

\begin{proof}
By the two-sided Gaussian bound and the forward gradient bound of \citet{MenozziPesceZhang2021} (see \eqref{eq:mpz} in Appendix~\ref{app:prop}), $p_\Delta$ is strictly positive and bounded, and $p_\Delta(\cdot|x;\theta)$ is Lipschitz with a constant that is uniform over $(x,\theta)$ in compact sets. Let $X^{x,\theta}$ denote the solution of \eqref{eq:sde} started at $x$ under $\theta$. Since the coefficients are Lipschitz in $x$ and smooth in $\theta$, with $|\mu(y;\theta)-\mu(y;\theta')|\le C|\theta-\theta'|\jb{y}$, the standard $L^2$ stability estimate gives $\E|X^{x,\theta}_\Delta-X^{x',\theta'}_\Delta|^2\to0$ as $(x',\theta')\to(x,\theta)$. Hence the law of $X^{x,\theta}_\Delta$ is weakly continuous in $(x,\theta)$. For a smooth mollifier $g_\eta$ with support in a ball of radius $\eta$, the map $(x,x',\theta)\mapsto\int g_\eta(x'-y)p_\Delta(y|x;\theta)\,dy$ is therefore continuous, and by the uniform Lipschitz bound it differs from $p_\Delta(x'|x;\theta)$ by at most $C\eta$, locally uniformly. A locally uniform limit of continuous functions is continuous.
\end{proof}

Only continuity in $\theta$ and strict positivity are used later: they make $\ell_n$ finite and continuous on the compact set $\Theta$, so that a measurable maximizer exists.

\section{Proof of Proposition~\ref{prop:unif}}\label{app:prop}

The proof follows the strategy of \citet[Theorems~3.2 and 3.3]{YCW2019}: write the error of an It\^o--Taylor kernel as a time integral by Dynkin's formula, then show that the integrand vanishes near the starting point to an order that compensates for the singularity of the Hermite kernels. It differs from theirs in three respects. Dynkin's formula is applied only to explicit functions, so no regularity of $p$ in the backward variables is needed. The Gaussian bounds are centered at the deterministic flow of the drift, which allows unbounded drift. And the bookkeeping of the vanishing orders is done uniformly on $\R^m$, which gives the polynomial weight and every order $K\ge1$.

Fix $\theta$, which we suppress, the starting point $x$, the endpoint $x'$, and $\Delta\in(0,1]$. Let $\mu_0=\mu_0(x;\theta)$ be the centering and $\nu_0=\nu(x;\beta)$. All constants below depend only on the bounds in Assumptions~\ref{A:smooth}--\ref{A:ell}, on $|\mu_0(x;\theta)|\le C\jb{x}$, and on the orders involved.

\paragraph{The It\^o--Taylor kernel.}
For $\tau>0$ and $y\in\R^m$ let
\[
F(\tau,y)=\tau^{-m/2}\phi(\varsigma;\nu_0),\qquad\varsigma=\frac{x'-y-\mu_0\tau}{\sqrt\tau}.
\]
This is the transition density from $y$ to $x'$ over time $\tau$ of a Brownian motion with constant drift $\mu_0$ and covariance $\nu_0$. Hence $(-\partial_\tau+\mathcal A_0)F=0$, where $\mathcal A_0=\mu_0\cdot\nabla+\frac12\nu_0:\nabla^2$ acts on $y$, and $\partial^h_yF=\tau^{-|h|/2}H_h(\varsigma;\nu_0)F$. Let $\mathcal A=\mu(y)\cdot\nabla+\frac12\nu(y):\nabla^2$ be the generator of \eqref{eq:sde}. Define functions $w_{N,h}(y)$, for $N\ge0$ and multi-indices $h$, by $w_{0,0}=1$, $w_{0,h}=0$ for $h\ne0$, and
\begin{equation}\label{eq:rec}
\begin{gathered}
w_{N+1,h}=\mathcal Aw_{N,h}+\sum_iA_iw_{N,h-e_i}+\frac12\sum_{i,j}\big(\nu_{ij}(y)-\nu_{0,ij}\big)w_{N,h-e_i-e_j},\\ A_if=(\mu_i(y)-\mu_{0,i})f+\sum_j\nu_{ij}(y)\partial_jf,
\end{gathered}
\end{equation}
with $w_{N,h}=0$ when $h$ has a negative entry. This is the recursion (18)--(20) of \citet{YCW2019} for time-homogeneous coefficients. Let $A_N=\sum_hw_{N,h}\,\partial^h_yF$. Because $\mathcal A_0$ has constant coefficients, $(-\partial_\tau+\mathcal A_0)\partial^h_yF=0$, and the product rule gives
\[
\begin{aligned}(-\partial_\tau+\mathcal A)(w\,\partial^h_yF)&=(\mathcal Aw)\,\partial^h_yF+\sum_i(\nu\nabla w)_i\,\partial^{h+e_i}_yF\\&\quad+w\Big[\sum_i(\mu_i-\mu_{0,i})\partial^{h+e_i}_yF+\tfrac12\sum_{i,j}(\nu_{ij}-\nu_{0,ij})\partial^{h+e_i+e_j}_yF\Big].\end{aligned}
\]
Summing over $h$ and collecting terms shows that $(-\partial_\tau+\mathcal A)A_N=A_{N+1}$. For $J\ge1$ let $P_J(\tau,y)=\sum_{N=0}^J\frac{\tau^N}{N!}A_N(\tau,y)$. The sum telescopes:
\begin{equation}\label{eq:psi}
(-\partial_\tau+\mathcal A)P_J=\psi_J:=\frac{\tau^J}{J!}A_{J+1}.
\end{equation}

The proof rests on two properties of the coefficient functions: they grow at most polynomially, and they vanish at the starting point to an order that increases with $|h|$. The next lemma makes both precise.

\begin{lemma}[coefficient functions]\label{lem:w}
For all $N$ and $h$:
\begin{enumerate}[label=(\alph*),nosep]
\item $w_{N,h}=0$ if $|h|>2N$, and $w_{N,0}=0$ if $N\ge1$;
\item $|\partial^aw_{N,h}(y)|\le C_{N,a}(\jb{x}+|y-x|)^N$ for every multi-index $a$;
\item $\partial^aw_{N,h}(x)=0$ whenever $|a|<2|h|-3N$. Consequently, $|\partial^aw_{N,h}(y)|\le C|y-x|^{(2|h|-3N-|a|)_+}(\jb{x}+|y-x|)^N$.
\end{enumerate}
\end{lemma}

\begin{proof}
All three claims follow by induction on $N$ from \eqref{eq:rec}. For (a), each step raises $|h|$ by at most two, and $\mathcal Aw_{N,0}=0$ when $w_{N,0}$ is constant. For (b), note that $|\mu(y)-\mu_0|\le C(\jb{x}+|y-x|)$, that $\mu$ appears in \eqref{eq:rec} at most once per step, either through $A_i$ or through $\mathcal A$, and that all derivatives of $\mu$ of order at least one, and $\nu$ with all its derivatives, are bounded. For (c), write $\mathrm{ord}(f)$ for the order to which $f$ vanishes at $x$. Then $\mathrm{ord}(fg)\ge\mathrm{ord}(f)+\mathrm{ord}(g)$, $\mathrm{ord}(\partial_jf)\ge\mathrm{ord}(f)-1$, and $\mathrm{ord}(\nu_{ij}-\nu_{0,ij})\ge1$. If $\mathrm{ord}(w_{N,h})\ge(2|h|-3N)_+$ for all $h$, the three terms of \eqref{eq:rec} vanish to order at least $(2|h|-3N)_+-2$, $(2|h|-2-3N)_+-1$ and $1+(2|h|-4-3N)_+$, respectively, and each of these is at least $2|h|-3(N+1)$. The bound follows from Taylor's theorem with integral remainder, applied to $\partial^aw_{N,h}$ at $x$, together with (b).
\end{proof}

Lemma~\ref{lem:w} holds for every constant centering $\mu_0$, which is what makes Proposition~\ref{prop:unif} independent of the centering. Part~(c) is Lemma~A.2 of \citet{YCW2019}; the growth bound (b) replaces their use of compactness.

\paragraph{Step 1: the kernel at the starting point.}
At $(\tau,y)=(\Delta,x)$ we have $\varsigma=z$ and $F=q$, so $P_J(\Delta,x)=q\sum_{N=0}^J\sum_h\frac1{N!}\Delta^{N-|h|/2}w_{N,h}(x)H_h(z;\nu_0)$. By Lemma~\ref{lem:w}(a,c), $w_{N,h}(x)\ne0$ only if $1\le|h|\le2N$ and $2|h|\le3N$, that is, with $k=2N-|h|$, only if $\lceil(k+1)/2\rceil\le N\le2k$. Grouping terms by $k$ shows that $P_J(\Delta,x)=q\big(1+\sum_{k=1}^{2J-1}\Delta^{k/2}C^{[J]}_k\big)$, where $C^{[J]}_k$ collects the terms of $C_k$ with $N\le J$. Take $J=2K$. Then $C^{[J]}_k=C_k$ for $k\le K$, and since $|w_{N,h}(x)|\le C\jb{x}^N$ by Lemma~\ref{lem:w}(b) and $\sup_z\phi(z;\nu_0)|H_h(z;\nu_0)|\le C$,
\begin{equation}\label{eq:step1}
\big|P_{2K}(\Delta,x)-p^{(K,\Delta)}(x'|x)\big|\le C\jb{x}^{2K}\Delta^{(K+1-m)/2}.
\end{equation}

\paragraph{Step 2: heat kernel bounds.}
Let $g_\lambda(t,v)=t^{-m/2}\exp(-\lambda|v|^2/t)$, and let $\chi_s(x)$ be the solution of $\dot\chi_s=\mu(\chi_s)$, $\chi_0=x$. By Gronwall's lemma, $|\chi_s(x)-x|\le Cs\jb{x}$ for $s\le1$. Under Assumptions~\ref{A:smooth}--\ref{A:ell}, the drift is Lipschitz with linear growth and $\sigma$ is bounded, uniformly elliptic and smooth with bounded derivatives. Theorem~1.2 of \citet{MenozziPesceZhang2021}, items~(i) and (iv), with time horizon $T=1$, then provides $\lambda>0$ and $C$, depending only on these bounds, such that for $s\in(0,1]$ the transition density is differentiable in the forward variable and
\begin{equation}\label{eq:mpz}
C^{-1}g_{1/\lambda'}\big(s,y-\chi_s(x)\big)\le p_s(y|x)\le C\,g_{2\lambda}\big(s,y-\chi_s(x)\big),\qquad|\nabla_yp_s(y|x)|\le Cs^{-1/2}g_{2\lambda}\big(s,y-\chi_s(x)\big).
\end{equation}
The first pair of bounds is their two-sided estimate (1.15), with $\lambda'=\lambda_0$; the second is their gradient estimate (1.19) in the forward variable. In what follows, $\alpha'$, $\beta'$ and the flow $\theta^{(1)}$ are in the notation of \citet{MenozziPesceZhang2021} and are unrelated to our parameters; we write $M_0,M_1,M_2$ for their constants $\kappa_0,\kappa_1,\kappa_2$, to avoid a clash with the mean-reversion matrix $\kappa$ of Section~\ref{sec:mc}. Their standing hypotheses are $(\mathbf H^\sigma_{\alpha'})$, namely $M_0^{-1}|\xi|^2\le\langle\sigma\sigma^\top(x)\xi,\xi\rangle\le M_0|\xi|^2$ and $|\sigma(x)-\sigma(y)|\le M_0|x-y|^{\alpha'}$ for some $\alpha'\in(0,1)$, and $(\mathbf H^b_{\beta'})$, namely $|\mu(0)|\le M_1$ and $|\mu(x)-\mu(y)|\le M_1(|x-y|^{\beta'}\vee|x-y|)$. Item~(iv) requires in addition $\beta'\in(0,1]$, $\|\nabla\sigma\|_\infty\le M_2$ and $|\nabla\sigma(x)-\nabla\sigma(y)|\le M_2|x-y|^{\alpha'}$. Under Assumptions~\ref{A:smooth}--\ref{A:ell} all of these hold with $\beta'=1$, any $\alpha'\in(0,1)$ and constants uniform in $\theta$: $\sigma$ and $\nabla\sigma$ are bounded with bounded derivatives, hence $\alpha'$-H\"older, and $|\mu(0;\alpha)|\le C$. Their bounds are centered at the flow $\theta^{(1)}_{s,0}(x)$ of a mollified drift. Since $\beta'=1>0$, their Remark~1.3 allows this flow to be replaced by the exact flow $\chi_s(x)$, which only changes the constants. Their constants depend on $(T,\alpha',\beta',M_0,M_1,m)$ and, for (iv), on $M_2$. For coefficients that are bounded and smooth with bounded derivatives, as for the localized process of Appendix~\ref{app:lin}, Gaussian bounds for the first derivatives of the transition density in both variables also follow from the classical parametrix construction \citep{Friedman1964}; the flow-centered estimates of \citet{MenozziPesceZhang2021} are what allow a drift of linear growth. In particular, $p_\Delta(\cdot|x)$ is continuous. Decreasing $\lambda$ if necessary, Assumption~\ref{A:ell} also gives $|\partial^h_yF(\tau,y)|\le C\tau^{-|h|/2}g_{2\lambda}(\tau,x'-\mu_0\tau-y)$ for $|h|\le2J+2$. We use two elementary facts: for $t\le1$ and $k,r\ge0$,
\begin{equation}\label{eq:elem}
|v|^k(1+|v|)^rg_{2\lambda}(t,v)\le C\,t^{k/2}g_\lambda(t,v),\qquad\int g_\lambda(\tau,a-y)\,g_\lambda(s,y-b)\,dy=(\pi/\lambda)^{m/2}g_\lambda(\tau+s,a-b).
\end{equation}
The factor $s^{(k-|a|)_+/2}$ below arises as follows: by Lemma~\ref{lem:w}(c), $|\partial^aw_{N,h}(y)|$ carries a factor $|y-x|^{(k-|a|)_+}$, and $|y-x|^j\le C(|y-\chi_s(x)|^j+(s\jb{x})^j)$; the first part is absorbed into the Gaussian by \eqref{eq:elem}, at the cost of $s^{j/2}$. Combining \eqref{eq:mpz} and \eqref{eq:elem} with Lemma~\ref{lem:w}(c), and using $|y-x|\le|y-\chi_s(x)|+Cs\jb{x}$ and $s^j\le s^{j/2}$ for $j\ge0$ and $s\le1$, we get for $|a|\le1$, $N=J+1$ and $k=(2|h|-3N)_+$,
\begin{equation}\label{eq:wp}\begin{gathered}
|\partial^aw_{N,h}(y)|\,p_s(y|x)\le C\jb{x}^{r}s^{(k-|a|)_+/2}g_\lambda\big(s,y-\chi_s(x)\big),\\ |w_{N,h}(y)|\,|\nabla_yp_s(y|x)|\le C\jb{x}^{r}s^{(k-1)/2}g_\lambda\big(s,y-\chi_s(x)\big),
\end{gathered}\end{equation}
the second bound for $k\ge1$.

\paragraph{Step 3: Duhamel representation.}
We claim that, for every $x'$,
\begin{equation}\label{eq:duhamel}
p_\Delta(x'|x)-P_J(\Delta,x)=\int_0^\Delta I(s)\,ds,\qquad I(s)=\int\psi_J(\Delta-s,y)\,p_s(y|x)\,dy,
\end{equation}
where, in the terms of $\psi_J$ with $|h|=2J+2$, the integral over $y$ is understood after one integration by parts. Precisely, for such $h$ we pick $i$ with $h_i\ge1$ and use
\[
\int w_{J+1,h}(y)\,\partial^h_yF\,p_s(y|x)\,dy=-\int\partial^{h-e_i}_yF\,\partial_{y_i}\big(w_{J+1,h}(y)\,p_s(y|x)\big)\,dy,
\]
which holds for $s\in(0,\Delta)$ because $F(\Delta-s,\cdot)$ and its derivatives decay like a Gaussian, while $w_{J+1,h}$ grows polynomially and $p_s$, $\nabla p_s$ are bounded by \eqref{eq:mpz}. Step~4 shows that $\int_0^\Delta|I(s)|\,ds<\infty$, with this convention.

To prove \eqref{eq:duhamel}, fix $\eta\in(0,\Delta)$ and apply It\^o's formula to $P_J(\Delta-s,X_s)$ for $s\in[0,\Delta-\eta]$, with $X_0=x$. On this time interval $P_J$ is smooth, and it has polynomially bounded derivatives, so the stochastic integral is a martingale. By \eqref{eq:psi},
\begin{equation}\label{eq:dynkin}
\int P_J(\eta,y)\,p_{\Delta-\eta}(y|x)\,dy-P_J(\Delta,x)=\int_0^{\Delta-\eta}I(s)\,ds .
\end{equation}
Multiply by $g\in C^\infty_c(\R^m)$, viewed as a function of $x'$, and integrate over $x'$. Since $F$ depends on $(x',y)$ through $x'-y$, $\partial^h_yF=(-1)^{|h|}\partial^h_{x'}F$, and integration by parts in $x'$ gives
\[
\int g(x')\,A_N(\eta,y)\,dx'=\sum_hw_{N,h}(y)\int\partial^hg(x')\,F(\eta,y)\,dx'.
\]
For $N\ge1$ this is bounded by $C(\jb{x}+|y-x|)^N$, so after multiplication by $\eta^N$ and integration against $p_{\Delta-\eta}(y|x)\,dy$ it tends to zero as $\eta\to0$. For $N=0$, $\int g(x')F(\eta,y)\,dx'=\E g(y+\mu_0\eta+\sqrt\eta\,\nu_0^{1/2}Z)$ with $Z\sim N(0,I)$, which converges to $g(y)$ uniformly in $y$. Hence the first term of \eqref{eq:dynkin}, integrated against $g$, tends to $\E g(X_\Delta)=\int g\,p_\Delta(\cdot|x)$. The right-hand side converges by Step~4 and dominated convergence. Therefore \eqref{eq:duhamel} holds after integration against every $g$. Both sides of \eqref{eq:duhamel} are continuous in $x'$: the left by Step~2, the right by the bounds of Step~4 and dominated convergence. So \eqref{eq:duhamel} holds for every $x'$.

\paragraph{Step 4: bounding the integral.}
Let $\tau=\Delta-s$. Consider first a term of $\psi_J$ with $|h|\le2J+1$, and let $k=(2|h|-3(J+1))_+$. By Step~2, \eqref{eq:wp} with $a=0$, and the convolution identity in \eqref{eq:elem},
\[
\begin{aligned}\frac{\tau^J}{J!}\int|w_{J+1,h}(y)|\,|\partial^h_yF(\tau,y)|\,p_s(y|x)\,dy&\le C\jb{x}^r\tau^{J-|h|/2}s^{k/2}\\&\qquad\times\int g_\lambda(\tau,x'-\mu_0\tau-y)\,g_\lambda(s,y-\chi_s(x))\,dy\\&\le C\jb{x}^r\Delta^{-m/2}\tau^{J-|h|/2}s^{k/2}.\end{aligned}
\]
The last step uses the convolution identity with $\tau+s=\Delta$ and $g_\lambda(\Delta,v)\le\Delta^{-m/2}$.
Since $J-|h|/2\ge-\frac12$, integration over $s\in(0,\Delta)$ gives a beta integral, and the term contributes at most $C\jb{x}^r\Delta^{-m/2}\Delta^{E_h}$ with $E_h=J+1-|h|/2+k/2$. For $J=2K$, if $|h|\le3K+1$ then $E_h\ge2K+1-(3K+1)/2=(K+1)/2$, and if $|h|\ge3K+2$ then $k=2|h|-6K-3$ and $E_h=|h|/2-K-\frac12\ge(K+1)/2$. For the terms with $|h|=2J+2$, after the integration by parts of Step~3, both bounds in \eqref{eq:wp}, with $|a|=1$ and $k=J+1$, give the integrand bound $C\jb{x}^r\Delta^{-m/2}\tau^{J-(2J+1)/2}s^{J/2}=C\jb{x}^r\Delta^{-m/2}\tau^{-1/2}s^{K}$. Integration over $s$ gives $\Delta^{K+1/2}\le\Delta^{(K+1)/2}$. Summing over the finitely many $h$,
\[
\big|p_\Delta(x'|x)-P_{2K}(\Delta,x)\big|\le C\jb{x}^r\Delta^{(K+1-m)/2}.
\]
Together with \eqref{eq:step1}, this proves Proposition~\ref{prop:unif}, with $r_K$ the larger of the two exponents. \qed

\medskip
The proof shows where each assumption enters. Uniform ellipticity and bounded, smooth volatility give the Gaussian bounds for the Hermite kernels and, through \citet{MenozziPesceZhang2021}, for $p_s$ and its forward gradient. The Lipschitz drift with linear growth keeps those bounds centered at a flow that moves by $O(s\jb{x})$, and it makes the coefficient functions grow polynomially. The vanishing orders of Lemma~\ref{lem:w}(c) do the rest.

\section{Proof of Theorem~\ref{thm:cons}}\label{app:cons}

The proof treats the two time scales separately. Step~1 identifies $\beta$ from the contrast normalized by $n$, and Step~2 identifies $\alpha$ from the contrast normalized by $n\Delta_n$, using the consistency of $\hat\beta_n$. For $K=0$ there is nothing to truncate, and $\Omega_n$ is the whole space. For $K\ge1$, Lemma~\ref{lem:max}(d) gives $\Pb(\Omega_n)\to1$. We work on $\Omega_n$, where $\ell_n^{(K)}=\sum_i[\log q_i+\mathcal L_i]$. Write $\nu_i(\beta)=\nu(X_{i-1};\beta)$, $\delta_i(\alpha)=\mu(X_{i-1};\alpha)-\mu(X_{i-1};\alpha_0)$ and $z_i=z_i(\alpha_0)$, so that $z_i(\alpha)=z_i-\eps\delta_i(\alpha)$ and $|\delta_i(\alpha)|\le C\jb{X_{i-1}}$.

\paragraph{Step 1: diffusion parameters.}
Let $Y_n(\alpha,\beta)=n^{-1}[\ell^{(K)}_n(\alpha,\beta)-\ell^{(K)}_n(\alpha,\beta_0)]$. By Lemmas~\ref{lem:struct}(d) and \ref{lem:mom}, the corrections contribute at most $2n^{-1}\sum_i\sup_\theta|\mathcal L_i|=O_\Pb(\eps)$, uniformly in $\theta$. The Gaussian part is
\[
-\frac1{2n}\sum_i\Big[\log\det\big(\nu_i(\beta)\nu_i(\beta_0)^{-1}\big)+z_i(\alpha)^\T\big(\nu_i(\beta)^{-1}-\nu_i(\beta_0)^{-1}\big)z_i(\alpha)\Big].
\]
By Lemma~\ref{lem:lln}, it converges uniformly in $(\alpha,\beta)$ to
\[
U_\beta(\beta)=-\frac12\int\big[\tr A(x)-\log\det A(x)-m\big]\,\pi(dx),\qquad A(x)=\nu(x;\beta)^{-1}\nu(x;\beta_0).
\]
The matrix $A(x)$ is similar to a symmetric positive definite matrix with eigenvalues $\lambda_j>0$, and $\tr A-\log\det A-m=\sum_j(\lambda_j-1-\log\lambda_j)\ge0$, with equality if and only if $A=I$. Hence $U_\beta\le0$ is continuous, and by Assumption~\ref{A:id}(a) it vanishes only at $\beta_0$. Since $\hat\theta^{(K)}_n$ maximizes $\ell^{(K)}_n$, $Y_n(\hat\alpha_n,\hat\beta_n)\ge0$. For $\eta>0$ let $u_\eta=\sup_{|\beta-\beta_0|\ge\eta}U_\beta(\beta)<0$. Then $\Pb(|\hat\beta_n-\beta_0|\ge\eta)\le\Pb(\sup_\theta|Y_n-U_\beta|\ge-u_\eta)+\Pb(\Omega_n^c)\to0$.

\paragraph{Step 2: drift parameters.}
Let $V_n(\alpha)=(n\Delta)^{-1}[\ell^{(K)}_n(\alpha,\hat\beta_n)-\ell^{(K)}_n(\alpha_0,\hat\beta_n)]$. The log-determinants cancel, and the Gaussian part equals
\[
\frac1{n\eps}\sum_i\delta_i(\alpha)^\T\nu_i(\hat\beta_n)^{-1}z_i-\frac1{2n}\sum_i\delta_i(\alpha)^\T\nu_i(\hat\beta_n)^{-1}\delta_i(\alpha).
\]
By Lemma~\ref{lem:mart} with $g=\nu(\cdot;\beta)^{-1}\delta(\cdot;\alpha)$, the first sum is $o_\Pb(1)$ uniformly in $(\alpha,\beta)$. By Lemma~\ref{lem:lln}, uniform continuity of the limit in $\beta$, and $\hat\beta_n\toP\beta_0$, the second converges uniformly in $\alpha$ to
\[
U_\alpha(\alpha)=-\frac12\int\delta(x;\alpha)^\T\nu(x;\beta_0)^{-1}\delta(x;\alpha)\,\pi(dx),\qquad\delta(x;\alpha)=\mu(x;\alpha)-\mu(x;\alpha_0).
\]
For the corrections, the mean-value theorem and Lemma~\ref{lem:struct}(d) give
\[
\frac1{n\Delta}\sum_i\big[\mathcal L_i(\alpha,\beta)-\mathcal L_i(\alpha_0,\beta)\big]
=(\alpha-\alpha_0)^\T\frac1n\sum_i\int_0^1G^{(1)}\big(z_i(\alpha_s);X_{i-1},\alpha_s,\beta\big)\,ds+O_\Pb(\eps),
\]
with $\alpha_s=\alpha_0+s(\alpha-\alpha_0)$, uniformly in $(\alpha,\beta)$. Lemma~\ref{lem:lln}, applied with $\vartheta=(s,\alpha,\beta)\in[0,1]\times\Theta$ and $\alpha'=\alpha_s$, shows that this converges uniformly to $(\alpha-\alpha_0)^\T\int_0^1\gamma(\alpha_s,\beta)\,ds$, where
\[
\gamma(\alpha,\beta)=\int\E_{Z\sim N(0,\nu(x;\beta_0))}\big[G^{(1)}(Z;x,\alpha,\beta)\big]\,\pi(dx).
\]
Because $G^{(1)}(\cdot;x,\alpha,\beta)\in\Hspan(x,\beta)$, we have $\gamma(\alpha,\beta_0)=0$ for every $\alpha$. By uniform continuity of $\gamma$ and $\hat\beta_n\toP\beta_0$, the corrections are $o_\Pb(1)$ uniformly in $\alpha$ at $\beta=\hat\beta_n$. Hence $\sup_\alpha|V_n(\alpha)-U_\alpha(\alpha)|\toP0$. Since $U_\alpha\le0$ is continuous with a unique zero at $\alpha_0$ by Assumption~\ref{A:id}(b), and $V_n(\hat\alpha_n)\ge0$, the argument of Step~1 gives $\hat\alpha_n\toP\alpha_0$. \qed

\begin{remark}[near-maximizers]\label{rem:near}
The proof uses the maximizing property only through $Y_n(\hat\theta^{(K)}_n)\ge0$ and $V_n(\hat\alpha_n)\ge0$. It therefore applies to any $\tilde\theta_n=(\tilde\alpha_n,\tilde\beta_n)$ with $\ell^{(K)}_n(\tilde\theta_n)\ge\sup_\Theta\ell_n^{(K)}-\eta_n$, where $\eta_n=o_\Pb(n\Delta_n)$: then $Y_n(\tilde\theta_n)\ge-\eta_n/n=-o_\Pb(1)$ and, since $(\alpha_0,\tilde\beta_n)\in\Theta$, $V_n(\tilde\alpha_n)\ge-\eta_n/(n\Delta_n)=-o_\Pb(1)$, which is all the argument needs. We use this in the proofs of Theorems~\ref{thm:mle} and \ref{thm:general}.
\end{remark}

\section{Proof of Theorem~\ref{thm:an} and Corollary~\ref{cor:regimes}}\label{app:an}

The proof has three parts: an expansion of the conditional bias of the score, a central limit theorem for its martingale part, and uniform convergence of the Hessian. Only the first part involves the expansion error, and it is where the bias vector $\bar B^{(K)}$ appears. Fix $c\ge8m\lambda_{\max}$ and set $r_n=\sqrt{c\log n}$. Fix also $\delta_0\in(0,a/2)$ with $\delta_0<a/(4r)$, where $r$ is the largest of the polynomial state-weight exponents that appear below: that of Proposition~\ref{prop:unif} at order $K+1$, and those of Lemma~\ref{lem:struct} for the score and for the tail term in Lemma~\ref{lem:bias}. Then $\jb{X_{i-1}}^r\le Cn^{a/4}$ on $\mathcal G_i$, and $\delta_0$ does not depend on $n$ or on the moment orders used in Appendix~\ref{app:moment}. Let $\mathcal G_i=\{|X_{i-1}|\le n^{\delta_0}\}$, which is $\cF_{i-1}$-measurable, and let $B_n=\{x':|z|\le r_n\}$, with $z$ evaluated at $\theta_0$. By Lemma~\ref{lem:max}(b), $\Pb(|z_i|>r_n\,|\,\cF_{i-1})\le Cn^{-3}$ on $\mathcal G_i$, and $\Pb(\max_i|z_i|\le r_n,\ \max_i|X_{i-1}|\le n^{\delta_0})\to1$. On $B_n$ and $\mathcal G_i$, Lemma~\ref{lem:struct}(a) gives $|R|\le C\eps n^{r\delta_0}r_n^{3K}\to0$, so $p^{(K,\Delta)}\ge\frac34q$ there for large $n$.

\begin{lemma}[bias of the score]\label{lem:bias}
Suppose the assumptions of Theorem~\ref{thm:an} hold, and let
\[
s_i=\partial_\theta\log p^{(K,\Delta)}(X_i|X_{i-1};\theta_0),\qquad\tilde s_i=s_i\,\mathbf 1\{|z_i|\le r_n\}\,\mathbf 1_{\mathcal G_i}.
\]
Then, uniformly in $i$, on $\mathcal G_i$,
\[
\Big|\eps^{-1}\E[\tilde s_{i,\alpha}|\cF_{i-1}]-\Delta^{\frac{K+1}2}b_\alpha(X_{i-1})\Big|+\Big|\E[\tilde s_{i,\beta}|\cF_{i-1}]-\Delta^{\frac{K+1}2}b_\beta(X_{i-1})\Big|\le C\Delta^{\frac{K+2}2}n^{r\delta_0}r_n^{r}+Cn^{-2}.
\]
\end{lemma}
Since $\mathcal G_i\in\cF_{i-1}$ and $\tilde s_i$ vanishes off $\mathcal G_i$, the bound holds on the whole probability space once $b_\alpha$, $b_\beta$ and the right-hand side are multiplied by $\mathbf 1_{\mathcal G_i}$.

\begin{proof}
Let $x=X_{i-1}$. By Lemma~\ref{lem:struct}(b), $\int\partial_\theta p^{(K,\Delta)}\,dx'=\partial_\theta1=0$. Hence
\[
\E[\tilde s_i|\cF_{i-1}]=\int_{B_n}\partial_\theta\log p^{(K,\Delta)}\,(p_\Delta-p^{(K,\Delta)})\,dx'-\int_{B_n^c}\partial_\theta p^{(K,\Delta)}\,dx'.
\]
We split $p_\Delta-p^{(K,\Delta)}=q\,\Delta^{(K+1)/2}C_{K+1}+(p_\Delta-p^{(K+1,\Delta)})$. By Proposition~\ref{prop:unif} at order $K+1$, the second piece is at most $C\jb{x}^r\Delta^{(K+2-m)/2}$. On $B_n$, Lemma~\ref{lem:struct} gives $|\partial_\beta\log p^{(K,\Delta)}|+\eps^{-1}|\partial_\alpha\log p^{(K,\Delta)}|\le C\jb{x}^r(1+|z|^{3K})$, and the Lebesgue measure of $B_n$ is $C\Delta^{m/2}r_n^m$. The contribution of the second piece is therefore within the stated bound. For the first piece, change variables to $z$, so that $q\,dx'=\phi(z;\nu)\,dz$ with $\nu=\nu(x;\beta_0)$. Since $|R|\le\frac14$ on $B_n$ and $\mathcal G_i$ (see the beginning of this appendix), Lemma~\ref{lem:struct}(d) applies there and gives
\[\begin{gathered}
\partial_\beta\log p^{(K,\Delta)}=S_\beta(z;x)+O\big(\eps\jb{x}^r(1+|z|^{3K})\big),\\ \eps^{-1}\partial_\alpha\log p^{(K,\Delta)}=(\partial_\alpha\mu)^\T\nu^{-1}z+O\big(\eps\jb{x}^r(1+|z|^{6K+2})\big),
\end{gathered}\]
where $S_{\beta_j}(z;x)=\frac12[z^\T\nu^{-1}\partial_{\beta_j}\nu\,\nu^{-1}z-\tr(\nu^{-1}\partial_{\beta_j}\nu)]$ is the Euler score. Integrating against $C_{K+1}\phi$, whose coefficients are bounded by $C\jb{x}^r$, over $\{|z|\le r_n\}$ gives $\Delta^{(K+1)/2}$ times $b(x)+O(\eps\jb{x}^r)$, up to Gaussian tails of order $n^{-2}$. Finally, $|\partial_\theta p^{(K,\Delta)}|\le q\,C\jb{x}^r(1+|z|^{3K+2})$, with an extra factor $\eps$ for the $\alpha$-component, so the integral over $B_n^c$ is at most $C\jb{x}^r\int_{|z|>r_n}(1+|z|^{3K+2})e^{-|z|^2/(2\lambda_{\max})}\,dz\le Cn^{-2}$ on $\mathcal G_i$, by the choice of $c$ and $\delta_0$.
\end{proof}

The lemma shows that the score of the delta-expansion contrast is a martingale up to a drift whose leading term is explicit and of order $\Delta^{(K+1)/2}$ per observation. The next lemma handles the martingale part. It also shows that this part does not depend on $K$, up to $o_\Pb(1)$, which is what makes the different estimators comparable.

\begin{lemma}[central limit theorem for the score]\label{lem:clt}
Under the assumptions of Theorem~\ref{thm:an}, $S^{(K)}_n:=D_n^{-1}\sum_{i=1}^n(\tilde s_i-\E[\tilde s_i|\cF_{i-1}])\toD N(0,\Gamma)$. Moreover, with $S_n:=S^{(0)}_n$ built from the Euler score, $S^{(K)}_n-S_n\toP0$ for every fixed $K$. In the rest of this appendix, $K$ is fixed and we write $S_n$ also for $S^{(K)}_n$, which is harmless by the second statement.
\end{lemma}

\begin{proof}
Write $\tilde s_{i,\alpha}=\eps a_i$ and $\tilde s_{i,\beta}=b_i$, so that $D_n^{-1}\sum_i\tilde s_i=n^{-1/2}\sum_i(a_i,b_i)$. By Lemma~\ref{lem:struct}, on $\{|z_i|\le r_n\}\cap\mathcal G_i$,
\[
a_i=(\partial_\alpha\mu)^\T\nu^{-1}z_i+\eps P_{i},\qquad
b_{i,j}=\tfrac12\big[z_i^\T\nu^{-1}\partial_{\beta_j}\nu\,\nu^{-1}z_i-\tr(\nu^{-1}\partial_{\beta_j}\nu)\big]+\eps P'_{i,j},
\]
with coefficients evaluated at $(X_{i-1};\theta_0)$ and polynomials $P_i$, $P'_{i,j}$ with weighted coefficients. By the Cauchy--Schwarz inequality and the conditional tail bound, removing the indicators from these polynomial expressions changes conditional moments of order up to four by $O(\jb{X_{i-1}}^rn^{-3/2})$. Lemma~\ref{lem:mom}(c) and Gaussian moment identities then give:
\begin{itemize}[nosep]
\item $\E[a_ia_i^\T|\cF_{i-1}]=(\partial_\alpha\mu)^\T\nu^{-1}\partial_\alpha\mu+O(\eps\jb{X_{i-1}}^r)$;
\item $\E[b_{i,j}b_{i,k}|\cF_{i-1}]=\frac12\tr(\nu^{-1}\partial_{\beta_j}\nu\,\nu^{-1}\partial_{\beta_k}\nu)+O(\eps\jb{X_{i-1}}^r)$, using $\mathrm{Cov}(\frac12Z^\T BZ,\frac12Z^\T B'Z)=\frac12\tr(B\nu B'\nu)$ for $Z\sim N(0,\nu)$;
\item $\E[a_ib_i^\T|\cF_{i-1}]=O(\eps\jb{X_{i-1}}^r)$, since odd Gaussian moments vanish.
\end{itemize}
The products of conditional means are $O(\Delta\jb{X_{i-1}}^r)$ by Lemma~\ref{lem:bias}, so by Lemma~\ref{lem:erg} and Assumption~\ref{A:erg}(ii) the conditional covariance of $S_n$ converges to $\Gamma$. The conditional Lyapunov condition $n^{-2}\sum_i\E[|(a_i,b_i)|^4|\cF_{i-1}]\le Cn^{-2}\sum_i\jb{X_{i-1}}^r=O_\Pb(n^{-1})$ follows from Lemma~\ref{lem:mom}. The martingale central limit theorem \citep[Cor.~3.1]{HallHeyde1980} gives the first claim. For the second, the increments of $S_n-S'_n$ are $n^{-1/2}\eps(P_i,P'_i)$ minus their conditional means, so the conditional variance of $S_n-S'_n$ is at most $C\Delta n^{-1}\sum_i\jb{X_{i-1}}^r=O_\Pb(\Delta)$.
\end{proof}

It remains to show that the normalized Hessian converges, so that the score equation can be linearized.

\begin{lemma}[Hessian]\label{lem:hess}
Under the assumptions of Theorem~\ref{thm:cons}, for any $K\ge0$ (for $K=0$, $R\equiv0$ and $\mathcal L\equiv0$), let $\mathbb H_n(\theta)=D_n^{-1}\partial_\theta^2\ell_n^{(K)}(\theta)D_n^{-1}$ on $\Omega_n$. There is a continuous function $\mathbb H(\theta)$ with $\mathbb H(\theta_0)=-\Gamma$ such that $\sup_{\theta\in\Theta}|\mathbb H_n(\theta)-\mathbb H(\theta)|\toP0$. In particular, $\sup_{|\theta-\theta_0|\le\eta_n}|\mathbb H_n(\theta)+\Gamma|\toP0$ for every sequence $\eta_n\to0$.
\end{lemma}

\begin{proof}
We use $z_i(\alpha)=z_i-\eps\delta_i(\alpha)$ and the derivatives of $\log q$, and treat each block separately.

\emph{Drift block.} From $\partial^2_\alpha\log q_i=-\Delta(\partial_\alpha\mu)^\T\nu^{-1}\partial_\alpha\mu+\eps(\partial^2_\alpha\mu)^\T\nu^{-1}z_i(\alpha)$ and Lemma~\ref{lem:struct}(d), $(n\Delta)^{-1}\partial^2_\alpha\ell^{(K)}_n$ equals
\[
-\frac1n\sum(\partial_\alpha\mu)^\T\nu^{-1}\partial_\alpha\mu+\frac1{n\eps}\sum(\partial^2_\alpha\mu)^\T\nu^{-1}z_i-\frac1n\sum(\partial^2_\alpha\mu)^\T\nu^{-1}\delta_i+\frac1n\sum G^{(2)}_i+O_\Pb(\eps).
\]
The second term vanishes uniformly by Lemma~\ref{lem:mart}. The others converge uniformly by Lemma~\ref{lem:lln}. The limit of the $G^{(2)}$ term vanishes at $\beta=\beta_0$ because $G^{(2)}\in\Hspan$, and the third term vanishes at $\alpha_0$. The limit at $\theta_0$ is $-\Gamma_\alpha$.

\emph{Diffusion block.} We have $n^{-1}\partial^2_\beta\ell^{(K)}_n=n^{-1}\sum_i\partial^2_\beta\log q_i+O_\Pb(\eps)$. By Lemma~\ref{lem:lln}, the first term converges uniformly to the $\pi$-integral of the Gaussian expectation of $\partial^2_\beta\log q$. A direct computation shows that this equals $-\Gamma_\beta$ at $\theta_0$.

\emph{Cross block.} Here $(n\eps)^{-1}\partial_\alpha\partial_\beta\ell^{(K)}_n=n^{-1}\sum_i(\partial_\alpha\mu)^\T\partial_\beta(\nu^{-1})z_i(\alpha)+O_\Pb(\eps)$. Since $z_i(\alpha)=z_i-\eps\delta_i(\alpha)$, Lemma~\ref{lem:mart} shows that this tends to zero uniformly in $\theta$.
\end{proof}

The convergence is uniform over all of $\Theta$, not only near $\theta_0$. The near-maximizer argument in the proof of Theorem~\ref{thm:mle} uses this.

\begin{proof}[Proof of Theorem~\ref{thm:an}]
Let $\Sigma_n=D_n^{-1}\partial_\theta\ell^{(K)}_n(\theta_0)$ be the normalized score. With probability tending to one, the following hold:
\begin{itemize}[nosep]
\item we are on $\Omega_n$, and $\tilde s_i=s_i$ for all $i$ (Lemma~\ref{lem:max}), so that $\Sigma_n=D_n^{-1}\sum_i\tilde s_i$;
\item $\hat\theta^{(K)}_n$ lies in the interior of $\Theta$ (Theorem~\ref{thm:cons}, since $\theta_0$ is interior), so $\partial_\theta\ell^{(K)}_n(\hat\theta^{(K)}_n)=0$.
\end{itemize}
By Lemma~\ref{lem:bias}, the normalized compensator is
\[
D_n^{-1}\sum_i\E[\tilde s_i|\cF_{i-1}]=\xi_n\,\frac1n\sum_ib(X_{i-1})\mathbf 1_{\mathcal G_i}+O_\Pb\big(\xi_n\,\eps\,n^{r\delta_0}r_n^r\big)+O(n^{-3/2}).
\]
Since $\delta_0<a/(4r)$, $\eps n^{r\delta_0}r_n^r\to0$. By Lemma~\ref{lem:erg} and Assumption~\ref{A:erg}(ii), $n^{-1}\sum_ib(X_{i-1})\mathbf 1_{\mathcal G_i}\toP\bar B^{(K)}$, because $b$ is $C^1$ with polynomially bounded gradient. Hence $\Sigma_n=S_n+\xi_n\bar B^{(K)}+o_\Pb(1+\xi_n)$. A Taylor expansion of the score gives $0=\Sigma_n+\bar{\mathbb H}_nD_n(\hat\theta^{(K)}_n-\theta_0)$, with $\bar{\mathbb H}_n=\int_0^1\mathbb H_n(\theta_0+u(\hat\theta^{(K)}_n-\theta_0))\,du$. By Lemma~\ref{lem:hess} and consistency, $\bar{\mathbb H}_n\toP-\Gamma$, so $D_n(\hat\theta^{(K)}_n-\theta_0)=\Gamma^{-1}\Sigma_n+o_\Pb(1)\Sigma_n$, which is the stated expansion because $S_n=O_\Pb(1)$ by Lemma~\ref{lem:clt}.

For the parity statement, recall that $H_h(-z;\Sigma)=(-1)^{|h|}H_h(z;\Sigma)$. In \eqref{eq:delta}, $C_{K+1}$ involves only $|h|=2N-K-1$, so it is an odd function of $z$ if $K$ is even and an even function if $K$ is odd. The Euler score $S_\beta$ is even in $z$ and the drift score $(\partial_\alpha\mu)^\T\nu^{-1}z$ is odd. Since $N(0,\nu)$ is symmetric, $b_\beta(x)=0$ for every $x$ when $K$ is even, and $b_\alpha(x)=0$ for every $x$ when $K$ is odd.
\end{proof}

\begin{proof}[Proof of Corollary~\ref{cor:regimes}]
Parts (a) and (b) follow from Theorem~\ref{thm:an}, Lemma~\ref{lem:clt} and Slutsky's lemma, since $\Gamma^{-1}\Gamma\Gamma^{-1}=\mathrm{diag}(\Gamma_\alpha^{-1},\Gamma_\beta^{-1})$. For (c), divide the expansion by $\xi_n\to\infty$: $D_n(\hat\theta^{(K)}_n-\theta_0)/\xi_n=\Gamma^{-1}\bar B^{(K)}+\Gamma^{-1}S_n/\xi_n+o_\Pb(1)\toP\Gamma^{-1}\bar B^{(K)}$.
\end{proof}

\paragraph{Proof of Theorem~\ref{thm:an} for $K=0$.}
For the Euler contrast, $R\equiv0$ and the score is the Euler score itself. In Lemma~\ref{lem:bias}, split $p_\Delta-q=q(\eps C_1+\Delta C_2)+(p_\Delta-p^{(2,\Delta)})$ instead. The term $\eps C_1$ contributes nothing, because $C_1$ involves only Hermite polynomials of degree three, which are orthogonal to both leading scores. The term $\Delta C_2$ contributes $\Delta b^{(1)}$, and the remainder is $O(\Delta^{3/2})$ up to factors $n^{r\delta_0}r_n^r$. Lemmas~\ref{lem:clt} and \ref{lem:hess} hold for $K=0$, and the rest of the proof is unchanged, with $\xi_n=\sqrt n\Delta_n$.

\paragraph{Justification of Remark~\ref{rem:parity}.}
The same computation applies to each block separately. In the block where $b^{(K)}$ vanishes identically, Lemma~\ref{lem:bias} leaves a compensator of order $\sqrt n\Delta_n^{(K+2)/2}$ up to factors $n^{r\delta_0}r_n^r$, which can be made smaller than any power of $n$. For $K=0$, expanding one order further, $p_\Delta-q=q(\eps C_1+\Delta C_2+\Delta^{3/2}C_3)+(p_\Delta-p^{(3,\Delta)})$; by parity, the score bias is $\Delta\,b^{(1)}_\beta$ in the diffusion block and $\Delta^{3/2}b^{(2)}_\alpha$ in the drift block. These block-wise statements about the score bias are proved.

\paragraph{Heuristic: the cross-block effect.}
The following statements concern the estimator rather than the score and are heuristic, since Theorem~\ref{thm:an} controls the non-leading block only up to $o_\Pb(\xi_n)$. We first explain where the cross-block coupling comes from. In the proof of Lemma~\ref{lem:hess}, the Euler part of $(n\eps)^{-1}\partial_\alpha\partial_\beta\ell^{(K)}_n(\theta_0)$ carries the prefactor $n^{-1}=\eps(n\eps)^{-1}$ and is only $O_\Pb(n^{-1/2}+\eps\Delta)$. For $K\ge1$, the coupling of order $\eps$ comes from the correction: by Lemma~\ref{lem:struct}(d), $(n\eps)^{-1}\sum_i\partial_\alpha\partial_\beta\mathcal L_i(\theta_0)=\eps n^{-1}\sum_i\partial_\beta G^{(1)}(z_i;X_{i-1},\theta_0)+O_\Pb(\eps^2)$. Although $G^{(1)}(\cdot;x,\alpha,\beta)$ has mean zero under $N(0,\nu(x;\beta))$ for every $\beta$, differentiating this identity in $\beta$ gives $\E[\partial_\beta G^{(1)}(Z;x,\theta_0)]=-\mathrm{Cov}\big(G^{(1)}(Z;x,\theta_0),S_\beta(Z;x)\big)$, which is in general nonzero when $C_2$ depends on $\alpha$ ($K\ge2$, including the Ornstein--Uhlenbeck model) or when $C_1\neq0$. This is the term that transmits the leading bias of one block to the other at order $\sqrt n\Delta_n^{(K+2)/2}$. For the Euler contrast, by contrast, the normalized cross-block Hessian at $\theta_0$ is only $O_\Pb(n^{-1/2}+\Delta)$, because $R\equiv0$, so the drift bias of the Euler estimator should be the pure score bias $\Delta^{3/2}b^{(2)}_\alpha$. This explains why the Euler and $K=2$ estimators have the same leading drift bias in Section~\ref{sec:mc}. For $K=1$ with $C_1\neq0$, the score bias and the cross-block effect combine at order $\Delta^{3/2}$, which is why that entry of the table in Remark~\ref{rem:parity} is an upper bound only.

\section{Proofs for Sections~\ref{sec:exact} and \ref{sec:other}}\label{app:mle}

All proofs in this appendix compare the exact log-likelihood, or an approximate one, with the explicit contrast $\ell^{(K)}_n$ of high order, and then use the results of Appendices~\ref{app:cons} and \ref{app:an} for that contrast.

\paragraph{Proof of Theorem~\ref{thm:mle}(a).}
Fix $K$ as in (a), and choose $\delta>0$ so small that $a(K+1)/2>1+(1+\delta)m\Lambda+r_K\delta$. Work on the intersection of $\Omega_n$, $\{\max_i|X_i|\le n^\delta\}$ and the event of Lemma~\ref{lem:max}(c), whose probability tends to one. There, $p^{(K,\Delta)}\ge\frac34q$, and $q\ge(2\pi\lambda_{\max}\Delta)^{-m/2}\exp(-|z|^2/(2\lambda_{\min}))$. By Proposition~\ref{prop:unif}, at every observation and for every $\theta$,
\[
\begin{aligned}\Big|\frac{p_\Delta}{p^{(K,\Delta)}}-1\Big|&\le\frac{C\jb{X_{i-1}}^{r_K}\Delta^{(K+1-m)/2}}{q}\le Cn^{r_K\delta}\Delta^{(K+1)/2}\exp\Big(\frac{|z_i(\alpha)|^2}{2\lambda_{\min}}\Big)\\&\le Cn^{r_K\delta}\Delta^{(K+1)/2}n^{(1+\delta)m\Lambda}.\end{aligned}
\]
With $\Delta_n\le Cn^{-a}$, the right-hand side is $o(n^{-1})$. Since $|\log(1+u)|\le2|u|$ for $|u|\le\frac12$, summing over $i$ gives $\eta_n:=\sup_{\theta}|\ell_n(\theta)-\ell^{(K)}_n(\theta)|=o_\Pb(1)$. \qed

\paragraph{Proof of Theorem~\ref{thm:mle}(b) and (c).}
The assumptions of Theorem~\ref{thm:an} hold for this $K$, and $\xi_n^2=n\Delta_n^{K+1}\le Cn^{1-a(K+1)}\to0$. Write $\ell_n(\hat\theta_n)\ge\sup_\Theta\ell_n-\rho_n$ with $\rho_n=o_\Pb(1)$; for an exact maximizer, $\rho_n=0$. Then
\[
\ell^{(K)}_n(\hat\theta_n)\ge\ell_n(\hat\theta_n)-\eta_n\ge\ell_n(\hat\theta^{(K)}_n)-\eta_n-\rho_n\ge\sup\ell^{(K)}_n-2\eta_n-\rho_n.
\]
By Remark~\ref{rem:near}, $\hat\theta_n\toP\theta_0$. Let $u_n=D_n(\hat\theta_n-\theta_0)$, let $\Sigma_n$ be the normalized score of $\ell^{(K)}_n$ at $\theta_0$, and let $f_n(u)=u^\T\Sigma_n-\frac12u^\T\Gamma u$. By the proof of Theorem~\ref{thm:an}, $\Sigma_n=S_n+o_\Pb(1)$. A second-order Taylor expansion and Lemma~\ref{lem:hess} give, uniformly over $\{u:|D_n^{-1}u|\le\eta'_n\}$ for any $\eta'_n\to0$,
\begin{equation}\label{eq:laq}
\ell^{(K)}_n(\theta_0+D_n^{-1}u)-\ell^{(K)}_n(\theta_0)=u^\T\Sigma_n-\tfrac12u^\T(\Gamma+o_\Pb(1))u .
\end{equation}
Choose $\eta'_n\to0$ with $\Pb(|\hat\theta_n-\theta_0|>\eta_n')\to0$. Evaluating \eqref{eq:laq} at $u_n$ and using $\ell^{(K)}_n(\hat\theta_n)\ge\ell_n^{(K)}(\theta_0)-2\eta_n-\rho_n$ gives $\frac12u_n^\T(\Gamma+o_\Pb(1))u_n\le|u_n||\Sigma_n|+2\eta_n+\rho_n$. Hence $|u_n|=O_\Pb(1)$. On compact sets of $u$, \eqref{eq:laq} shows that $u_n$ is an $o_\Pb(1)$-maximizer of $f_n$. The point $\Gamma^{-1}\Sigma_n$ is $O_\Pb(1)$, so it also lies in the region where \eqref{eq:laq} holds. Strong concavity gives $|u-\Gamma^{-1}\Sigma_n|^2\le2\lambda_{\min}(\Gamma)^{-1}(\max f_n-f_n(u))$, so $u_n=\Gamma^{-1}\Sigma_n+o_\Pb(1)=\Gamma^{-1}S_n+o_\Pb(1)$. Lemma~\ref{lem:clt} completes the proof of (b). For (c), let $K'\ge1$ satisfy $n\Delta_n^{K'+1}\to0$. By Theorem~\ref{thm:an}, $D_n(\hat\theta^{(K')}_n-\theta_0)=\Gamma^{-1}S^{(K')}_n+o_\Pb(1)$, where $S^{(K')}_n$ is the martingale part built from order $K'$. By Lemma~\ref{lem:clt}, $S^{(K')}_n$ and $S_n$ both differ from the Euler martingale part by $o_\Pb(1)$. Subtracting the two expansions gives (c). \qed

\paragraph{Proof of Corollary~\ref{cor:lan}.}
Take $K$ as in Theorem~\ref{thm:mle}(a) and let $u_n$ be bounded. Then $|D_n^{-1}u_n|\to0$, so by part (a) and \eqref{eq:laq},
$\ell_n(\theta_0+D_n^{-1}u_n)-\ell_n(\theta_0)=\ell^{(K)}_n(\theta_0+D_n^{-1}u_n)-\ell^{(K)}_n(\theta_0)+o_\Pb(1)=u_n^\T S_n-\frac12u_n^\T\Gamma u_n+o_\Pb(1)$, and $S_n\toD N(0,\Gamma)$ by Lemma~\ref{lem:clt}. This is LAN at $\theta_0$ with rate matrix $D_n$ and information $\Gamma$. By Theorem~\ref{thm:mle}(b,c), the exact MLE and every delta-expansion estimator with $n\Delta_n^{K'+1}\to0$ are asymptotically linear in the central sequence $S_n$, with influence matrix $\Gamma^{-1}$. In a LAN experiment, estimators of this form are regular by Le Cam's third lemma, and they attain the bounds of the convolution and local asymptotic minimax theorems. These results hold for LAN experiments with matrix norming such as $D_n$ \citep[Ch.~6]{LeCamYang2000}; see also \citet{Jeganathan1982}. \qed

\paragraph{Proof of Theorem~\ref{thm:general}(b).}
We show that $\sup_\Theta|\tilde\ell_n-\ell_n|\toP0$ under either condition. Fix $\delta>0$ small, to be chosen below, and work on the intersection $\mathcal E'_n$ of the events $\{\max_i|X_i|\le n^\delta\}$, $\{\max_i\sup_\alpha|z_i(\alpha)|^2\le2(1+\delta)m\lambda_{\max}\log n\}$ and $\Omega_n$. Its probability tends to one by Lemma~\ref{lem:max}. On $\mathcal E'_n$ we also have $\max_i|X_i-X_{i-1}|/\eps\le C\sqrt{\log n}$, because $X_i-X_{i-1}=\eps z_i(\alpha_0)+\Delta\mu(X_{i-1};\alpha_0)$.

\emph{Condition (i).} Take $K$ with $K+1>2(1+m\Lambda)/a$, and write $p^{(K)}=p^{(K,\Delta)}$. As in the proof of Theorem~\ref{thm:mle}(a), on $\mathcal E'_n$ we have $p^{(K)}\ge\frac34q$ and $1/q\le C\Delta^{m/2}n^{(1+\delta)m\Lambda}$ at every observation and for every $\theta$. Hence, by condition~(i) and Proposition~\ref{prop:unif},
\[
\Big|\frac{\tilde p_n}{p^{(K)}}-1\Big|+\Big|\frac{p_\Delta}{p^{(K)}}-1\Big|\le Cn^{r\delta}\big(\Delta^{\gamma}+\Delta^{(K+1)/2}\big)n^{(1+\delta)m\Lambda}
\]
at every observation, uniformly in $\theta$. Since $\Delta_n\le Cn^{-a}$ and both $a\gamma$ and $a(K+1)/2$ exceed $1+m\Lambda$, the right-hand side is $o(n^{-1})$ for $\delta$ small. In particular $\tilde p_n>0$ at all observations, and since $|\log(1+u)|\le2|u|$ for $|u|\le\frac12$, $|\log\tilde p_n-\log p_\Delta|=o(n^{-1})$ at every observation, uniformly in $\theta$. Summing over $i$ gives $\sup_\Theta|\tilde\ell_n-\ell_n|\toP0$.

\emph{Condition (ii).} On $\mathcal E'_n$, the right-hand side of condition~(ii) is at most $Cn^{r\delta}\Delta^\gamma(\log n)^{\varrho/2}$ at every observation, uniformly in $\theta$. Summing over $i$ gives at most $Cn^{1+r\delta-a\gamma}(\log n)^{\varrho/2}\to0$ for $\delta$ small, because $a\gamma>1$. \qed

\paragraph{Proof of Theorem~\ref{thm:general}(a).}
Take $K$ as in Theorem~\ref{thm:mle}(a), so that $\eta_n=\sup_\Theta|\ell_n-\ell^{(K)}_n|=o_\Pb(1)$, and let $\tilde\rho_n=o_\Pb(1)$ be the optimization slack of $\tilde\theta_n$. Write $g_n=\sup_\Theta|\tilde\ell_n-\ell_n|=o_\Pb(n\Delta_n)$ and $h_n=\sup_{|\theta-\theta_0|\le\eta}|\tilde\ell_n-\ell_n|=o_\Pb(1)$.

\emph{Consistency.} For every $\theta$, $\ell^{(K)}_n(\tilde\theta_n)\ge\tilde\ell_n(\tilde\theta_n)-g_n-\eta_n\ge\tilde\ell_n(\theta)-g_n-\eta_n-\tilde\rho_n\ge\ell^{(K)}_n(\theta)-2(g_n+\eta_n)-\tilde\rho_n$. Hence $\ell^{(K)}_n(\tilde\theta_n)\ge\sup\ell^{(K)}_n-o_\Pb(n\Delta_n)$, and Remark~\ref{rem:near} gives $\tilde\theta_n\toP\theta_0$.

\emph{Rate and limit.} With probability tending to one, $\tilde\theta_n$, $\theta_0$ and $\theta_0+D_n^{-1}\Gamma^{-1}\Sigma_n$ all lie in $\{|\theta-\theta_0|\le\eta\}$. For $\theta$ in this set, the same chain of inequalities with $h_n$ in place of $g_n$ gives $\ell^{(K)}_n(\tilde\theta_n)\ge\ell^{(K)}_n(\theta)-2(h_n+\eta_n)-\tilde\rho_n$. The proof of Theorem~\ref{thm:mle}(b) compares $\hat\theta_n$ only with $\theta_0$ and with $\theta_0+D_n^{-1}\Gamma^{-1}\Sigma_n$. It therefore applies verbatim to $\tilde\theta_n$, with $o_\Pb(1)$ slack, and gives $D_n(\tilde\theta_n-\theta_0)=\Gamma^{-1}S_n+o_\Pb(1)$. Subtracting the expansion of the exact MLE gives $D_n(\tilde\theta_n-\hat\theta_n)\toP0$, and efficiency follows as in Corollary~\ref{cor:lan}. \qed

\paragraph{Proof of Corollary~\ref{cor:others}.}
(a) Proposition~\ref{prop:unif} holds for every centering with $|\mu_0(x;\theta)|\le C\jb{x}$, and it is condition~(i) with $\gamma=(K+1)/2$.

(b) For each $\beta$, $Y=F(X;\beta)$ is a diffusion with unit volatility and drift $\mu_Y(y;\theta)=(\mu/\sigma-\frac12\partial_x\sigma)(F^{-1}(y;\beta);\theta)$ under $\Pb_\theta$, and the change-of-variables formula gives $p_\Delta(x'|x;\theta)=p_{Y,\Delta}(F(x';\beta)|F(x;\beta);\theta)/\sigma(x';\beta)$ exactly. Hence $\tilde p_n-p_\Delta=(p^{\mathrm{LSL}}_Y-p_{Y,\Delta})/\sigma(x';\beta)$, evaluated at $y=F(x;\beta)$ and $y'=F(x';\beta)$. By assumption, $p^{\mathrm{LSL}}_Y$ is the delta expansion of $Y$ of order $2K'$ with $\mu_0=0$, so Proposition~\ref{prop:unif}, applied to $Y$ with $m=1$, bounds this difference by $C\jb{y}^r\Delta^{K'}/\lambda_{\min}^{1/2}$. Since $F$ is bi-Lipschitz, $\jb{y}\le C\jb{x}$. This is condition~(i) with $\gamma=K'+\frac12=(2K'+1)/2$ and $m=1$, and $\gamma>(1+\Lambda)/a$ is the stated condition on $K'$. \qed

\section{Proof of Theorem~\ref{thm:moment}}\label{app:moment}

Write $U_n=D_n(\tilde\theta_n-\theta_0)$ for $\tilde\theta_n\in\{\hat\theta_n,\hat\theta^{(K)}_n\}$. Since $U_n\toD N$ by Theorem~\ref{thm:mle}(b) and Corollary~\ref{cor:regimes}(a), it suffices to show that $\sup_n\E|U_n|^p<\infty$ for every $p$, because uniform integrability of $f(U_n)$ follows. Because $\Theta$ is bounded, $|U_n|\le C\sqrt n$ always. It is therefore enough to find, for every $L>0$, events $\mathcal E_n$ with $\Pb(\mathcal E_n^c)=O(n^{-L})$ on which $|U_n|\le C(1+|\Sigma_n|)$, where $\Sigma_n$ is a normalized score with $\sup_n\E|\Sigma_n|^p<\infty$. Then $\E|U_n|^p\le C\E(1+|\Sigma_n|)^p+Cn^{p/2-L}$, which is bounded for $L\ge p/2$. We say that an event holds \emph{with overwhelming probability} (w.o.p.) if its complement has probability $O(n^{-L})$ for every $L$.

\paragraph{Step 1: polynomial deviation bounds.}
We show that each convergence in probability used in Appendices~\ref{app:prelim}--\ref{app:an} has a version with overwhelming probability. Three tools suffice.
\begin{enumerate}[label=(T\arabic*),leftmargin=*,nosep]
\item \emph{Martingales.} By Burkholder's inequality and the Sobolev argument in the proof of Lemma~\ref{lem:mart}, normalized martingale sums with weighted increments satisfy
\[\E\sup_\theta|n^{-1/2}M_n(\theta)|^p\le C_p\qquad\text{for every }p.\]
\item \emph{Ergodic averages.} For $g(x;\vartheta)$ with $|g|+|\nabla_xg|+|\partial_\vartheta g|+|\nabla_x\partial_\vartheta g|\le C\jb{x}^r$, Assumption~\ref{A:rate} applied to $g(\cdot;\vartheta)$ and $\partial_\vartheta g(\cdot;\vartheta)$, the discretization bound in the proof of Lemma~\ref{lem:erg} taken in $L^p$, and the Sobolev embedding in $\vartheta$ give $\E\sup_\vartheta|n^{-1}\sum_ig(X_{i-1};\vartheta)-\int g\,d\pi|^p\le C_p((n\Delta_n)^{-p/2}+\Delta_n^{p/2})$.
\item \emph{Remainders.} Terms bounded by $\eps^c$ times an average of $\jb{X_{i-1}}^r(1+\sup_\alpha|z_i(\alpha)|^r)$ have $L^p$ norms $O(\eps^c)$ by Lemma~\ref{lem:mom}(d) and Assumption~\ref{A:erg}(ii). Terms bounded by powers of $\eps n^{r\delta}\log n$ are controlled on the events of Lemma~\ref{lem:max}(a,c), which hold with overwhelming probability when $c=c(L)$ is large in part (c).
\end{enumerate}
Because $\Delta_n\le Cn^{-a}$ and $n\Delta_n\ge cn^{1-a_1}$, every bound in (T1)--(T3) is polynomially small, and Markov's inequality with $p$ large turns it into a statement with overwhelming probability. Lemma by lemma:
\begin{itemize}[nosep]
\item Lemma~\ref{lem:max}(d), and hence $\Omega_n$: (T3).
\item Lemma~\ref{lem:erg}: (T2).
\item Lemma~\ref{lem:lln}: the replacement of $z_i(\alpha')$ by $z_i$ and the error in the conditional mean are (T3); the martingale part is (T1); the average of $\bar f(X_{i-1};\vartheta)$ is (T2), applied to $\bar f$ and $\partial_\vartheta\bar f$, which have polynomially bounded $x$-gradients. In our applications $f$ is smooth in $\vartheta$, so the Sobolev argument replaces the $\eta$-net.
\item Lemma~\ref{lem:mart}: (T1), plus Lemma~\ref{lem:mom}(b) for the drift term.
\item Lemma~\ref{lem:hess}: a combination of Lemmas~\ref{lem:lln} and \ref{lem:mart} and $O(\eps)$ remainders of type (T3); the limit $\mathbb H$ is deterministic and uniformly continuous, so no further probability bound is needed.
\end{itemize}
Consequently, every statement ``$\sup_\theta|A_n(\theta)|\toP0$'' used in the proofs of Theorems~\ref{thm:cons}--\ref{thm:mle} holds in the form ``$\sup_\theta|A_n(\theta)|\le\eta$ with overwhelming probability'' for every fixed $\eta>0$. The lower bound on $\Delta_n$ enters only here, but Steps~2--4 rely on it through this step.

\paragraph{Step 2: consistency with overwhelming probability.}
The proof of Theorem~\ref{thm:cons} uses only uniform laws of large numbers and $\Omega_n$, so for every $\eta>0$, $|\hat\theta^{(K)}_n-\theta_0|\le\eta$ w.o.p. For the exact MLE, fix $L$, take $c=c(L)$ in Lemma~\ref{lem:max}(c), and choose $K_L$ so large that the proof of Theorem~\ref{thm:mle}(a), with $c\log n$ in place of $2(1+\delta)m\lambda_{\max}\log n$, gives $\sup_\theta|\ell_n-\ell^{(K_L)}_n|\le1$ outside an event of probability $O(n^{-L})$. Remark~\ref{rem:near} with $\eta_n=2$ then gives $|\hat\theta_n-\theta_0|\le\eta$ outside an event of probability $O(n^{-L})$.

\paragraph{Step 3: the bound on $U_n$.}
Let $K^*=K$ for the delta-expansion estimator and $K^*=K_L$ for the exact MLE, let $\tilde s_i$ be as in Lemma~\ref{lem:bias} for the contrast $\ell^{(K^*)}_n$, with $r_n=\sqrt{c(L)\log n}$, let $\Sigma_n=D_n^{-1}\sum_i\tilde s_i$, which equals the normalized score $D_n^{-1}\partial_\theta\ell^{(K^*)}_n(\theta_0)$ on the event $\mathcal E_n$ defined below, and let $\lambda_\Gamma$ be the smallest eigenvalue of $\Gamma$. By Lemma~\ref{lem:hess}, continuity of $\mathbb H$ and Step~1, there is $\eta>0$ such that $\sup_{|\theta-\theta_0|\le\eta}|\mathbb H_n(\theta)+\Gamma|\le\lambda_\Gamma/4$ w.o.p. Let $\mathcal E_n$ be the intersection of this event with the events of Step~2 and Lemma~\ref{lem:max} used above, and with $\{\tilde s_i=s_i\text{ for all }i\}$. The exponent $\delta_0$ in $\mathcal G_i$ is fixed as in Appendix~\ref{app:an}, independently of $L$; the moment order in Lemma~\ref{lem:max}(a) is then taken large relative to $L/\delta_0$, so that $\Pb(\bigcap_i\mathcal G_i)=1-O(n^{-L})$, while $\eps n^{r\delta_0}r_n^r\to0$ for each fixed $L$. On $\mathcal E_n$:
\begin{itemize}[nosep]
\item for $\tilde\theta_n=\hat\theta^{(K)}_n$, which is then an interior critical point, the Taylor expansion of the proof of Theorem~\ref{thm:an} gives $|U_n|\le(4/(3\lambda_\Gamma))|\Sigma_n|$;
\item for $\tilde\theta_n=\hat\theta_n$, the inequality $\ell^{(K^*)}_n(\hat\theta_n)\ge\ell^{(K^*)}_n(\theta_0)-2$ and the Taylor expansion of $\ell^{(K^*)}_n$ around $\theta_0$ give $U_n^\T\Sigma_n-\frac38\lambda_\Gamma|U_n|^2\ge-2$, hence $|U_n|\le C(1+|\Sigma_n|)$.
\end{itemize}

\paragraph{Step 4: moments of the score.}
Write $\Sigma_n=S_n+D_n^{-1}\sum_i\E[\tilde s_i|\cF_{i-1}]$, which holds on the whole probability space. Burkholder's inequality and the conditional moment bounds in the proof of Lemma~\ref{lem:clt} give $\E|S_n|^p\le C_p$. By Lemma~\ref{lem:bias}, in the form valid on the whole space, the compensator is bounded by $C\xi_nn^{-1}\sum_i\jb{X_{i-1}}^r+C$, whose moments are bounded because $\xi_n\to0$. Hence $\sup_n\E|\Sigma_n|^p<\infty$. \qed

\medskip
The argument does not need a polynomial-type large deviation inequality for the contrast, as in the quasi-likelihood analysis of \citet{Yoshida2011}, because the parameter space is compact and the tail events are controlled directly. The price is the two-sided polynomial condition on $\Delta_n$ and Assumption~\ref{A:rate}.

Two extensions follow from the same argument, and we note them without stating them formally. First, Step~3 for the delta-expansion estimator and Step~4 apply equally to the Euler estimator, $K=0$, under $n\Delta_n^2\to0$, by the case $K=0$ of Theorem~\ref{thm:an}. Second, Steps~2 and 3 for $\hat\theta_n$ use only $\ell_n(\hat\theta_n)\ge\sup_\Theta\ell_n$ to obtain $\ell^{(K_L)}_n(\hat\theta_n)\ge\sup_\Theta\ell^{(K_L)}_n-2$. They therefore also apply to a near-maximizer whose slack $\rho_n$ satisfies $\Pb(\rho_n>C)=O(n^{-L})$ for every $L$, with the constant $2$ replaced by $2+C$, such as a maximizer computed to a fixed numerical tolerance.

\section{Proofs for Section~\ref{sec:lin}}\label{app:lin}

The proof of Proposition~\ref{prop:lin} has three steps. We rescale the process around its starting point, so that its coefficients become bounded near the origin, uniformly in the starting point and the parameter (Lemma~\ref{lem:resc}). We then observe that the delta expansion is equivariant under this rescaling, so that Proposition~\ref{prop:unif} can be applied to a modified process with globally bounded coefficients. Finally, we show that the modification changes the density only by an exponentially small amount (Lemmas~\ref{lem:exit}--\ref{lem:loc}). Throughout this appendix, Assumption~\ref{A:lin} holds. Fix $x$ and $\theta$, and write $\varpi=\jb{x}$ for the local scale. All constants are independent of $x$, $\theta$ and $\Delta$.

\begin{lemma}[rescaled coefficients]\label{lem:resc}
Let $Y_t=(X_t-x)/\varpi$, with $X_0=x$. Then $Y_0=0$ and $dY_t=\bar\mu(Y_t)\,dt+\bar\sigma(Y_t)\,dW_t$, where
\[
\bar\mu(y)=\varpi^{-1}\mu(x+\varpi y;\alpha),\qquad\bar\sigma(y)=\varpi^{-1}\sigma(x+\varpi y;\beta).
\]
For $|y|\le\frac12$ and all multi-indices $a$ (in $y$) and $b$ (in $\theta$), $|\partial^a_y\partial^b_\theta\bar\mu(y)|+|\partial^a_y\partial^b_\theta\bar\sigma(y)|\le C_{a,b}$ and $\frac14\lambda_{\min}I\le\bar\sigma\bar\sigma^\T(y)\le\frac94\lambda_{\max}I$, with constants that do not depend on $x$ or $\theta$.
\end{lemma}

\begin{proof}
For $|y|\le\frac12$, the triangle inequality for $\jb{\cdot}$ gives $\frac12\varpi\le\jb{x+\varpi y}\le\frac32\varpi$. By the chain rule and Assumption~\ref{A:lin}, $|\partial^a_y\partial^b_\theta\bar\mu(y)|=\varpi^{|a|-1}|\partial^a_x\partial^b_\theta\mu(x+\varpi y;\alpha)|\le C_{a,b}\varpi^{|a|-1}\jb{x+\varpi y}^{1-|a|}\le C'_{a,b}$, and the same holds for $\bar\sigma$. The ellipticity bounds follow from $\lambda_{\min}\jb{x+\varpi y}^2\le\nu\le\lambda_{\max}\jb{x+\varpi y}^2$.
\end{proof}

The lemma is where the derivative bounds of Assumption~\ref{A:lin} are used: each derivative of a coefficient loses one power of the local scale, and the rescaling restores it.

\paragraph{The localized process.}
Let $\Pi:\R^m\to\R^m$ be a fixed smooth map with bounded derivatives, $\Pi(y)=y$ for $|y|\le\frac14$ and $|\Pi(y)|\le\frac12$ for all $y$. Let $\tilde Y$ solve $d\tilde Y_t=\tilde\mu(\tilde Y_t)\,dt+\tilde\sigma(\tilde Y_t)\,dW_t$, $\tilde Y_0=0$, with $\tilde\mu=\bar\mu\circ\Pi$ and $\tilde\sigma=\bar\sigma\circ\Pi$. By Lemma~\ref{lem:resc}, $\tilde Y$ satisfies Assumptions~\ref{A:smooth} and \ref{A:ell}, with constants that do not depend on $x$ or $\theta$, and its coefficients agree with those of $Y$ on $\mathcal D=\{|y|<\frac14\}$.

\paragraph{Equivariance of the expansion.}
The delta expansion of $X$ from $x$ with centering $\mu_0$ and the delta expansion of $\tilde Y$ from $0$ with centering $\mu_0/\varpi$ satisfy
\begin{equation}\label{eq:equiv}
p^{(K,\Delta)}_X(x'|x)=\varpi^{-m}\,p^{(K,\Delta)}_{\tilde Y}\big((x'-x)/\varpi\,\big|\,0\big).
\end{equation}
To see this, note that the Gaussian kernels satisfy $F_X=\varpi^{-m}F_{\tilde Y}$ under $y\mapsto x+\varpi y$, and that the generator of $X$ becomes the generator of $Y$. Hence $P_J$ of Appendix~\ref{app:prop} transforms in the same way, for every $J$, and so does its regrouping by powers of $\Delta^{1/2}$. Moreover, $P_J$ at the starting point involves only the derivatives of the coefficients at the starting point (Lemma~\ref{lem:w}), and $\tilde Y$ and $Y$ have the same derivatives at $0$. The same change of variables gives $p_{\Delta,X}(x'|x)=\varpi^{-m}p_{\Delta,Y}((x'-x)/\varpi\,|\,0)$. A by-product of \eqref{eq:equiv} is that the coefficients $w_{N,h}(x)$ of $X$ equal $\varpi^{|h|}$ times those of $\tilde Y$. So, as polynomials in $\tilde z=z/\varpi$, the corrections $C_k$ have coefficients that are bounded uniformly in $x$, together with their $\theta$-derivatives.

It remains to compare the densities of $Y$ and $\tilde Y$ near the origin. They differ only through paths that leave $\mathcal D$, and the next lemma shows that this is exponentially unlikely, uniformly in the starting point.

\begin{lemma}[exit from the localization region]\label{lem:exit}
Let $\tau_{\mathcal D}=\inf\{t\ge0:Y_t\notin\mathcal D\}$. There is $c>0$, depending only on the constants in Assumption~\ref{A:lin}, such that $\Pb(\tau_{\mathcal D}\le t)\le2m\,e^{-c/t}$ for all $t\in(0,1]$, $x$ and $\theta$. The same bound holds for the exit time of $\tilde Y$, which coincides with $\tau_{\mathcal D}$.
\end{lemma}

\begin{proof}
By Lemma~\ref{lem:resc}, $|\bar\mu|\le C_0$ and $\bar\sigma\bar\sigma^\T\le\frac94\lambda_{\max}I$ on $\mathcal D$. Let $t_0=1/(8C_0)$. For $t\le t_0$, the drift moves $Y$ by at most $\frac18$ before $\tau_{\mathcal D}\wedge t$, so on $\{\tau_{\mathcal D}\le t\}$ the stopped martingale $M_u=\int_0^{u\wedge\tau_{\mathcal D}}\bar\sigma(Y_v)\,dW_v$ satisfies $\sup_{u\le t}|M_u|\ge\frac18$, and hence some coordinate exceeds $1/(8\sqrt m)$. Each coordinate has quadratic variation at most $\frac94\lambda_{\max}t$, so the exponential inequality for continuous martingales gives $\Pb(\tau_{\mathcal D}\le t)\le2m\exp\{-1/(288\,m\lambda_{\max}t)\}$. For $t\in(t_0,1]$, the bound holds trivially after decreasing $c$. Since $Y$ and $\tilde Y$ are driven by the same Brownian motion and have the same coefficients on $\mathcal D$, pathwise uniqueness gives $Y_t=\tilde Y_t$ for $t\le\tau_{\mathcal D}$, so the two exit times coincide.
\end{proof}

\begin{lemma}[the excursion measure]\label{lem:excur}
For Borel sets $A\subset\mathcal D$ let $\rho_t(A)=\Pb(Y_t\in A,\ \tau_{\mathcal D}\le t)$ for $t>0$, and $\rho_t=0$ for $t\le0$; define $\tilde\rho_t$ in the same way with $\tilde Y$. Then:
\begin{enumerate}[label=(\alph*),nosep]
\item $\Pb(Y_t\in A)-\Pb(\tilde Y_t\in A)=\rho_t(A)-\tilde\rho_t(A)$ for all Borel $A\subset\mathcal D$ and $t>0$;
\item $\rho$ solves $\partial_t\rho=\mathcal A_Y^*\rho$ in the sense of distributions on $\mathcal D\times\R$, where $\mathcal A_Y$ is the generator of $Y$, and the same holds for $\tilde\rho$;
\item $\rho$ and $\tilde\rho$ have smooth, nonnegative densities on $\mathcal D\times\R$, with $\int_{\mathcal D}\rho(t,y)\,dy\le2me^{-c/t}$ for $t\in(0,1]$ and $\int_{\mathcal D}\rho(t,y)\,dy\le1$ for all $t>0$, and likewise for $\tilde\rho$.
\end{enumerate}
\end{lemma}

\begin{proof}
(a) Since $Y_t=\tilde Y_t$ on $\{t<\tau_{\mathcal D}\}$ by Lemma~\ref{lem:exit}, $\Pb(Y_t\in A,\tau_{\mathcal D}>t)=\Pb(\tilde Y_t\in A,\tau_{\mathcal D}>t)$, and (a) follows. (b) For $\varphi\in C_c^\infty(\mathcal D\times(0,\infty))$, $\varphi(u,Y_u)=0$ whenever $Y_u\notin\mathcal D$, in particular at $u=\tau_{\mathcal D}$. It\^o's formula applied between $t\wedge\tau_{\mathcal D}$ and $t$ therefore gives $\E[\varphi(t,Y_t)\mathbf 1\{\tau_{\mathcal D}\le t\}]=\int_0^t\E[(\partial_u+\mathcal A_Y)\varphi(u,Y_u)\mathbf 1\{\tau_{\mathcal D}\le u\}]\,du$. Taking $t$ beyond the time support of $\varphi$, so that $\varphi(t,\cdot)=0$, yields $\int_0^\infty\rho_u\big((\partial_u+\mathcal A_Y)\varphi(u,\cdot)\big)\,du=0$. For $\varphi\in C_c^\infty(\mathcal D\times\R)$ whose support meets $\{t=0\}$, apply this to $\varphi(u,y)\chi(u/t_0)$, where $\chi$ is smooth, vanishes on $(-\infty,\frac12]$ and equals one on $[1,\infty)$. The extra term is $\int\rho_u\big(\varphi(u,\cdot)\big)t_0^{-1}\chi'(u/t_0)\,du$, which is at most $C\sup_{u\le t_0}\rho_u(\mathcal D)\le2mCe^{-c/t_0}$ by Lemma~\ref{lem:exit}, and the remaining terms converge as $t_0\downarrow0$ by dominated convergence. Since $\rho_u=0$ for $u\le0$, this gives the identity for all $\varphi\in C_c^\infty(\mathcal D\times\R)$. (c) In $\mathcal D$, $\mathcal A_Y^*u=\partial_i\big(\frac12a_{ij}\partial_ju+(\frac12\partial_ja_{ij}-\bar\mu_i)u\big)$, with $a=\bar\sigma\bar\sigma^\T$, is a uniformly parabolic operator in divergence form with smooth, bounded coefficients (Lemma~\ref{lem:resc}). By parabolic hypoellipticity \citep{Hormander1967}, the distributional solution in (b) has a smooth density; it is nonnegative because $\rho_t$ is a nonnegative measure. The first mass bound is Lemma~\ref{lem:exit}, and the second holds because $\rho_t(\mathcal D)\le\Pb(Y_t\in\mathcal D)\le1$.
\end{proof}

\begin{lemma}[localization]\label{lem:loc}
For $|y'|\le\frac18$ and $\Delta\in(0,1]$,
\begin{equation}\label{eq:loc}
\big|p_{\Delta,Y}(y'|0)-p_{\Delta,\tilde Y}(y'|0)\big|\le C\Delta^{-m/2}e^{-c/(3\Delta)},
\end{equation}
with $C$ and $c$ independent of $x$, $\theta$ and $\Delta$.
\end{lemma}

\begin{proof}
By Lemma~\ref{lem:excur}(a) and continuity of all densities involved, $p_{\Delta,Y}-p_{\Delta,\tilde Y}=\rho(\Delta,\cdot)-\tilde\rho(\Delta,\cdot)$ on $\mathcal D$, so it suffices to bound $\rho(\Delta,y')$; the same argument applies to $\tilde\rho$. Let $r=\sqrt\Delta$ and suppose first that $\Delta\le\Delta_0:=1/2000$, so that $4r<\frac18$ and the closed ball $\overline{B_{4r}(y')}$ lies inside $\mathcal D$ for $|y'|\le\frac18$. The parabolic Harnack inequality for nonnegative solutions of divergence-form equations with bounded lower-order terms \citep{Moser1964,Moser1967,AronsonSerrin1967}, applied on the cylinder $B_{4r}(y')\times(t-4r^2,t+4r^2)$, gives
\[
\rho(t,y')\le\sup_{B_r(y')\times(t-r^2,t]}\rho\le C\inf_{B_r(y')\times(t+r^2,t+2r^2)}\rho\le\frac{C}{|B_r|\,r^2}\int_{t+r^2}^{t+2r^2}\!\!\int_{\mathcal D}\rho\le Cr^{-m}\Pb(\tau_{\mathcal D}\le t+2r^2).
\]
The Harnack constant $C$ is uniform because the radius is at most $1$, the lower-order coefficients are bounded, and the ellipticity constants are uniform (Lemma~\ref{lem:resc}). With $t=\Delta=r^2$ and Lemma~\ref{lem:exit}, the right-hand side is at most $C\Delta^{-m/2}e^{-c/(3\Delta)}$. For $\Delta\in(\Delta_0,1]$, the same argument with $r=\sqrt{\Delta_0}$ and the mass bound $\int_{\mathcal D}\rho(t,y)\,dy\le1$, valid for all $t>0$ by Lemma~\ref{lem:excur}(c), shows that $\rho$ and $\tilde\rho$ are bounded by a constant on $\{|y'|\le\frac18\}$, so \eqref{eq:loc} holds after enlarging $C$.
\end{proof}

\paragraph{Proof of Proposition~\ref{prop:lin}.}
Let $y'=(x'-x)/\varpi$, so that $|y'|\le\frac18$. By \eqref{eq:equiv}, the identity $p_{\Delta,X}=\varpi^{-m}p_{\Delta,Y}$, Proposition~\ref{prop:unif} applied to $\tilde Y$ at the starting point $0$, and Lemma~\ref{lem:loc},
\[
\begin{aligned}\varpi^m\big|p^{(K,\Delta)}_X(x'|x)-p_{\Delta,X}(x'|x)\big|&\le\big|p^{(K,\Delta)}_{\tilde Y}(y'|0)-p_{\Delta,\tilde Y}(y'|0)\big|+\big|p_{\Delta,\tilde Y}(y'|0)-p_{\Delta,Y}(y'|0)\big|\\&\le C_K\Delta^{(K+1-m)/2}+C\Delta^{-m/2}e^{-c/(3\Delta)}.\end{aligned}
\]
The constant $C_K$ is uniform in $x$ and $\theta$ because Proposition~\ref{prop:unif} depends on $\tilde Y$ only through the constants in Assumptions~\ref{A:smooth} and \ref{A:ell}, which are uniform by Lemma~\ref{lem:resc}, and the centering $\mu_0/\varpi$ is bounded. The last term is at most $C_K\Delta^{(K+1-m)/2}$ for $\Delta\le1$. Dividing by $\varpi^m$ gives the claim. \qed

\paragraph{Regularity of the exact density.}
Lemma~\ref{lem:reg} remains true under Assumption~\ref{A:lin}. Local boundedness and a local Lipschitz bound for $p_\Delta(\cdot|x;\theta)$ near $x'$, uniform over $(x,\theta)$ in compact sets, follow in two cases. If $|x'-x|\le\jb{x}/8$, they follow from the decomposition $p_{\Delta,Y}=p_{\Delta,\tilde Y}+\rho-\tilde\rho$ above, the bounds of \citet{MenozziPesceZhang2021} for $\tilde Y$, and interior estimates for $\rho$ and $\tilde\rho$. Otherwise $|x'-x|>\jb{x}/8$, so $\jb{x}<8|x'-x|$ and hence $\jb{x'}\le\jb{x}+|x'-x|<9|x'-x|$, that is, $|x'-x|>\jb{x'}/9$. The ball of radius $\jb{x'}/18$ around $x'$ then excludes $x$, and the law of $X_t$, extended by zero for $t\le0$, solves the forward equation on that ball. The Harnack inequality with radius $r\le\min(\jb{x'}/72,\sqrt\Delta/2)$ and total mass at most one then bounds $p_\Delta$ near $x'$, and interior estimates give the Lipschitz bound. Weak continuity of the law of $X_\Delta$ in $(x,\theta)$ holds as before, because the coefficients are Lipschitz in $x$, so the mollification argument in the proof of Lemma~\ref{lem:reg} applies. The only change is that $|\sigma(y;\beta)-\sigma(y;\beta')|\le C|\beta-\beta'|\jb{y}$ now grows with $y$; the Gronwall estimate for $\E|X^{x,\theta}_t-X^{x',\theta'}_t|^2$ still applies, because the starting points are deterministic and range over a compact set, so that $\sup_{t\le\Delta}\E\jb{X^{x',\theta'}_t}^2$ is bounded there, by the standard moment estimate for SDEs whose coefficients grow at most linearly. No moment assumption on the data is involved. Strict positivity follows from the support theorem for diffusions with smooth, nondegenerate coefficients \citep{StroockVaradhan1972}, combined with a Harnack chain.

\paragraph{Proof of Theorem~\ref{thm:lin}.}
The proofs in Appendices~\ref{app:prelim}--\ref{app:moment} use Assumptions~\ref{A:smooth}--\ref{A:ell} through four facts:
\begin{enumerate}[label=(\roman*),nosep]
\item the expansion error bound of Proposition~\ref{prop:unif};
\item polynomial bounds on the coefficients of $C_k$;
\item conditional moment and tail bounds for the standardized increments $z_i$ (Lemmas~\ref{lem:max} and \ref{lem:mom});
\item moments of $X$, used to average polynomial weights.
\end{enumerate}
Under Assumption~\ref{A:lin}, measure increments on the local scale, $\tilde z_i=z_i/\jb{X_{i-1}}$, and replace each weight $\jb{X_{i-1}}^r$ by a constant. Then:
\begin{itemize}[nosep]
\item (i) is replaced by Proposition~\ref{prop:lin}. It is used only at observed transitions, and in Lemma~\ref{lem:bias} on $\{|\tilde z|\le r_n\}$. The observed transitions lie in $\{|x'-x|\le\jb{x}/8\}$ with probability tending to one, because $\max_i|\tilde z_i|=O_\Pb(\sqrt{\log n})$ and $\eps\sqrt{\log n}\to0$. The integration domain $B_n=\{x':|\tilde z(\alpha_0)|\le r_n\}$ in Lemma~\ref{lem:bias} lies there deterministically once $n$ is large: since $x'-x=\eps\jb{x}\tilde z(\alpha_0)+\Delta_n\mu(x;\alpha_0)$ and $|\mu(x;\alpha_0)|\le C\jb{x}$, every $x'\in B_n$ has $|x'-x|\le\jb{x}(\eps r_n+C\Delta_n)\le\jb{x}/8$. The integral over $B_n^c$ in that lemma involves only $\partial_\theta p^{(K,\Delta)}$, not $p_\Delta$. By \eqref{eq:equiv}, and because $\varpi=\jb{x}$ does not depend on $\theta$, the substitution $\tilde z=z/\jb{x}$ removes all factors of $\varpi$, and the Gaussian tail integral over $\{|\tilde z|>r_n\}$ is at most $Cn^{-2}$, uniformly in $X_{i-1}$.
\item (ii) holds with bounded coefficients in $\tilde z$, by the equivariance remark above. Lemma~\ref{lem:struct} holds with $\jb{x}^r$ replaced by $1$ and $z$ by $\tilde z$: the Euler scores $\eps(\partial_\alpha\mu)^\T\nu^{-1}z$ and $\frac12[z^\T\nu^{-1}\partial_{\beta_j}\nu\,\nu^{-1}z-\tr(\nu^{-1}\partial_{\beta_j}\nu)]$ are polynomials in $\tilde z$ with bounded coefficients, since $|\nu^{-1/2}\partial_\alpha\mu|$, $|\nu^{-1}\partial_\beta\nu|$ and $|\nu^{-1/2}(\mu(x;\alpha)-\mu(x;\alpha_0))|$ are bounded.
\item (iii) The conditional moment bounds of Lemma~\ref{lem:mom} need no stopping. Under Assumption~\ref{A:lin} the coefficients are globally Lipschitz with linear growth, so the Burkholder--Davis--Gundy inequality and Gronwall's lemma, applied to the unstopped process on $[t_{i-1},t_i]$, give $\E[\sup_{s\in[t_{i-1},t_i]}|X_s-x|^p\,|\,X_{i-1}=x]\le C_p\Delta^{p/2}\jb{x}^p$. Dividing by $\varpi^p=\jb{x}^p$ yields the bounds of Lemma~\ref{lem:mom} for $\tilde z_i$, uniformly in $X_{i-1}$ and on the whole probability space, including the event that the rescaled process leaves $\mathcal D$. In particular, $\E[|\tilde z_i-\varpi^{-1}\sigma(X_{i-1})\upsilon_i|^p\,|\,\cF_{i-1}]\le C_p\Delta^{p/2}$ and $|\E[\tilde z_i|\cF_{i-1}]|\le C\Delta^{3/2}$, with no weight. For the latter, $|\mathcal A\mu(y)|\le C\jb{y}$, and $\E[\jb{X_u}\,|\,X_0=x]\le C\jb{x}$ for $u\le1$. In Lemma~\ref{lem:max}(b,c), stop the rescaled process at its exit from $\{|y|<\eta\}$, with $\eta=\eta(\delta)$ small. Before that time, the quadratic variation of the rescaled martingale is at most $(1+\eta)^2\lambda_{\max}\Delta$, and the factor $(1+\eta)^2$ is absorbed into $\delta$, which preserves the constant $2(1+\delta)m\lambda_{\max}$. The exit probability is at most $2me^{-c\eta^2/\Delta}$, which is super-polynomially small. Lemma~\ref{lem:max}(a) is no longer needed, and the term $O(n^{-L})$ in Lemma~\ref{lem:max}(c) now comes from this exit probability. Stopping is thus used only for the exponential martingale inequality, which needs bounded quadratic variation; all $L^p$ bounds, including those used to remove truncation indicators by the Cauchy--Schwarz inequality, come from the unstopped estimate.
\item (iv) is no longer needed. Every function averaged over the sample in Lemmas~\ref{lem:erg}--\ref{lem:mart}, \ref{lem:clt} and \ref{lem:hess} becomes, after the substitution $z=\jb{x}\tilde z$, a function $g(x)$ times a polynomial in $\tilde z$, with $|g(x)|+\jb{x}|\nabla g(x)|\le C$. For such $g$, $|g(X_s)-g(X_{i-1})|\le C\min(1,|X_s-X_{i-1}|/\jb{X_{i-1}})$: if $|X_s-X_{i-1}|\le\frac12\jb{X_{i-1}}$, every point $\xi$ on the segment has $\jb{\xi}\ge\frac12\jb{X_{i-1}}$, and the mean value theorem applies; otherwise use $|g|\le C$. The conditional expectation of the right-hand side given $\cF_{i-1}$ is at most $C\eps$, uniformly in $X_{i-1}$, by (iii). Taking expectations afterwards, which requires no moments of the initial law, shows that the discretization error in Lemma~\ref{lem:erg} is $O(\eps)$, and Assumption~\ref{A:erg}(i) gives the ergodic averages for bounded $g$. In Lemma~\ref{lem:lln}, the Gaussian average $\bar f(x;\vartheta)=\int\tilde f(x,\tilde z;\vartheta)\phi(\tilde z;\tilde\nu(x;\beta_0))\,d\tilde z$, with $\tilde\nu=\jb{x}^{-2}\nu$, satisfies the same bound $|\bar f|+\jb{x}|\nabla_x\bar f|\le C$, because $\jb{x}|\nabla_x\tilde\nu(x;\beta_0)|\le C$ by Assumption~\ref{A:lin}.
\end{itemize}
The threshold in Theorem~\ref{thm:mle}(a) is unchanged. Indeed, $q\ge(2\pi\lambda_{\max}\jb{x}^2\Delta)^{-m/2}\exp(-|\tilde z|^2/(2\lambda_{\min}))$, and Proposition~\ref{prop:lin} carries the same factor $(\jb{x}^2\Delta)^{-m/2}$, so the relative error at an observation is at most $C\Delta^{(K+1)/2}n^{(1+\delta)m\Lambda}$. The proof of Theorem~\ref{thm:general}(b) uses (i)--(iii) in the same way, which gives part (b). For part (c), Step~1 of Appendix~\ref{app:moment} uses Assumption~\ref{A:rate} only for the functions $g$ just described, and (T3) needs no moments of $X$. \qed

\section{Expansion coefficients used in Section~\ref{sec:mc}}\label{app:mc}

\paragraph{Ornstein--Uhlenbeck model.}
With $\nu=\mathrm{diag}(s_1^2,s_2^2)$ constant, $C_1\equiv0$. Write $E_1=-\frac12(\kappa\nu+\nu\kappa^\T)$, $E_2=\frac16(\kappa^2\nu+2\kappa\nu\kappa^\T+\nu(\kappa^\T)^2)$, $P_j=\nu^{-1}E_j$ and $a(x)=\frac12\kappa^2x$. Expanding the exact Gaussian transition density at fixed $z$ gives
\[
\begin{aligned}
C_2&=\tfrac12\big[z^\T\nu^{-1}E_1\nu^{-1}z-\tr P_1\big],\qquad C_3=a(x)^\T\nu^{-1}z,\\
C_4&=-\tfrac12\big[\tr P_2-\tfrac12\tr P_1^2\big]-\tfrac12z^\T\big(\nu^{-1}E_1\nu^{-1}E_1\nu^{-1}-\nu^{-1}E_2\nu^{-1}\big)z+\tfrac12C_2^2 .
\end{aligned}
\]
The degree-one term $C_3$ equals $\frac12\sum_iw_{2,e_i}H_{e_i}$ with $w_{2,e_i}(x)=(\mathcal A\mu_i)(x)$. It is responsible for the drift bias of the Euler and $K=2$ estimators. In particular, $\bar B^{(2)}_\alpha=\int(\partial_\alpha\mu)^\T\nu^{-1}a\,d\pi$, and $\bar B^{(1)}_{s_j}=\E[S_{s_j}(Z)C_2(Z)]=-\kappa_{jj}/s_j$, where $S_{s_j}(z)=(z_j^2/s_j^2-1)/s_j$ is the Euler score for $s_j$. These are the constants used in Table~\ref{tab:biaspred}.

\paragraph{General $C_1$ and $C_2$.}
Writing the recursion of \citet{YCW2019} in tensor form, with $H_I=(-1)^{|I|}\phi^{-1}\partial_I\phi$ for $\phi=\phi(\cdot;\nu(x))$ and sums over all ordered indices, we obtain
\[
\begin{aligned}
C_1&=\tfrac14\textstyle\sum\nu_{ab}\,\partial_a\nu_{ij}\,H_{bij},\\
C_2&=\tfrac12\textstyle\sum T_{bi}H_{bi}+\tfrac16\sum\big(W'_{ijkl}H_{ijkl}+W''_{bdij}H_{bdij}\big)+\tfrac1{24}\sum V_{bdijkl}H_{bdijkl},
\end{aligned}
\]
where all coefficients are evaluated at $x$:
\[
\begin{aligned}
T_{bi}&=\textstyle\sum_a\nu_{ab}\partial_a\mu_i+\tfrac12\mathcal A\nu_{bi},\qquad
W'_{ijkl}=\tfrac14\textstyle\sum_{a,b}\nu_{ab}\,\partial_a\nu_{ij}\,\partial_b\nu_{kl},\\
W''_{bdij}&=\tfrac12\textstyle\sum_a\nu_{ab}\,\partial_a\big(\sum_c\nu_{cd}\,\partial_c\nu_{ij}\big),\\
V_{bdijkl}&=\tfrac14\textstyle\sum_{a,c}\nu_{ab}\nu_{cd}\big(2\,\partial_c\nu_{ij}\,\partial_a\nu_{kl}+\partial_a\nu_{ij}\,\partial_c\nu_{kl}\big).
\end{aligned}
\]
After symmetrization, $\frac1{24}\sum V_{bdijkl}H_{bdijkl}=\frac1{32}\sum\nu_{ab}\nu_{cd}\,\partial_a\nu_{ij}\,\partial_c\nu_{kl}\,H_{bdijkl}$, so the degree-six part of $C_2$ equals that of $\frac12C_1^2$. We verified these formulas, and those for the Ornstein--Uhlenbeck model up to $C_4$, against a direct numerical implementation of the recursion (18)--(20) of \citet{YCW2019} in truncated Taylor arithmetic, at randomly drawn $(x,z,\theta)$. They agree to machine precision. Replication code for all computations is provided as supplementary material.

\end{document}